\documentclass[final]{amsart}
\usepackage{amsmath,amssymb,amsthm,graphicx,amscd}
\usepackage{hyperref}
  \newtheorem{theorem}{Theorem}[section]
 \newtheorem{remark}{Remark}[section]

\begin{document}

\title{Numerical Study of Blowup in the Davey-Stewartson 
System}

\author{C.~Klein}
\address{Institut de Math\'ematiques de Bourgogne,
		Universit\'e de Bourgogne, 9 avenue Alain Savary, 21078 Dijon
		Cedex, France}
   \email{christian.klein@u-bourgogne.fr}

\author{B.~Muite}
\address{Department of Mathematics, University of Michigan, 
2074 East Hall, 530 Church Street, MI 48109, USA}
\email{muite@umich.edu}

\author{K.~Roidot}
\address{Institut de Math\'ematiques de Bourgogne,
		Universit\'e de Bourgogne, 9 avenue Alain Savary, 21078 Dijon
		Cedex, France}
    \email{kristelle.roidot@u-bourgogne.fr}

\begin{abstract}
    Nonlinear dispersive partial differential 
       equations such as  the nonlinear 
	   Schr\"odinger equations can have solutions that blow-up. 
	   We numerically study the long time behavior and potential 
	   blowup of  solutions to the 
	   focusing Davey-Stewartson II equation by analyzing 
	   perturbations of  the lump and the Ozawa exact
	   solutions. It is shown in this way that the lump is unstable to both 
	   blowup and dispersion, and that blowup in the  
	   Ozawa solution is generic.
\end{abstract}


\thanks{We thank J.-C.~Saut for helpful discussions and hints. CK 
and KR were supported by the project FroM-PDE funded 
by the European Research Council through the Advanced Investigator 
Grant Scheme, the Conseil R\'egional de Bourgogne via a FABER grant, the Marie-Curie IRSES program RIMMP 
and the ANR via the program ANR-09-BLAN-0117-01. BM is grateful to 
the Institut de Math\'ematiques de Bourgogne for their hospitality 
and financial support by the CNRS where part of this work was completed. 
Part of the computational requirements of this research were 
supported by an allocation of advanced computing resources provided by 
the National Science Foundation. 
This work was granted access to the HPC resources of CCRT/IDRIS under 
the allocation 2011106628 made by GENCI, Kraken at the National 
Institute for Computational Sciences and Hopper at NERSC. 
Computational support was also provided by CRI de Bourgogne and 
SCREMS NSF DMS-1026317.    
}
\maketitle

\section{Introduction}
The Davey-Stewartson (DS) system  models the evolution of weakly nonlinear water waves that
travel predominantly in one direction for which the wave amplitude
is modulated slowly in two horizontal directions \cite{DS}, 
\cite{DR}. It is also used
in plasma physics \cite{NAS,NAS94}, to describe the evolution of a plasma under the
action of a magnetic field. The DS system can be written in the form
\begin{equation}
    \label{DSII}
\begin{array}{ccc}
i
\partial_{t}u+\partial_{xx}u-\alpha\partial_{yy}u+2\rho\left(\Phi+\left|u\right|^{2}\right)u & = & 0,
\\
\partial_{xx}\Phi+\beta\partial_{yy}\Phi+2\partial_{xx}\left|u\right|^{2} & = & 0,
\end{array}
\end{equation}
where $\alpha$, $\beta$ and $\rho$ take the values $\pm1$, and where $\Phi$ is a mean field.
The DS equations can be seen as a two-dimensional nonlinear 
Schr\"odinger (NLS) equation with a nonlocal term if the equation for 
$\Phi$ can be solved for given boundary conditions. They are 
classified \cite{GS} according to the ellipticity or hyperbolicity of the 
operators in the first and second line in (\ref{DSII}). The case $\alpha=\beta$ is 
completely integrable  \cite{AH} and thus provides a $2+1$-dimensional 
generalization of the integrable NLS equation in $1+1$ dimensions 
with a cubic nonlinearity. The integrable cases are
elliptic-hyperbolic  called DS I, and the hyperbolic-elliptic 
called DS II. For both there is a focusing ($\rho=-1$) and a 
defocusing ($\rho=1$) version. The complete integrability of the DS 
equation implies that it has an infinite number of conserved 
quantities, for instance the $L_{2}$ norm  and the energy
\begin{align}\label{energy}
    E[u(t)] & := \frac{1}{2} \int_{\mathbb{T}^2} \bigg[ |\partial_x u(t,x,y)|^2 - 
|\partial_y u(t,x,y)|^2     \nonumber\\
& \left.-\rho\left(|u(t,x,y)|^{4}-\frac{1}{2}\left(\Phi(t,x,y)^{2}+(\partial_{x}^{-1}\partial_{y}\Phi(t,x,y))^{2}\right)\right) 
  \right] d x d y.
    \nonumber
\end{align}

DS reduces to the cubic NLS in one 
dimension if the potential is independent of $y$, and if $\Phi$ satisfies 
certain boundary conditions (for instance rapidly decreasing at 
infinity or periodic).
In the following, we 
will only consider the case DS II ($\alpha=1$) since the mean field 
$\Phi$ is then obtained by inverting an elliptic operator. The 
non-integrable elliptic-elliptic DS is very similar to the $2+1$ dimensional NLS 
equation, see for instance \cite{GS} and \cite{BMS} for numerical 
simulations, and is therefore not studied 
here.

There exist many explicit solutions for the integrable cases which 
thus allow to address the question about the long time behavior of 
solutions for given initial data. For the famous Korteweg-de Vries 
(KdV) equation, 
it is known that general initial data are either radiated away or 
asymptotically decompose 
into solitons. The DS II system and the two-dimensional 
integrable generalization of KdV known as the Kadomtsev-Petviashvili 
I (KP I) equation have so-called lump solutions, a two-dimensional soliton which is 
localized in all spatial directions with an algebraic falloff 
towards infinity. For KP I it was shown 
\cite{AF} that small initial data asymptotically decompose into 
radiation and lumps. It is conjectured that this is also true for 
general initial data. 

For the defocusing DS II global existence in time was shown by Fokas and Sung \cite{FS} 
for solutions of certain classes of Cauchy problems. These initial data will 
simply disperse. The situation is more involved for the focusing 
case. Pelinovski and Sulem \cite{PS2} showed that the lump solution is 
spectrally unstable. In addition the focusing NLS equations in $2+1$ 
dimensions with cubic nonlinearity
have the critical dimension, i.e., solutions from smooth initial data 
can have \emph{blowup}. This means that the solutions lose after 
finite time the regularity of the initial data, a norm of the 
solution or of one of its derivatives becomes infinite.  For focusing 
NLS equations in $2+1$ dimensions, it is known that 
blowup is possible if the energy of the initial data is greater than 
the energy of the ground state solution, see e.g.~\cite{DSnls} and 
references therein,  and 
\cite{MR} for an asymptotic description of the blowup profile. For 
the focusing DS II equation Sung \cite{Sun} gave a smallness condition on the 
Fourier transform $\mathcal{F}[u_{0}]$ of the initial data 
to establish global existence in time for solutions to Cauchy 
problems
\begin{equation}
    ||\mathcal{F}[u_{0}]||_{L_{1}}||\mathcal{F}[u_{0}]||_{L_{\infty}}
    <\frac{\pi^{3}}{2}\left(\frac{\sqrt{5}-1}{2}\right)^{2}\sim 5.92\ldots.
    \label{sungcond}
\end{equation}
 with 
initial data $u_{0}\in L_{p}$, $1\leq p < 2$ with a Fourier 
transform $\mathcal{F}[u_{0}]\in L_{1}\cap L_{\infty}$. 

It is not known whether there is generic blowup for initial data not 
satisfying this condition, nor whether the condition is optimal. 
Since the initial 
data studied in this paper are not in this class, we cannot 
provide further insight into this question.
An explicit solution with blowup for lump-like initial data was given 
by Ozawa \cite{Oza}. It has  an 
$L_{\infty}$ blowup in one point 
($x_{c},y_{c},t_{c}$) and is analytic for all other values of $(x,y,t)$. 
It is unknown whether this is the typical blowup behavior for the 
focusing DS II equation. 

From the point of view of applications, a blowup of a solution does 
not mean that the studied equation is not relevant in this context. It just indicates  the limit of the used approximation. It is 
thus of particular interest, not only in mathematics,
but also in physics, since 
it shows the limits of the applicability of the studied model. This 
breakdown of the model will also in general indicate how to amend the 
used approximations.

In view of the open analytical questions concerning blowup in DS II solutions,  
we study the issue in the present paper numerically, 
which is a highly non-trivial problem for several reasons: first 
DS is a purely dispersive equation which means that the introduction 
of numerical dissipation has to be avoided as much as possible to 
preserve dispersive effects such as rapid oscillations. This
makes the use of spectral methods attractive since they are known for 
minimal numerical dissipation and for their excellent approximation properties for smooth 
functions. But the algebraic falloff of both the lump and the Ozawa 
solution
leads to strong Gibbs phenomena at the boundaries of the 
computational domain if the solutions are periodically continued there. We will nonetheless use Fourier spectral 
methods because they also allow for efficient time integration algorithms 
which should be ideally of high order to avoid a pollution of the 
Fourier coefficients due to numerical errors in the time integration. 

An additional problem is the modulational instability of the focusing 
DS II equation, i.e., a self-induced amplitude modulation of a continuous 
wave propagating in a nonlinear medium, with subsequent generation of 
localized structures, see for instance \cite{Agr,CH,FL} for the NLS 
equation.
Thus to address  numerically questions of stability and blowup of 
its solutions, high resolution is needed which cannot be 
achieved on single processor computers. Therefore we use  
parallel computers to study the related questions. 
The use of Fourier spectral method is also very convenient in this context,
since for a parallel spectral code only existing
optimized serial FFT algorithms are necessary. In addition such codes 
are not memory intensive, in contrast to other approaches 
such as finite difference or finite element methods.
The first numerical 
studies of DS were done by White and Weideman \cite{WW} using 
Fourier spectral methods for the spatial coordinates and a second 
order time splitting scheme. Besse, Mauser and Stimming \cite{BMS} 
used essentially a parallel version of this code to study the Ozawa 
solution and blowup in the focusing elliptic-elliptic DS equation. McConnell, Fokas and 
Pelloni \cite{MFP} used Weideman's code to study 
numerically  DS I and DS II, but did not have enough resolution to 
get conclusive results for the blowup  in perturbations of the lump 
in the focusing DS II case.
In this paper we repeat some of their computations with 
considerably higher resolution.

We use a parallelized version of a fourth order time splitting scheme 
which was studied for DS in \cite{KR}. Obviously it is non-trivial to 
decide numerically whether a solution blows up or whether it just has 
a strong maximum. To allow to make nonetheless reliable statements, 
we perform a series of tests for the numerics. First we test the code on 
known exact solutions with algebraic falloff, the lump and the Ozawa 
solution. We establish that energy conservation can be used to judge 
the quality of the numerics if a sufficient spatial resolution is 
givem. It is shown that the splitting code 
continues to run in general beyond a potential blowup which makes it 
difficult to decide whether there is blowup. 
We argue at examples for the  quintic NLS in $1+1$ dimensions 
(which is known to have blowup solutions) and the Ozawa solution
that energy conservation is a reliable indicator in this case 
since the energy of the solution changes completely after a blowup, whereas it will be in 
accordance with the numerical accuracy after a strong maximum. Thus 
we reproduce well known blowup cases in this way and establish with 
the energy conservation a criterion to ensure the accuracy of the 
numerics also in unknown cases.  Then 
we study perturbations of the lump and the Ozawa solution to see when 
blowup is actually observed.

The paper is organized as follows: in section 2, we describe the 
numerical code and its parallelization, and study as an example the lump 
solution. In section 3 we 
numerically study blowup in the focusing 1+1-dimensional quintic NLS and the Ozawa solution for DS II. In section 4 we discuss 
perturbations of the lump, and in section 5 perturbations of the Ozawa 
solution. In section 6 we give some concluding remarks. 

\section{Numerical methods}
In this paper we are interested in the numerical study of solutions 
to the focusing DS II equation for initial data with algebraic 
falloff towards infinity in all spatial directions. This algebraic 
decrease of the initial data and consequently of the solution for all 
times is in principle not an ideal setting for the use of Fourier 
methods since the periodic continuation of the function at the 
boundaries of the computational domain will lead to Gibbs phenomena. 

Nonetheless there are several reasons for the use of Fourier methods 
in this case: 
First the DS equation is a 
purely dispersive PDE, and we are interested in dispersive effects. 
Thus it important to use numerical methods that introduce as little 
numerical dissipation as possible, and spectral methods are 
especially effective in this context. Furthermore the discrete Fourier transform 
 can be  
efficiently computed with a \emph{fast Fourier transform} (FFT). In addition  Fourier 
methods allow to use  splitting techniques for the time integration 
as explained below in a very efficient way. Last but not least the 
focusing DS II equation is known to have a modulation instability 
which makes the use of high resolution necessary, especially close to 
the blowup situations we want to study. This instability leads to an 
artificial increase of the high wave numbers which eventually breaks 
the code, if not enough spatial resolution is provided (see for 
instance \cite{ckkdvnls} for the focusing NLS equation). It is not 
possible to reach the necessary resolution on single processors which 
makes a parallelization of the code obligatory. As explained below, 
this can be conveniently done for 2-dimensional (even for 3-dimensional) Fourier transformations where 
the task of the 1-dimensional FFTs is performed simultaneously by several processors. 
This reduces also the memory requirements per processor 
over alternative approaches.

\subsection{Splitting Methods}

Splitting methods are very efficient if an equation can be split into two or
more equations which can be directly integrated. They are 
unconditionally stable.
The  motivation for these methods is the Trotter-Kato formula
\cite{TK,Ka}
\begin{equation}\label{e11}
 {\lim}_{n\rightarrow\infty}\left(e^{-tA/n}e^{-tB/n}\right)^{n}=e^{-t\left(A+B\right)}
\end{equation}
where $A$ and $B$ are certain unbounded linear operators, for details 
see \cite{Ka}.
In particular this includes the cases studied by Bagrinovskii and
Godunov in \cite{BG} and by Strang \cite{ST}. For hyperbolic equations,
first references are Tappert \cite{Tap} and Hardin and Tappert \cite{HT}
who introduced the split step method for the NLS
equation.

The idea of these methods for an equation of the form $u_{t}=\left(A+B\right)u$ is to write the solution in the
form
\[
u(t)=\exp(c_{1}tA)\exp(d_{1}tB)\exp(c_{2}tA)\exp(d_{2}tB)\cdots\exp(c_{k}tA)\exp(d_{k}tB)u(0)
\]
where $(c_{1},\,\ldots,\, c_{k})$ and $(d_{1},\,\ldots,\, d_{k})$
are sets of real numbers that represent fractional time steps.  
Yoshida \cite{Y} gave an approach 
which produces split step methods of any even order. 
The DS equation can be split into
\begin{eqnarray}
i\partial_{t} u=(-\partial_{xx}u+\alpha \partial_{yy}u),\,\,\, 
\partial_{xx}\Phi+\alpha\partial_{yy}\Phi+2\partial_{xx}\left(\left|u\right|^{2}\right)=0 ,
\\
i\partial_{t} u= -2\rho\left(\Phi+\left|u\right|^{2}\right)u,  
\label{2.4}              
\end{eqnarray}
which are explicitly integrable, the first two in Fourier space, 
equation (\ref{2.4}) in physical space since 
$|u|^{2}$ and thus $\Phi$ is constant in time for this equation. Convergence of second order 
splitting along these lines was studied in \cite{BMS}. We use here 
fourth order splitting as given in \cite{Y} and already studied in \cite{KR} for the DS II equation.
In the latter reference, it was shown that this scheme is very 
efficient in this context.  The method is convenient for
parallel computing, because of  easy coding (loops) and low memory 
requirements. 

Notice that the splitting method in the form (\ref{2.4}) conserves 
the $L_{2}$ norm: the first equation implies that 
its  solution in Fourier space is just the initial condition (from the 
last time step) multiplied by a factor $e^{i\phi}$ with 
$\phi\in\mathbb{R}$. Thus the $L_{2}$ norm is constant for solutions 
to this equation because of Parseval's theorem. The second equation 
as mentioned conserves the $L_{2}$ norm exactly. Thus the used 
splitting scheme has the conservation of the $L_{2}$ norm 
implemented. As we will show in the following, this does not 
guarantee the accuracy of the numerical solution since other 
conserved quantities as the energy the conservation of which is not implemented 
might not be numerically conserved. In fact we will show that 
the numerically computed energy provides a valid 
indicator of the quality of the numerics.

\subsection{Parallelization of the code}

Since high resolution is needed to numerically examine the focusing
DS II equation, the code is parallelized to reduce the wall clock
time required to run the simulation and to allow the problem to fit
in memory. The runs typically used $N_x=N_y=2^{15}$, where
$N_x$ and $N_y$ denote the number of Fourier modes in $x$ and $y$
respectively. The parallelization is done by a slab domain
decomposition. The grid points are given by
$$x_n=\frac{2\pi n L_x}{N_x},\quad y_m=\frac{2\pi m L_y}{N_y},$$
so that the numerical solution is in the computational domain
$$x\times y \in [-L_x\pi,L_x\pi]\times[-L_y\pi,L_y\pi].$$
In the computations, $L_x=L_y$ is chosen large enough so that the
numerical solution is small at the boundaries, and hence a numerical
solution on a periodic domain can be considered as a good
approximation to the solution on an unbounded domain. The approximate
solution $u$ is represented by an $N_x\times N_y$ matrix, which is
distributed among the MPI processes (note that each MPI process uses
a single core). For programming ease and for the efficiency of the
Fourier transform, $N_x$ and $N_y$ are chosen to be powers of two.
The number of MPI processes, $n_p$ is chosen to divide $N_x$ and
$N_y$ perfectly, so that each process holds $N_x\times N_y/n_p$
elements of $u$, for example process $i$ holds the elements
$$u(1:N_x,(i-1)\frac{N_y}{n_p}+1:i\frac{N_y}{n_p})$$
in the global array. To avoid performing global Fourier transforms
which are inefficient, the array is transposed once all the one
dimensional Fourier transforms in the $x$ direction have been done.
Since the data is evenly distributed among the MPI processes, this
transpose is efficiently implemented using {\tt MPI\_ALLTOALL} 
\cite{MPI}. After the transpose, the Fourier
transform $\hat{u}$ is distributed on the processes so that process
$i$ holds the elements corresponding to 
$$\hat{u}((i-1)\frac{N_x}{n_p}+1:\frac{iN_x}{n_p},1:Ny),$$
on which the second set of one dimensional FFTs can be done. The one
dimensional FFTs were done using FFTW 3.0, FFTW 3.1 and FFTW 3.2\footnote{
\url{http://www.fftw.org}} which are close to optimal on x86
architectures and allow the resulting program to be portable but
still simple.

\subsection{Lump solution of the focusing DS II equation}
To test the performance of the code, we first propagate initial data 
from known exact solutions and compare the numerical and the exact 
solution at a later time.
 
The focusing DS II equation has solitonic solutions which are regular 
for all $x,y,t$, and which are
localized with an algebraic falloff towards infinity, known as lumps 
\cite{APP}. 
The single lump  is given by
\begin{equation}\label{lump}
 u(x,y,t) = 2c \frac{\exp \left( -2i(\xi x - \eta y + 
2(\xi^{2}-\eta^{2} )t)\right)}{|x + 4\xi t + i(y + 4\eta t) +
 z_{0}|^2+|c|^2}
\end{equation}
 where $(c,z_{0})\in \mathbb{C}^2$ and $(\xi,\eta)\in\mathbb{R}^2$ 
 are constants. The lump moves with constant 
 velocity $(-2\xi, -2\eta)$ and decays as $(x^2+y^2)^{-1}$ for 
 $x,y\to\infty$.

\begin{figure}[htb!]
\centering
\includegraphics[width=0.45\textwidth]{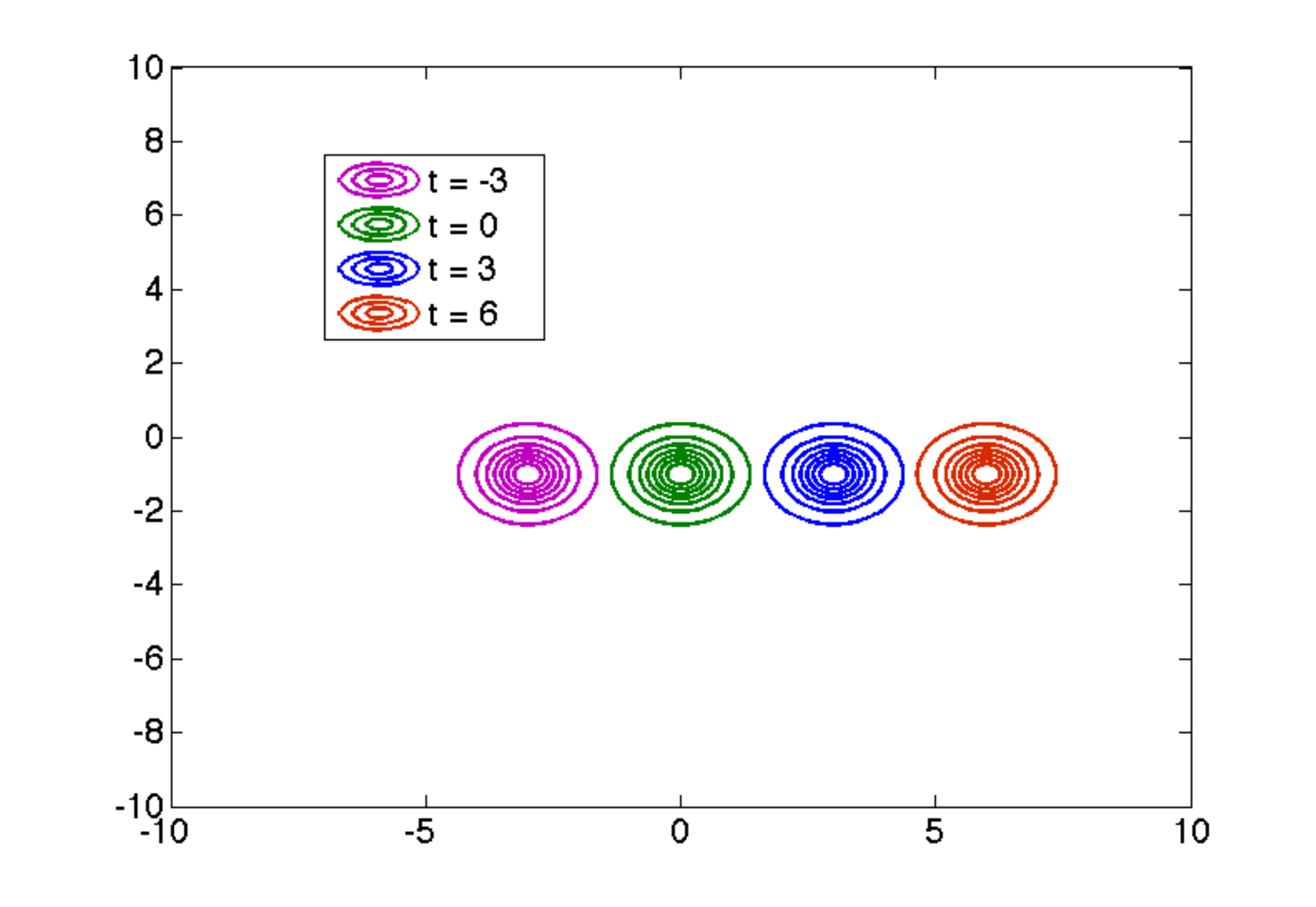}
\includegraphics[width=0.45\textwidth]{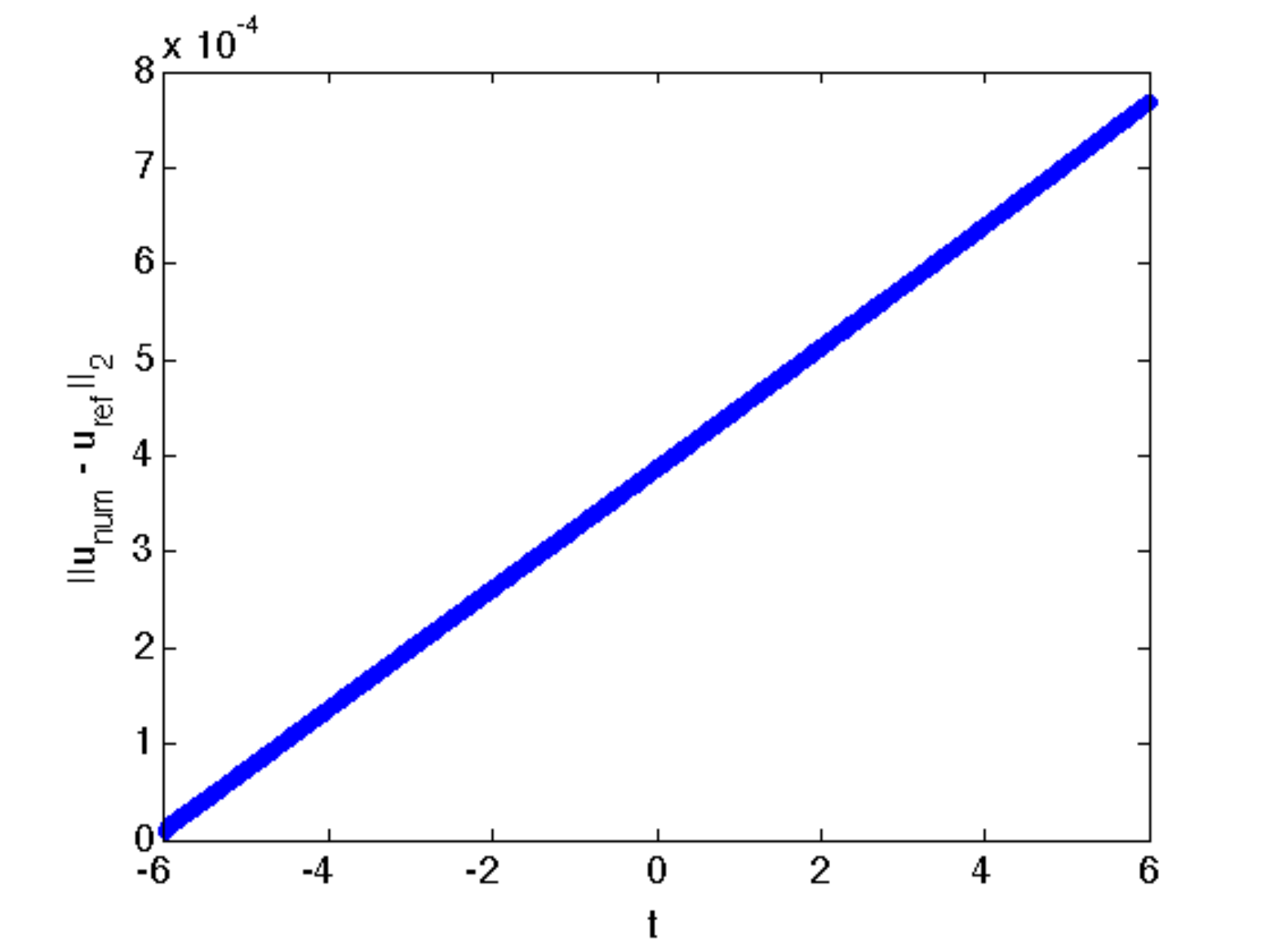}
\caption{Contours of $|u|^2$ on the left and a plot of $||u_{exact} - 
u_{num}||_2$ on the right in 
dependence of time for the solution to the focusing DS II equation 
(\ref{DSII}) for lump initial data (\ref{lump}).} 
\label{ratcont}
\end{figure}

We choose  $N_{x}=N_{y}=2^{14}$ and $L_{x}=L_{y}=50$,
with $\xi=0, \eta=-1, z_{0}=1$ and $c=1$. The 
large values for $L_{x}$ and $L_{y}$ are necessary to ensure that the 
solution is small at the boundaries of the computational domain to 
reduce Gibbs phenomena. The difference for the mass of the lump and 
the computed mass on this periodic setting is of the order of 
$6*10^{-5}$. The initial data for $t=-6$ are propagated 
with $N_{t}=1000$ time steps until $t=6$.  
In  Fig.~\ref{ratcont} contours of the solution
 at different times are shown. Here and in the following we always 
 show closeups of the solution. The actual computation is done on the 
 stated much larger domain. In this paper we will always show the 
 square of the modulus  of the complex solution for ease of 
 presentation. The time dependence of the $L_{2}$ norm of the difference between the 
numerical and the exact solution can be also seen there.   

\begin{figure}[htb!]
\centering
\includegraphics[width=0.45\textwidth]{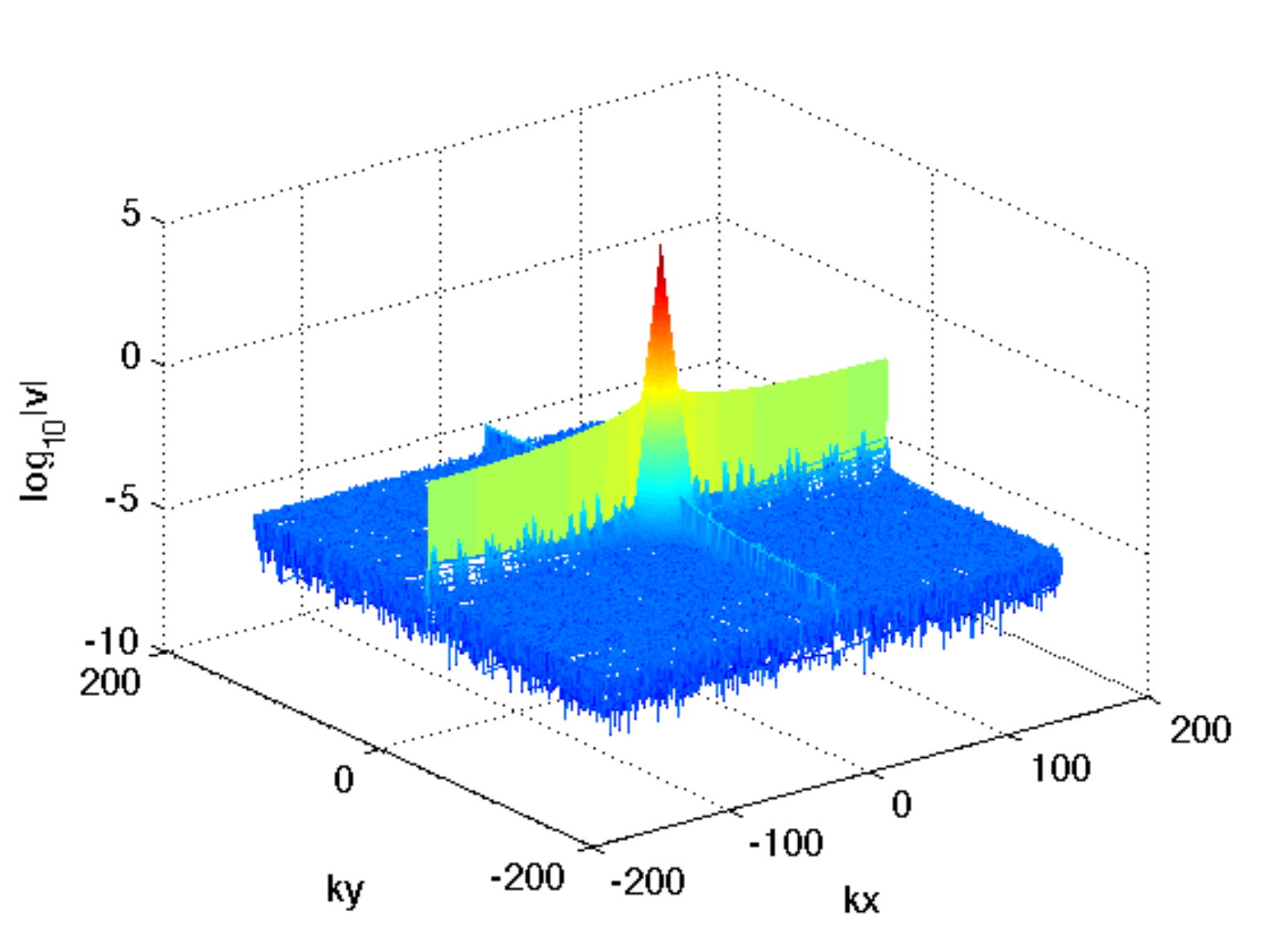}
\includegraphics[width=0.45\textwidth]{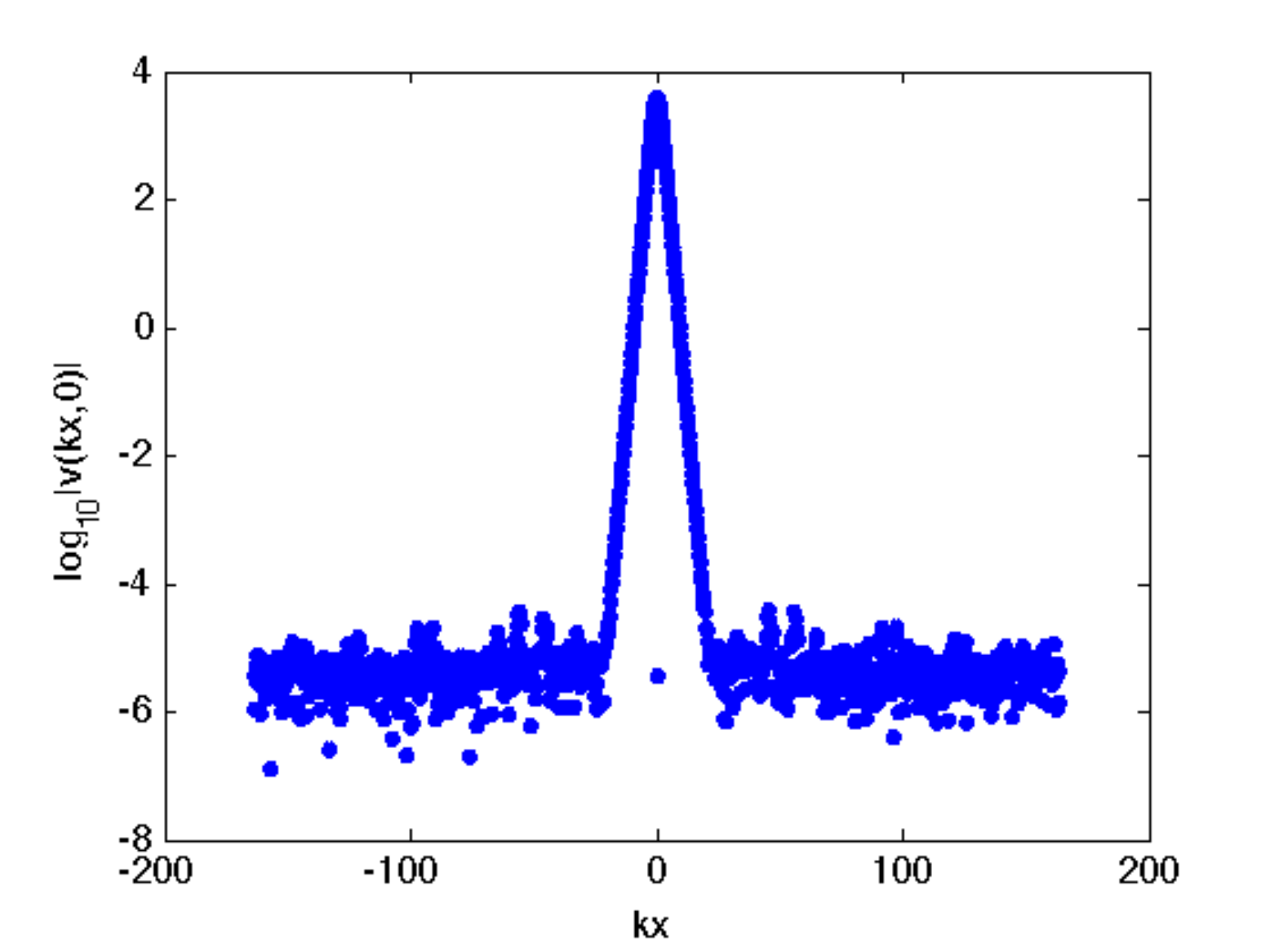}
\caption{Fourier coefficients for the situation in Fig.~\ref{ratcont} at  $t=6$.} 
\label{ratcf}
\end{figure}

The numerical error 
is here mainly due to the lack of resolution in time. Since the 
increase in the number of time steps is computationally expensive, 
a fourth order scheme is very useful in this context. The spatial 
resolution can be seen from the modulus of the
 Fourier coefficients at the final time of computation 
$t=6$ in Fig. \ref{ratcf}. It decreases to $10^{-6}$, thus 
essentially the value for the initial data. For 
computational speed considerations we always use double precision 
which because of finite precision arithmetic give us a range of 15 
orders of magnitude. Since function values computed using the split 
step method were for most of the computation  of order 1, and less 
than 5,000, rounding errors allow for a precision of $10^{-14}$ when 
less than $2^{15}\times 2^{15}$ Fourier modes are used. When more modes than 
$2^{15}\times 2^{15}$ were used, we found a reduction in 
precision.  Despite the 
algebraic falloff of the solution we have a satisfactory spatial 
resolution because of the large computational 
domain and the high resolution. The modulational instability does not 
show up in this and later examples before blowup.

\section{Blowup for the quintic NLS in $1+1$ dimensions and the focusing DS II}
It is known that focusing NLS equations can have solutions with 
blowup, if the nonlinearity exceeds a 
critical value depending on the spatial dimension. For the $1+1$ dimensional case, 
the critical nonlinearity is quintic, for the $2+1$ dimensional it is 
cubic, see for instance \cite{DSnls} and references therein. Thus the 
focusing DS II equation can have solutions with blowup. In this 
section we will first study numerically blowup for the $1+1$ dimensional quintic NLS 
equation, and then numerically evolve initial data for a known exact blowup solution to 
the focusing DS II equation due to Ozawa \cite{Oza}. We discuss some 
peculiarities of the fourth order splitting scheme in this context.

\subsection{Blowup for the quintic one-dimensional NLS}
The focusing quintic NLS in $1+1$ dimensions has the form
\begin{equation}
    i\partial_{t} u+\partial_{xx}u+ |u|^{4}u=0
    \label{quintic},
\end{equation}
where $u\in \mathbb{C}$ depends on $x$ and $t$ (we consider again 
solutions periodic in $x$). This equation is not 
completely integrable, but assuming the solution is in 
$L_{2}$, has conserved $L_{2}$ norm and, provided the 
solution $u \in H_{2}$,
a conserved energy,
\begin{equation}
    E[u] = 
    \int_{\mathbb{R}}^{}\left(\frac{1}{2}|\partial_{x}u|^{2}-\frac{1}{6}|u|^{6}\right)dx
    \label{quintener}.
\end{equation}
It is known that initial data with negative energy blow up for this 
equation in finite time, and that the behavior close to blowup is given in terms 
of a solution to an ODE, see \cite{MR}.

As discussed in sect.~2.1, the splitting scheme we are using here has the property that the
$L_{2}$ norm is conserved. Thus the quality of the 
numerical conservation of the $L_{2}$ norm gives no indication on the 
accuracy of the numerical solution. However as discussed in 
\cite{ckkdvnls}, conservation of the numerically computed energy gives a 
valid criterion for the quality of the numerics: in general it 
overestimates the $L_{\infty}$ numerical error by two orders of 
magnitude at the typically aimed at precisions.

\begin{figure}[htb!]
\centering
\includegraphics[width=0.45\textwidth]{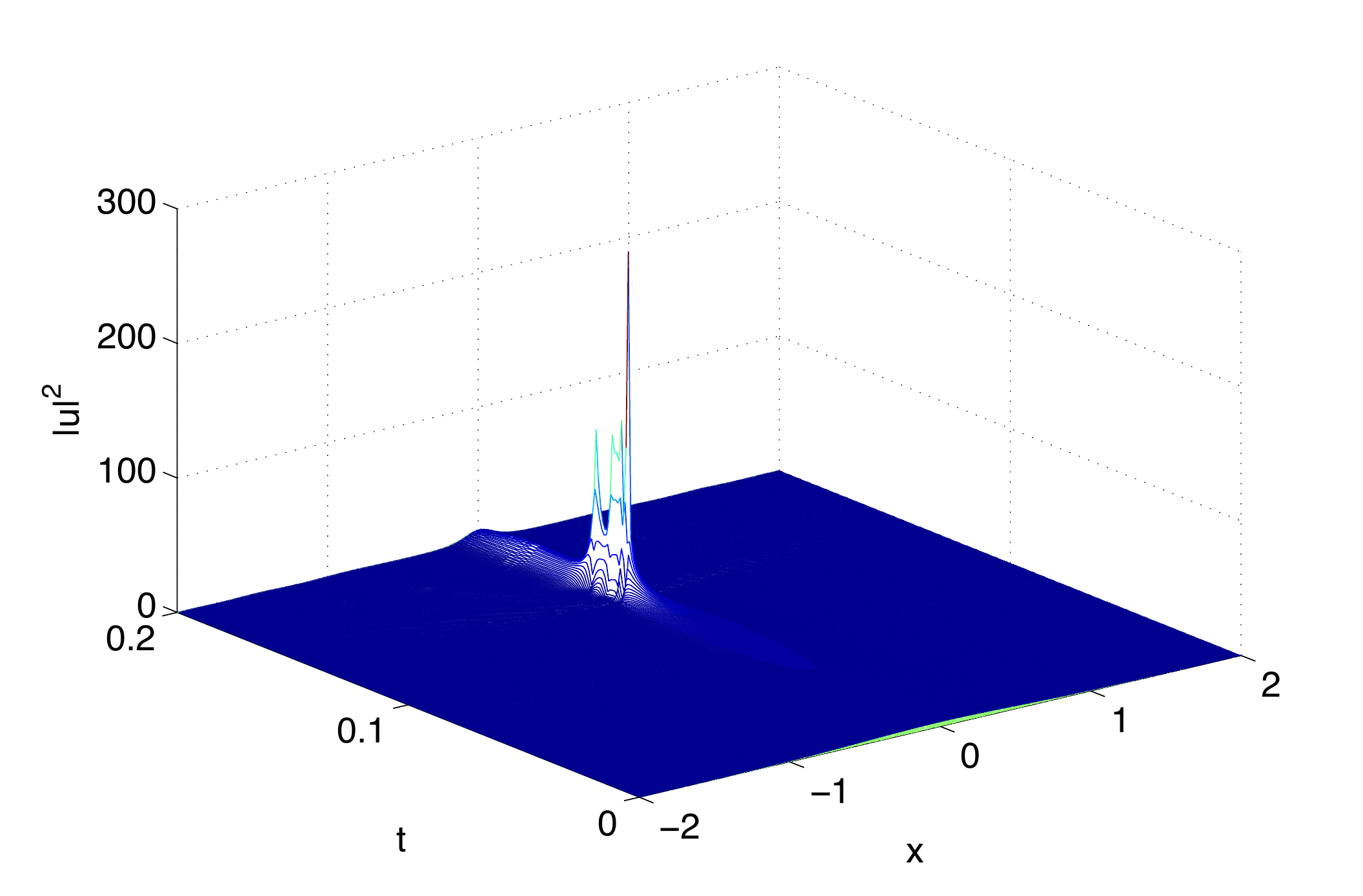}
\includegraphics[width=0.45\textwidth]{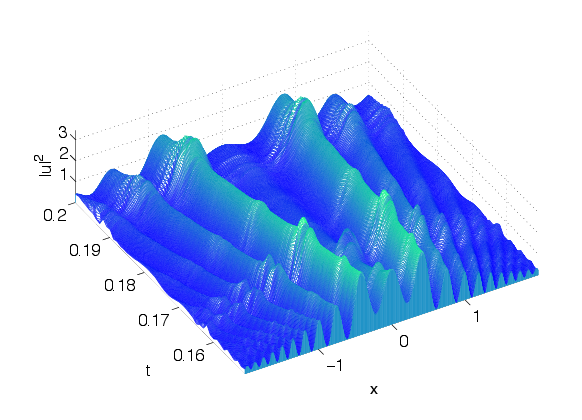}
\caption{Solution to the focusing quintic NLS (\ref{quintic}) for 
the initial data $u_{0}=1.8i\exp(-x^{2})$ with $N=2^{12}$ on the left 
and $N=2^{15}$ on the right for $t>t_{c}$.} 
\label{quintsol}
\end{figure}

\begin{figure}[htb!]
\centering
\includegraphics[width=0.45\textwidth]{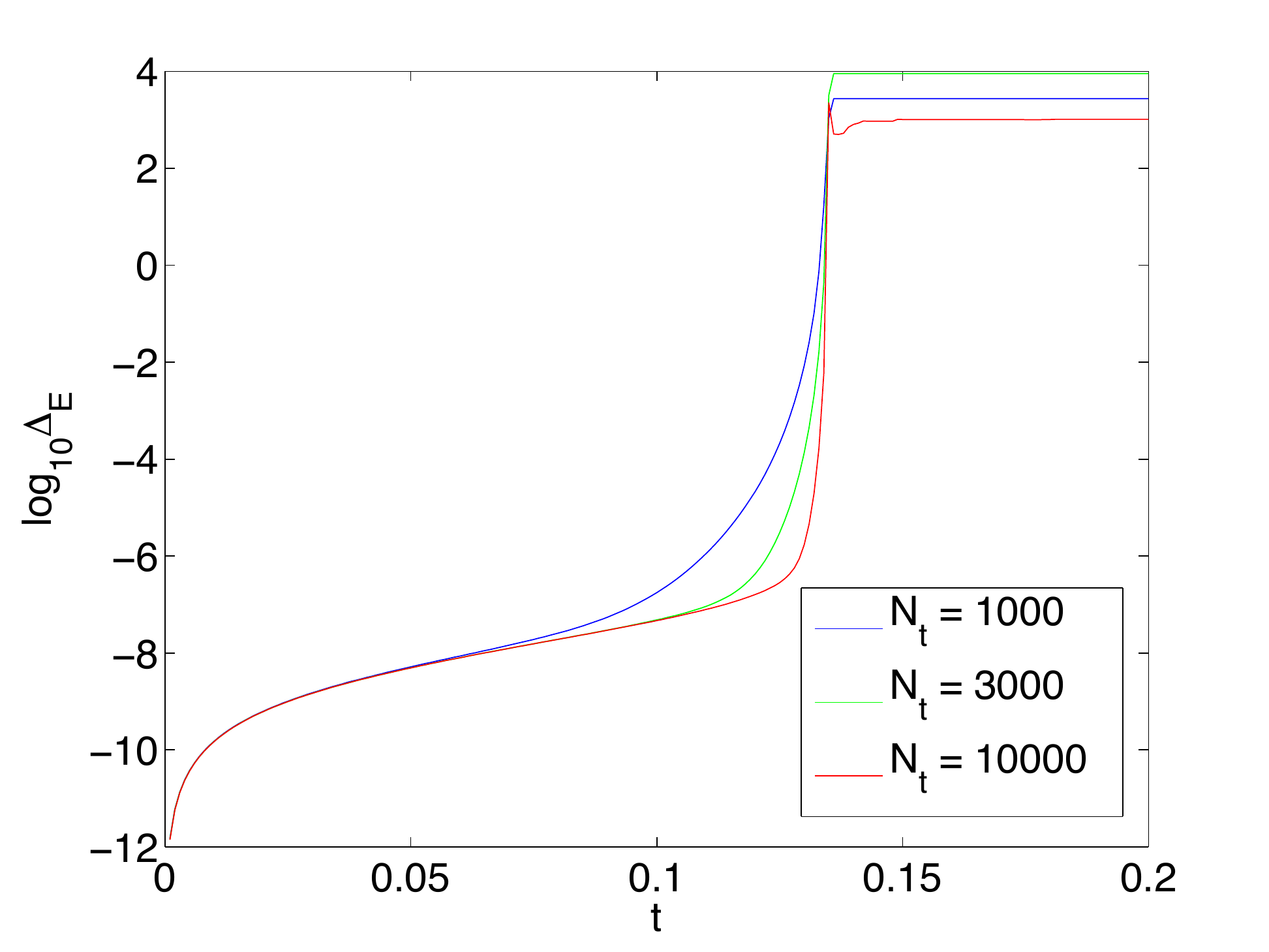}
\includegraphics[width=0.45\textwidth]{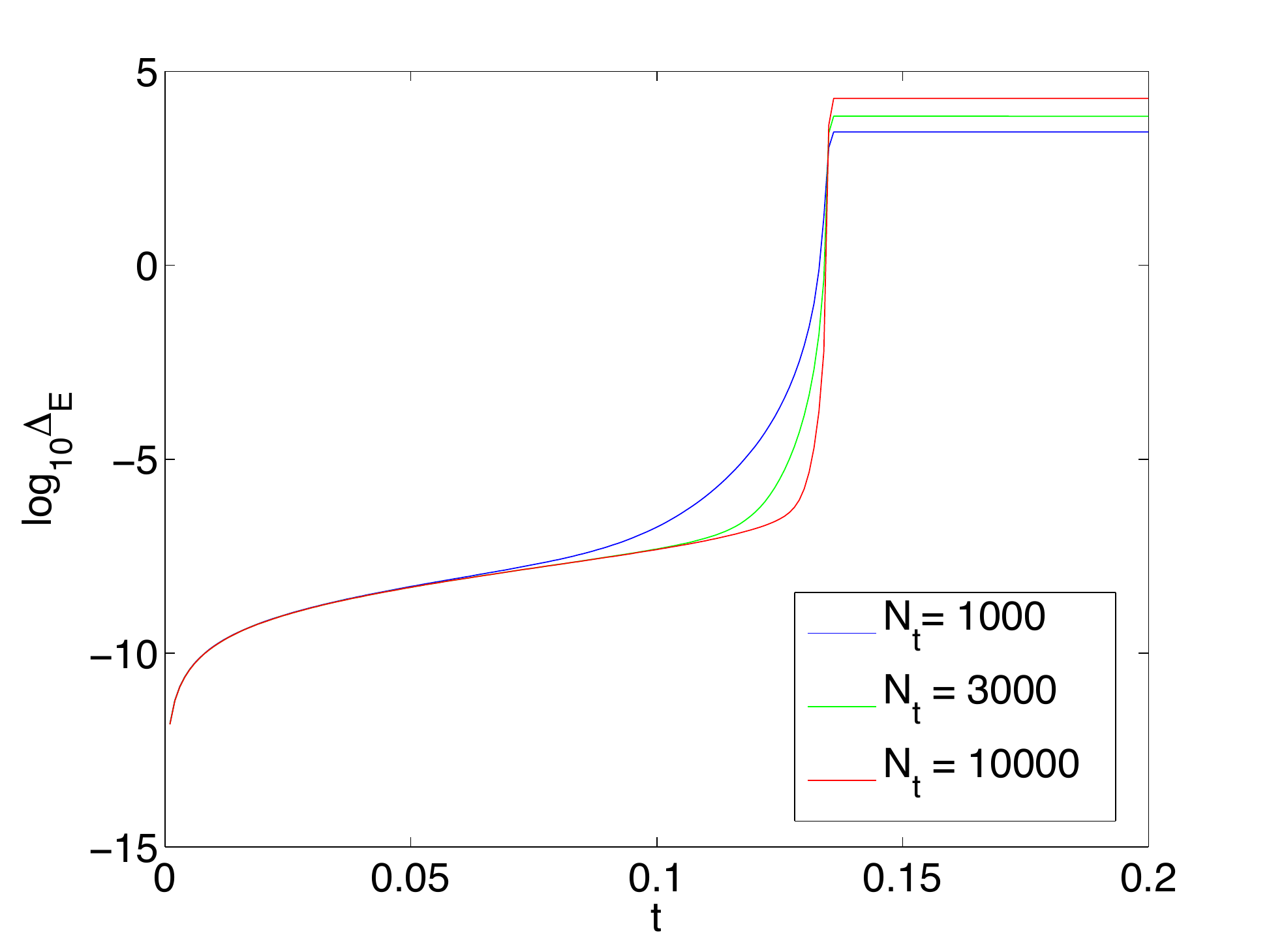}
\caption{Numerically computed energy for the situation studied in 
Fig.~\ref{quintsol} for $N=2^{12}$ on the left and $N=2^{15}$ on the 
right for several values of $N_{t}$. At the blowup, the energy jumps. } 
\label{nlsquintE}
\end{figure}

If we consider as in \cite{Sti} for the quintic NLS the initial data 
$u_{0}(x)=1.8i\exp(-x^{2})$, the energy is negative. We compute the 
solution with $L_{x}=5$ and $N_{x}=2^{15}$ with $N_{t}=10^{4}$ time 
steps. The result can be seen in Fig.~\ref{quintsol} (to obtain more 
structure in the solution after the blow up due to a less pronounced 
maximum, the plot on the left was generated 
with the lower spatial resolution $N=2^{12}$).
The initial data clearly get focused to a strong maximum, but the 
code does not break. We note that this is in contrast to other fourth order 
schemes tested for $1+1$ dimensional NLS equations in 
\cite{ckkdvnls}, which typically produce an overflow close to the 
blowup. But clearly the 
solution shows spurious oscillations after the time $t_{c}\sim0.155$.
In fact the numerically computed energy, which will always be 
time-dependent due to unavoidable numerical errors, will be 
completely changed after this time. We consider 
\begin{equation}
    \Delta_{E}=\left|1-\frac{E(t)}{E(0)}\right|
    \label{DeltaE},
\end{equation}
where $E(t)$ is the numerically computed 
energy (\ref{quintener}) and get for the example in 
Fig.~\ref{quintsol} the behavior shown in Fig.~\ref{nlsquintE}. At 
the presumed blowup at $t_{c}\sim0.155$ as in \cite{Sti}, the energy jumps to a completely different 
value. Thus this jump can and will be used to indicate 
blowup. To illustrate the effects of a lower resolution in time and 
space imposed by hardware limitations for the DS 
computations, we show this quantity for several resolutions in 
Fig.~\ref{nlsquintE}. If a lower resolution in time is used as in some of the 
DS examples in this paper, the jump is slightly smoothed out. But the 
plateau is still reached at essentially the same time which indicates 
blowup. Thus  a lack of resolution in time in the 
given limits will not be an obstacle to identify a possible 
singularity. The reason for this is the use of a fourth order 
scheme that allows to take larger time steps. We will present 
computations with different resolutions to illustrate the steepening 
of the energy jump as above if this is within the limitations imposed 
by the hardware.

\begin{figure}[htb!]
\centering
\includegraphics[width=0.45\textwidth]{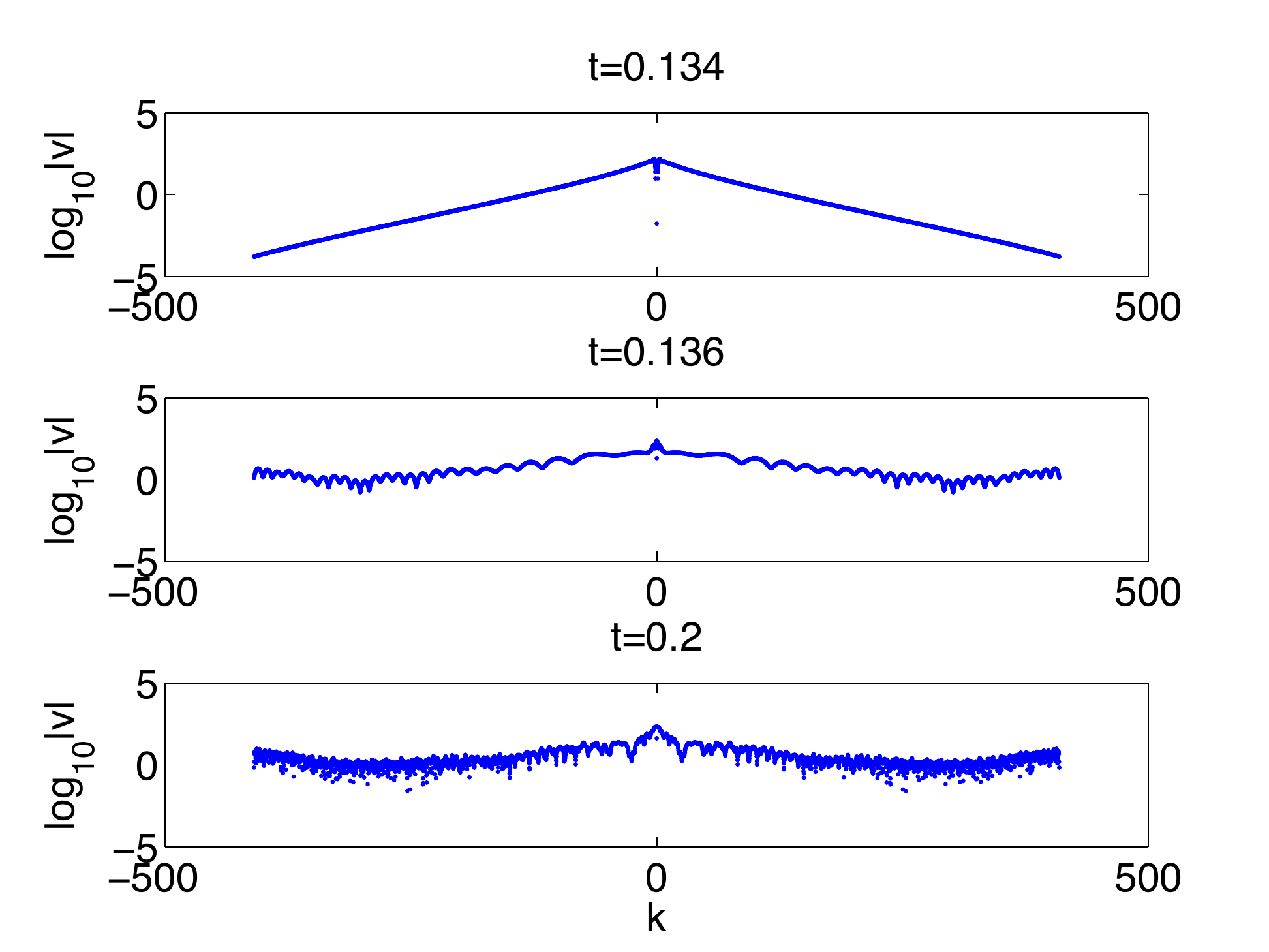}
\includegraphics[width=0.45\textwidth]{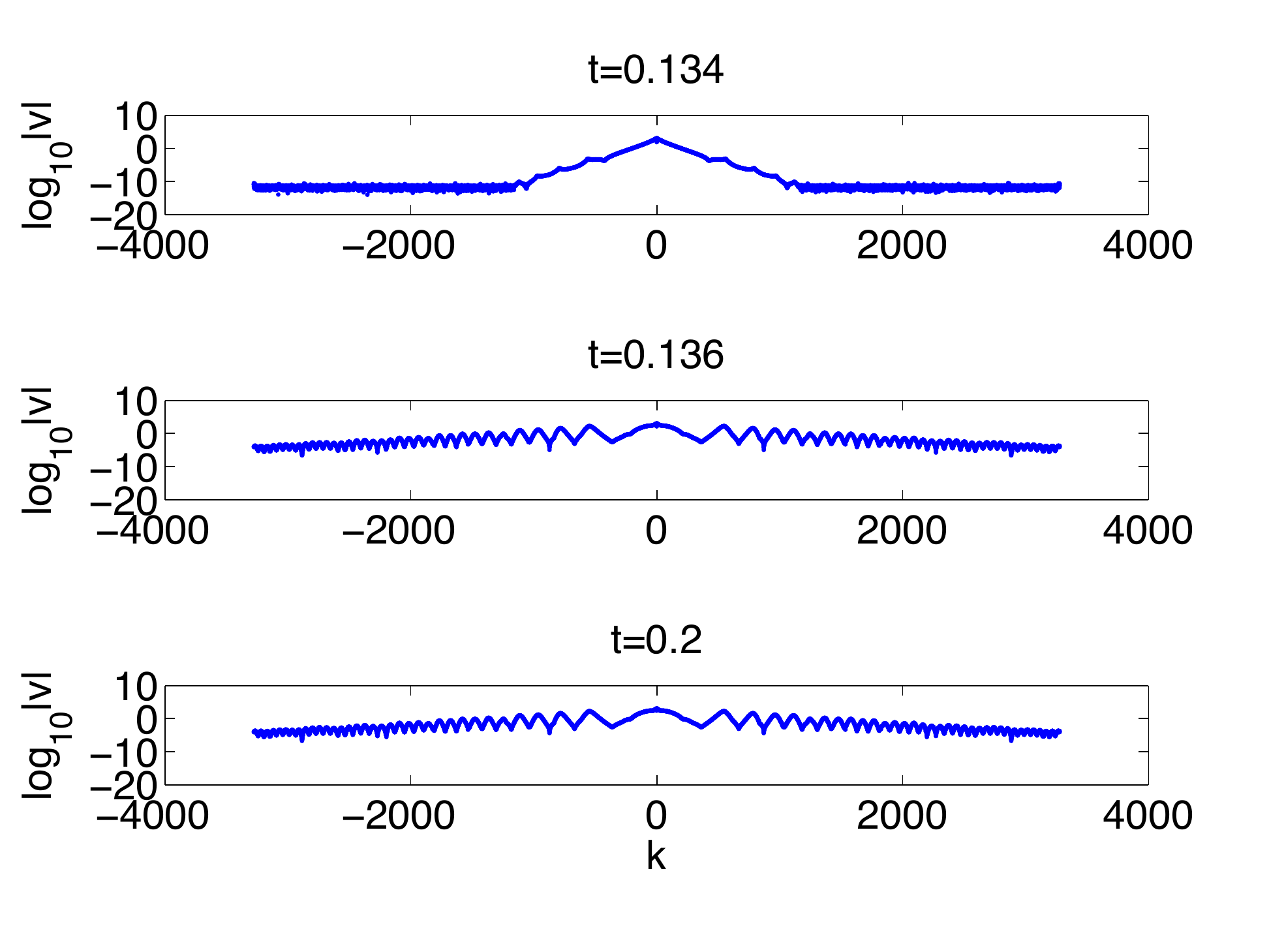}
\caption{Fourier coefficients for the solution in 
Fig.~\ref{quintsol} close to the critical and at a later time for $N=2^{12}$ on the left and $N=2^{15}$ on the 
right for $N_{t}=10^{4}$.} 
\label{nlsquintfourier}
\end{figure}

We show the modulus of the Fourier coefficients for $N=2^{12}$ and $N=2^{15}$ 
before and after the critical time in Fig.~\ref{nlsquintfourier}. It can be 
seen that the solution is well resolved before blowup in the latter 
case, and that the
singularity leads to oscillations in the Fourier 
coefficients. A lack of spatial resolution as for $N=2^{12}$ in 
Fig.~\ref{nlsquintfourier} triggers the modulation instability close 
to the blowup and at later times as can be seen from the modulus of 
the Fourier 
coefficients that increase for larger wavenumbers. Therefore we 
always aim at a sufficient resolution in space even for times close 
to a blowup. After this time the modulation instability will be 
present in the spurious solution produced by the splitting scheme as 
we will show for an example.

\begin{remark}
    Stinis \cite{Sti} has recently computed singular solutions to the 
    focusing quintic nonlinear Schr\"odinger equation in $1+1$ 
    dimensions. This equation has solutions in $L_{\infty}L_{2}$ 
    that  may not be unique for given smooth initial data and that 
    may exhibit blowup of the $L_{\infty}H_{1}$ norm. Following Tao 
    \cite{Tao}, Stinis \cite{Sti} has used a selection criteria to pick a 
    solution after the blow up time of the $L_{\infty}H_{1}$ norm. 
    They suggest that `mass' is ejected (which means that the $L_{2}$ 
    norm is changed) at times where the 
    $L_{\infty}H_{1}$ norm blows up. The splitting scheme studied 
    here in contrast produces a weak solution with a different energy 
    since the $L_{2}$ norm conservation is built in. 
\end{remark}

\subsection{Blowup in the Ozawa solution}
For the focusing DS II equation, an exact solution was given by Ozawa 
\cite{Oza} which is in $L_{2}$ for all times with an $L_{\infty}$ blowup in finite time. 
We can  summarize his results as follows:
\begin{theorem}[Ozawa] 
Let   $ab<0$ and $T=-a/b$. 
Denote by $u (x,y,t)$ the function defined by
\begin{equation}
 u (x,y,t) = \exp \left( i \frac{b}{4(a+bt)} (x^2 - y^2)
\right) \frac{v(X,Y)}{a+bt}
\label{soluoz}
\end{equation}
 where 
\begin{equation}
 v(X,Y) = \frac{2}{1+X^2+Y^2}, \,\, X=\frac{x}{a+bt},
\,\, Y=\frac{y}{a+bt}
\end{equation}
Then, $u$ is a solution of (\ref{DSII}) with 
\begin{equation}
 \|u(x,y,t)\|_2 = \|v(X,Y)\|_2 = 2\sqrt{\pi}
\end{equation}
and
\begin{equation}
 |u(t)|^2 \rightarrow 2\pi \delta 
 \,\, \mbox{when} \,\, t \rightarrow T.
\end{equation}
where $\delta$ is the Dirac measure.
\end{theorem}

We thus consider initial data of the form
\begin{equation}
    u(x,y,0) = 2\frac{ \exp \left( -i (x^2 - y^2)
    \right)}{1+x^2 +y^2}
    \label{ozawaini}
\end{equation}
($a=1$ and $b=-4$ in (\ref{soluoz})). As for the 
quintic NLS in $1+1$ dimensions, we always 
trace the conserved energy for DS II (\ref{DSII}).

\begin{figure}[htb!]
\centering
\includegraphics[width=0.45\textwidth]{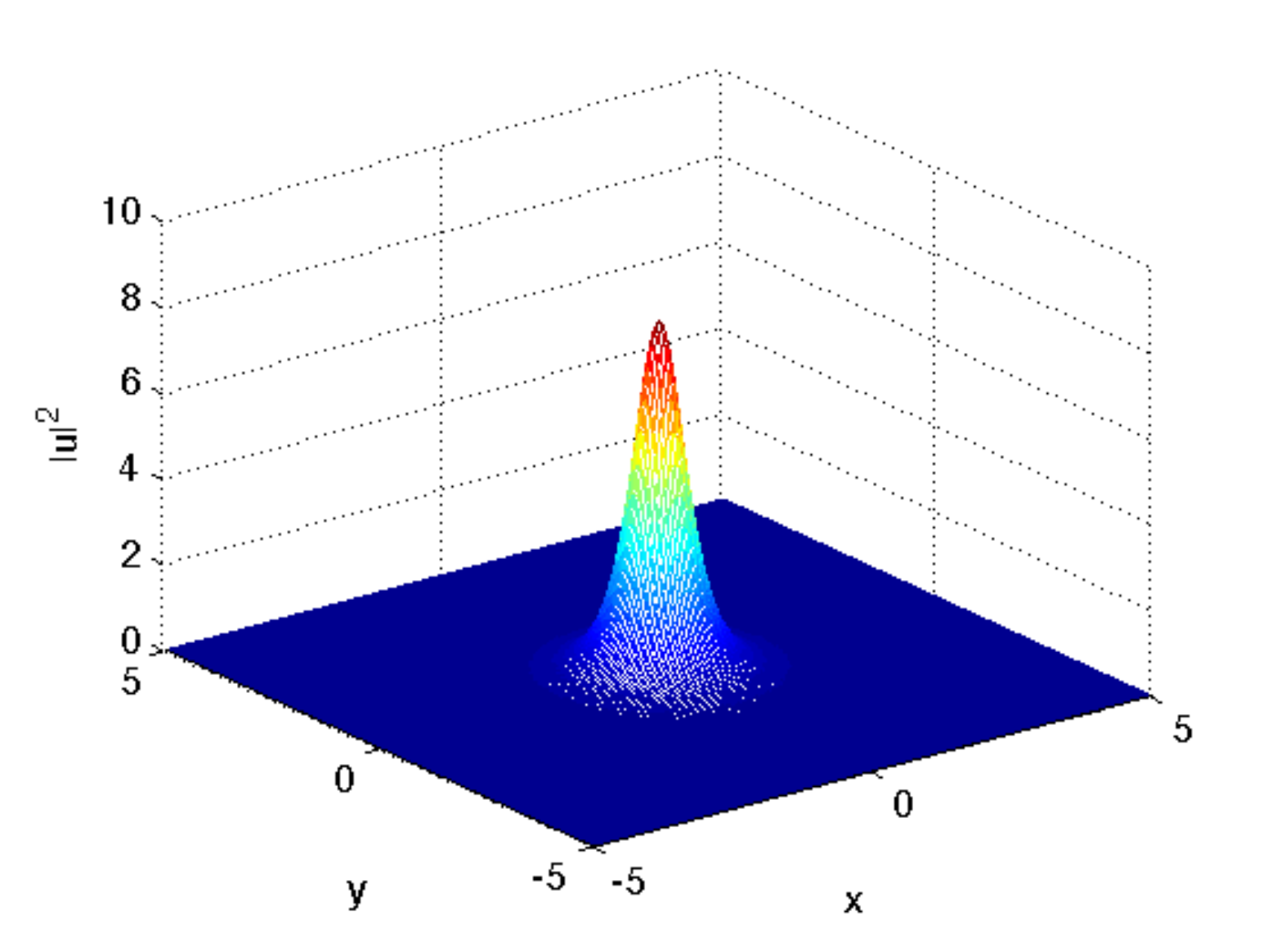}
\includegraphics[width=0.45\textwidth]{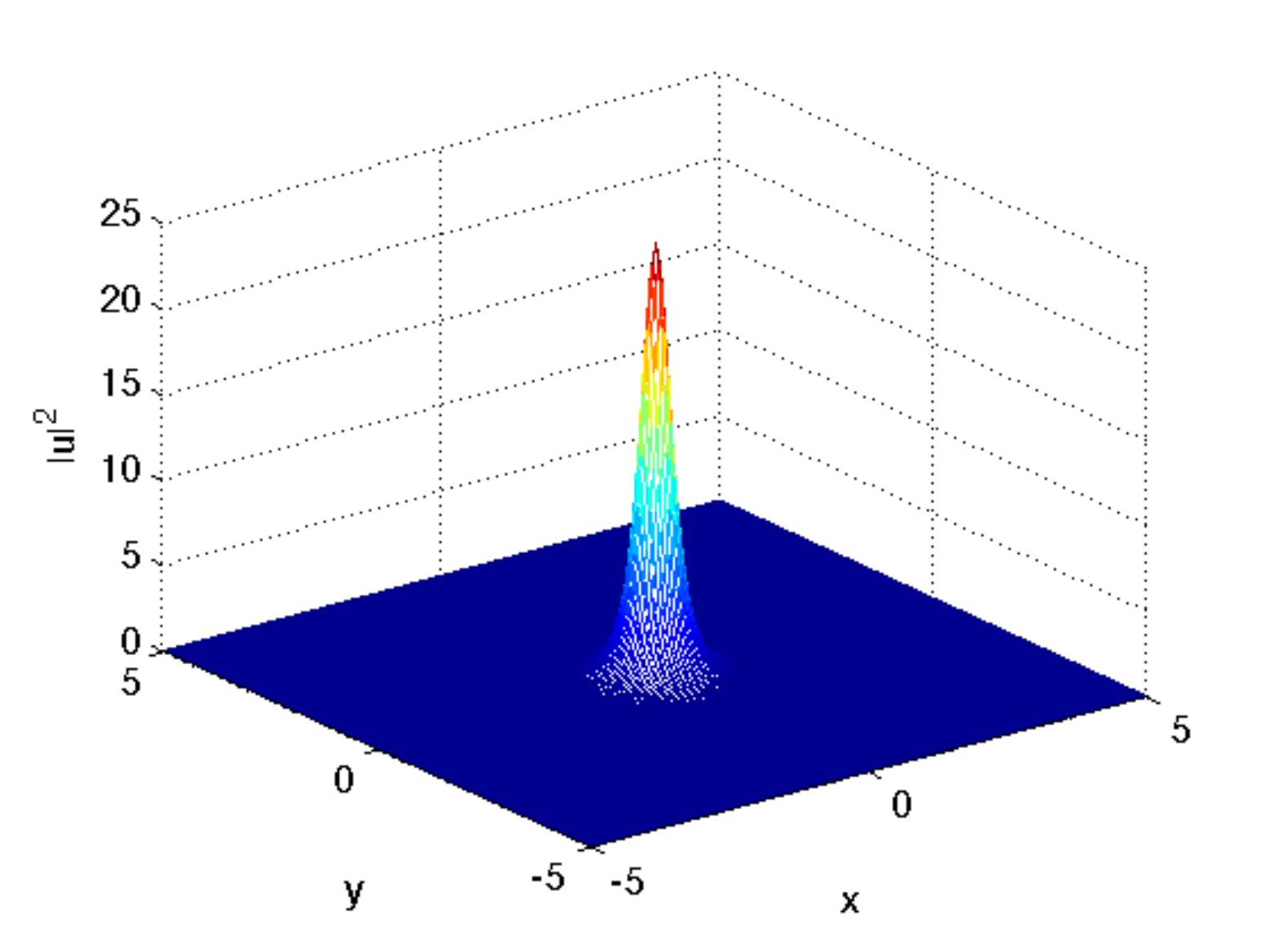}\\
\includegraphics[width=0.45\textwidth]{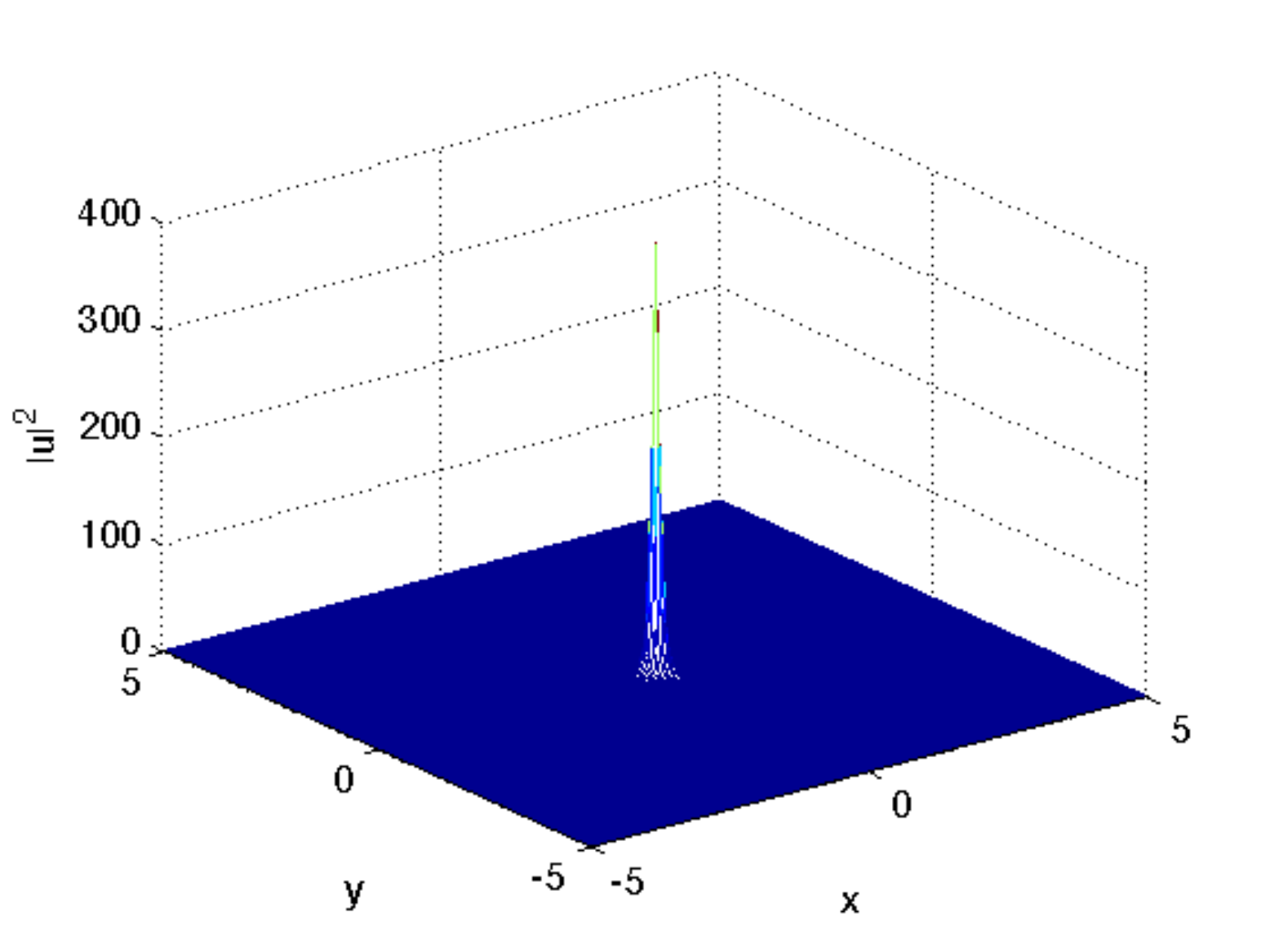}
\includegraphics[width=0.45\textwidth]{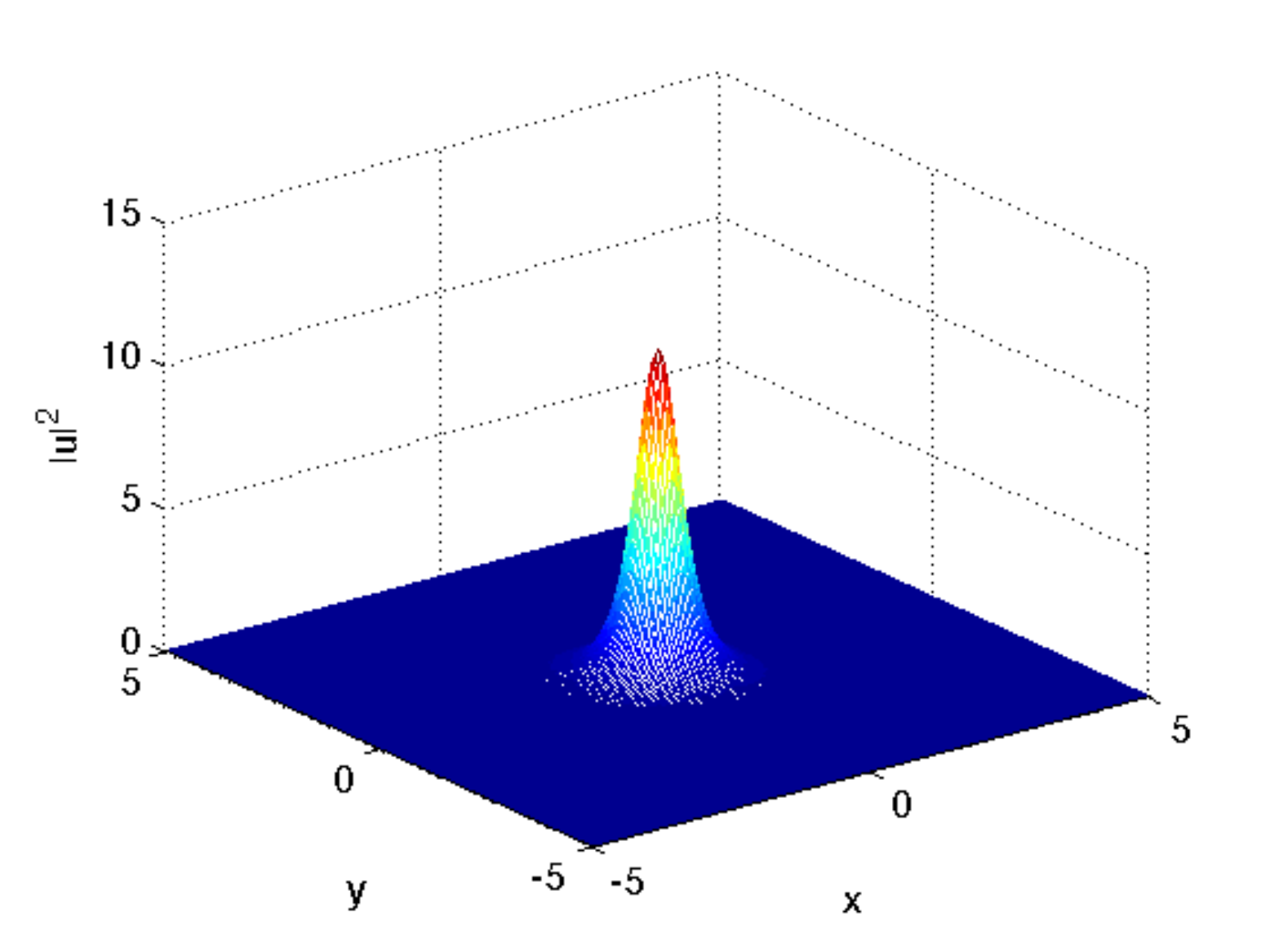}
\caption{Solution to the focusing DS II equation (\ref{DSII}) for 
$t=0.075$ and $t=0.15$ in the first row and $t=0.225$ and $t=0.3$ 
below 
for an initial condition of the form (\ref{ozawaini}).} 
\label{uoz}
\end{figure}

The computation is carried out with $N_{x}=N_{y}=2^{15}$, 
$L_{x}=L_{y}=20$, and $N_t=1000$ respectively $N_{t}=4000$;
we show the solution at different times
in Fig. \ref{uoz}. The difference of the Ozawa mass and the computed 
$L_{2}$ norm on the periodic setting is of the order of $9*10^{-5}$.

\begin{figure}[htb!]
\centering
\includegraphics[width=0.45\textwidth]{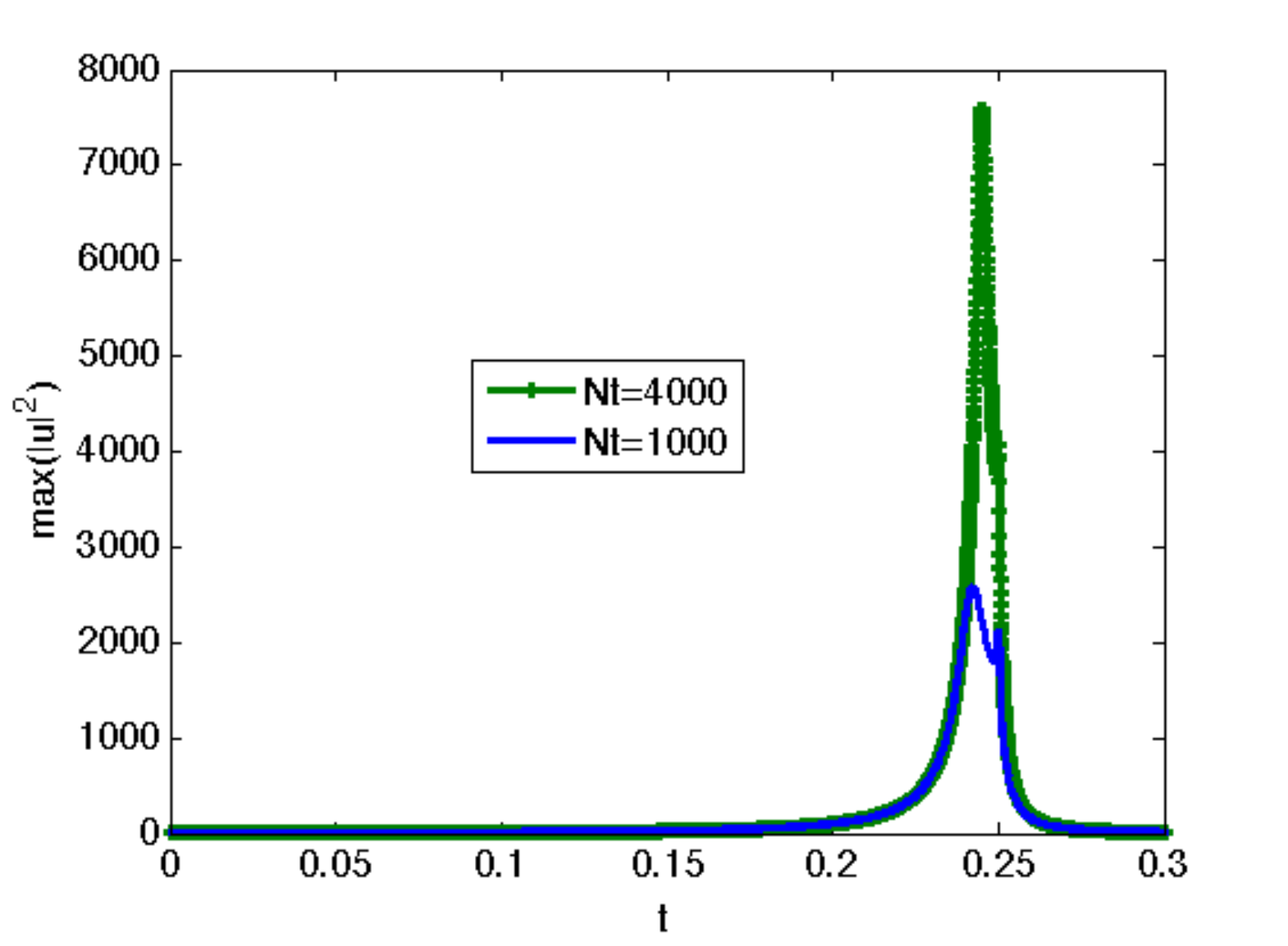}
\includegraphics[width=0.45\textwidth]{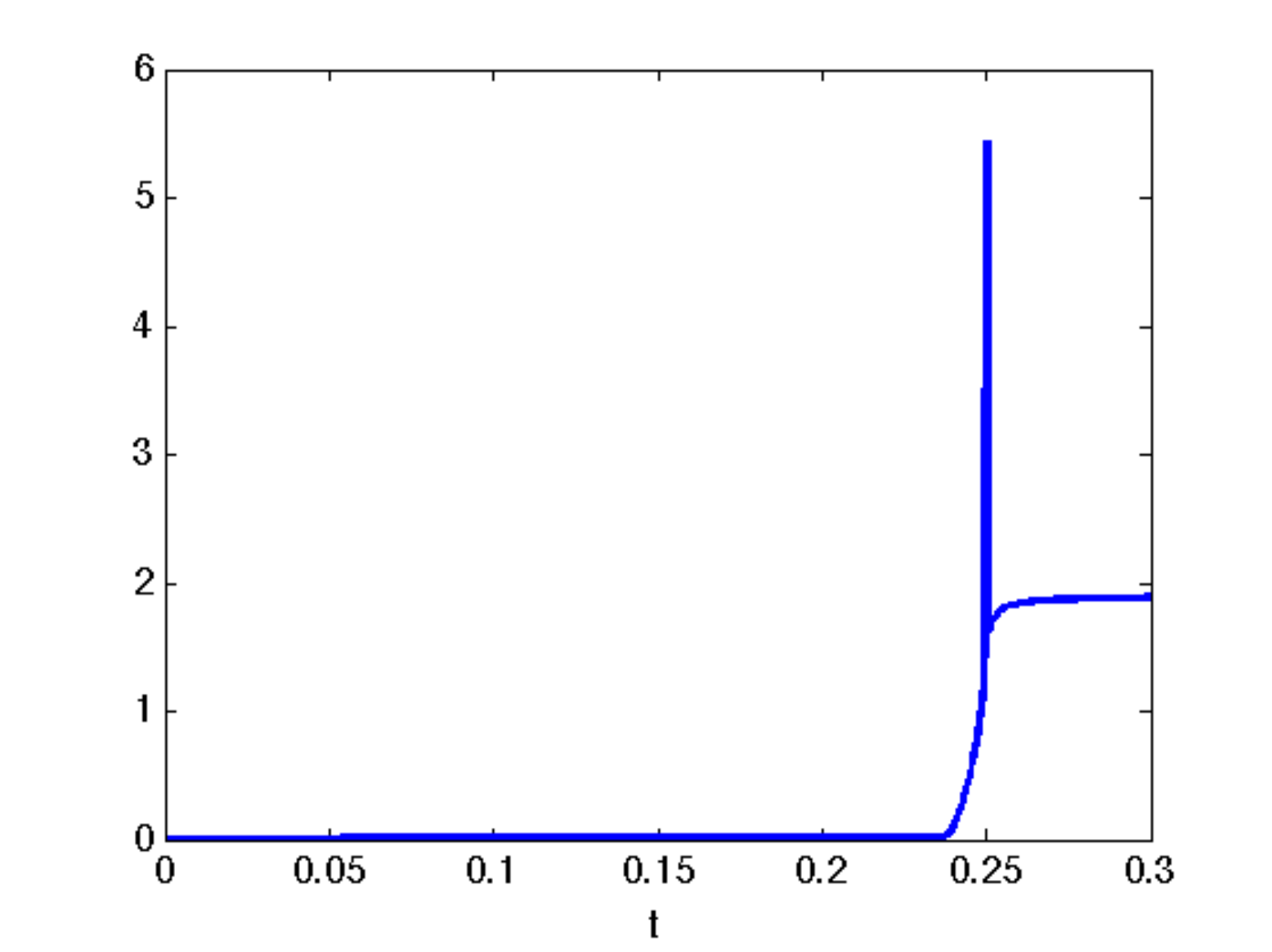}
\caption{Time evolution of $ max(| u_{num} |^2)$ and 
of $\| u_{num}-u_{exact} \|_2$ for the situation in Fig.~\ref{uoz}.}
\label{ampluoz}
\end{figure}

\begin{figure}[htb!]
\centering
\includegraphics[width=0.45\textwidth]{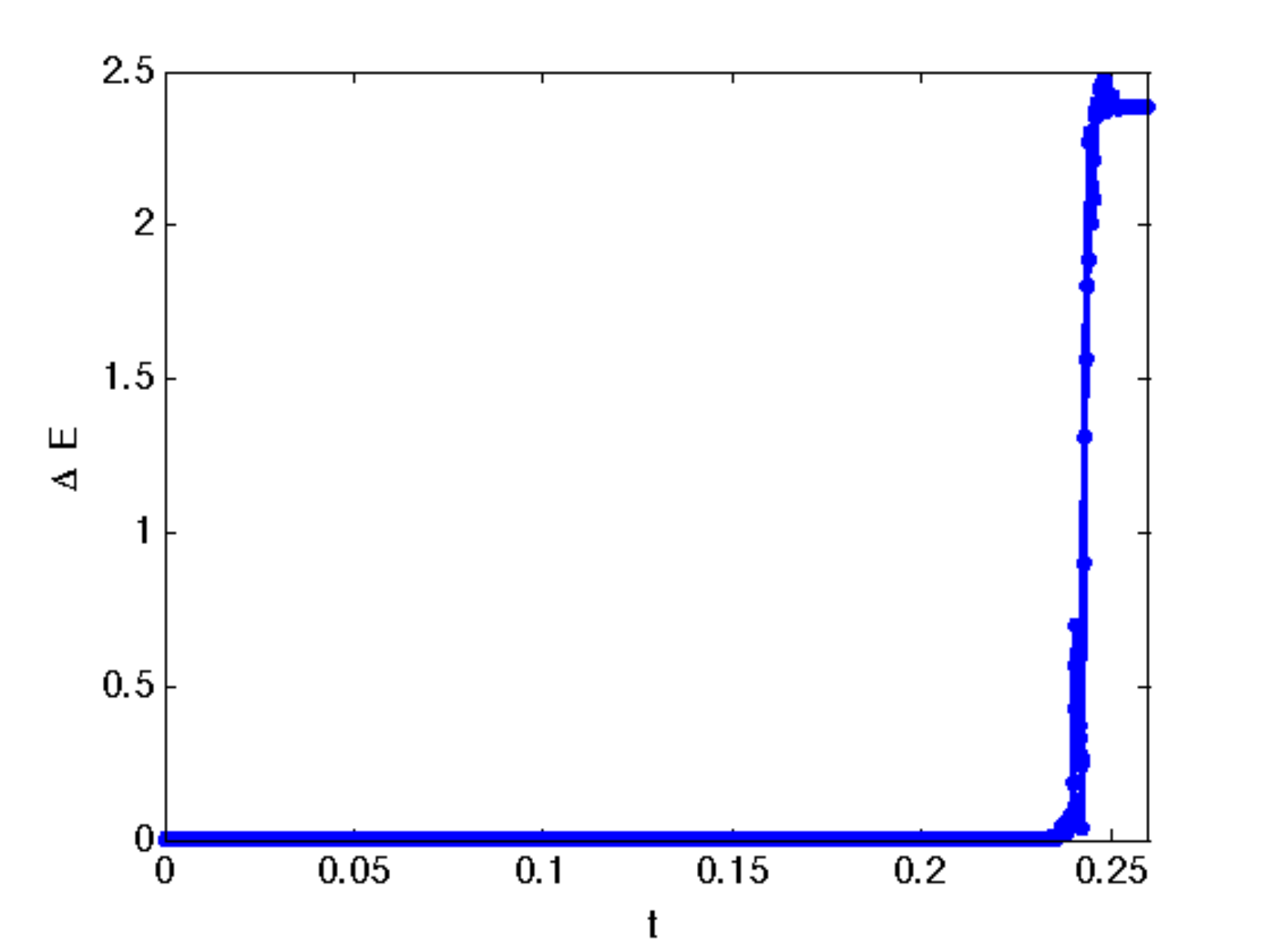}
\includegraphics[width=0.45\textwidth]{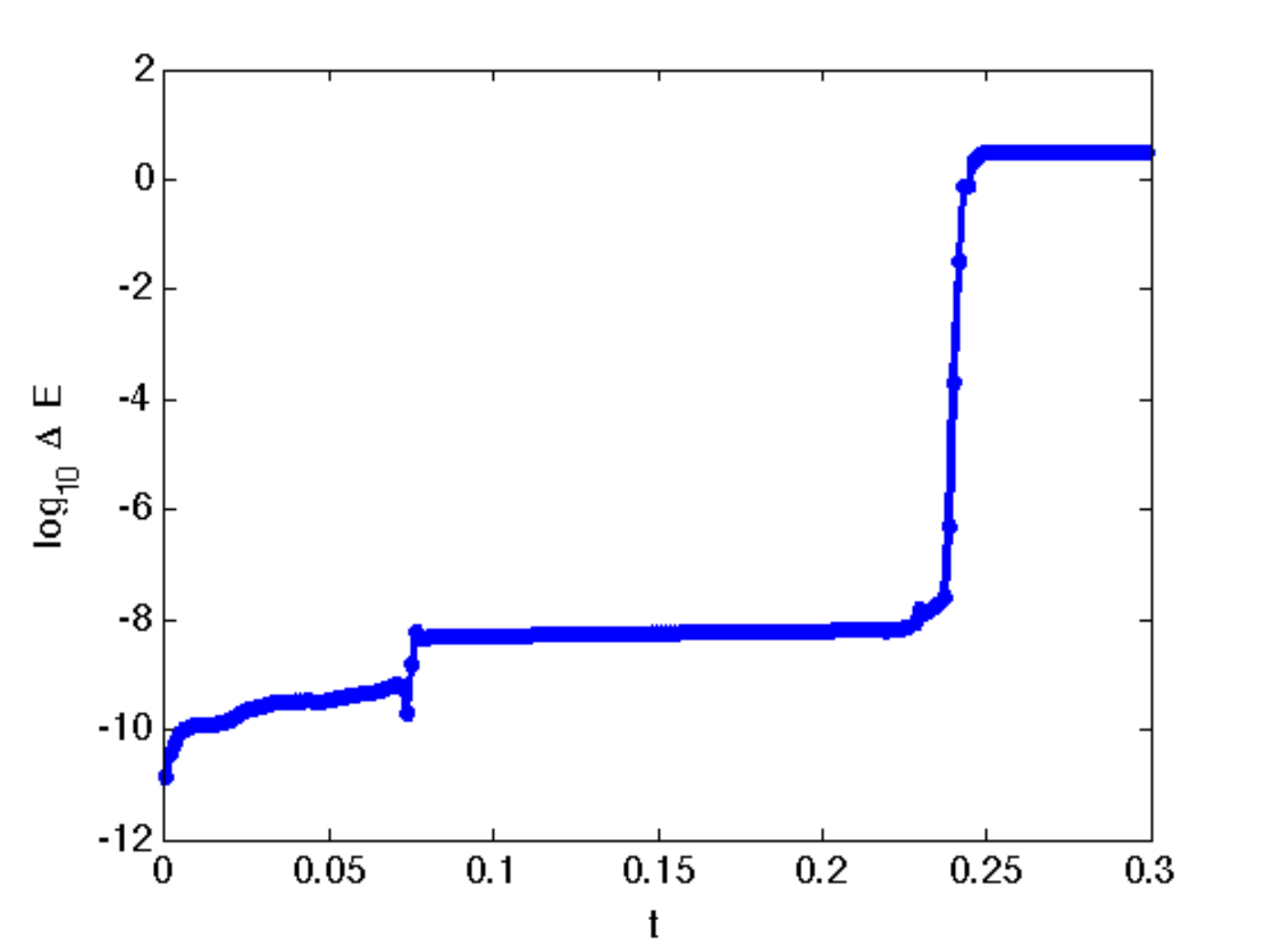}
\caption{Numerically computed energy $E(t)$ and $\Delta E = 
|1-E(t)/E(0)|$ (\ref{DeltaE})  for the situation in Fig.~\ref{uoz}.}
\label{uozenergy}
\end{figure}

\begin{figure}[htb!]
\centering
\includegraphics[width=0.45\textwidth]{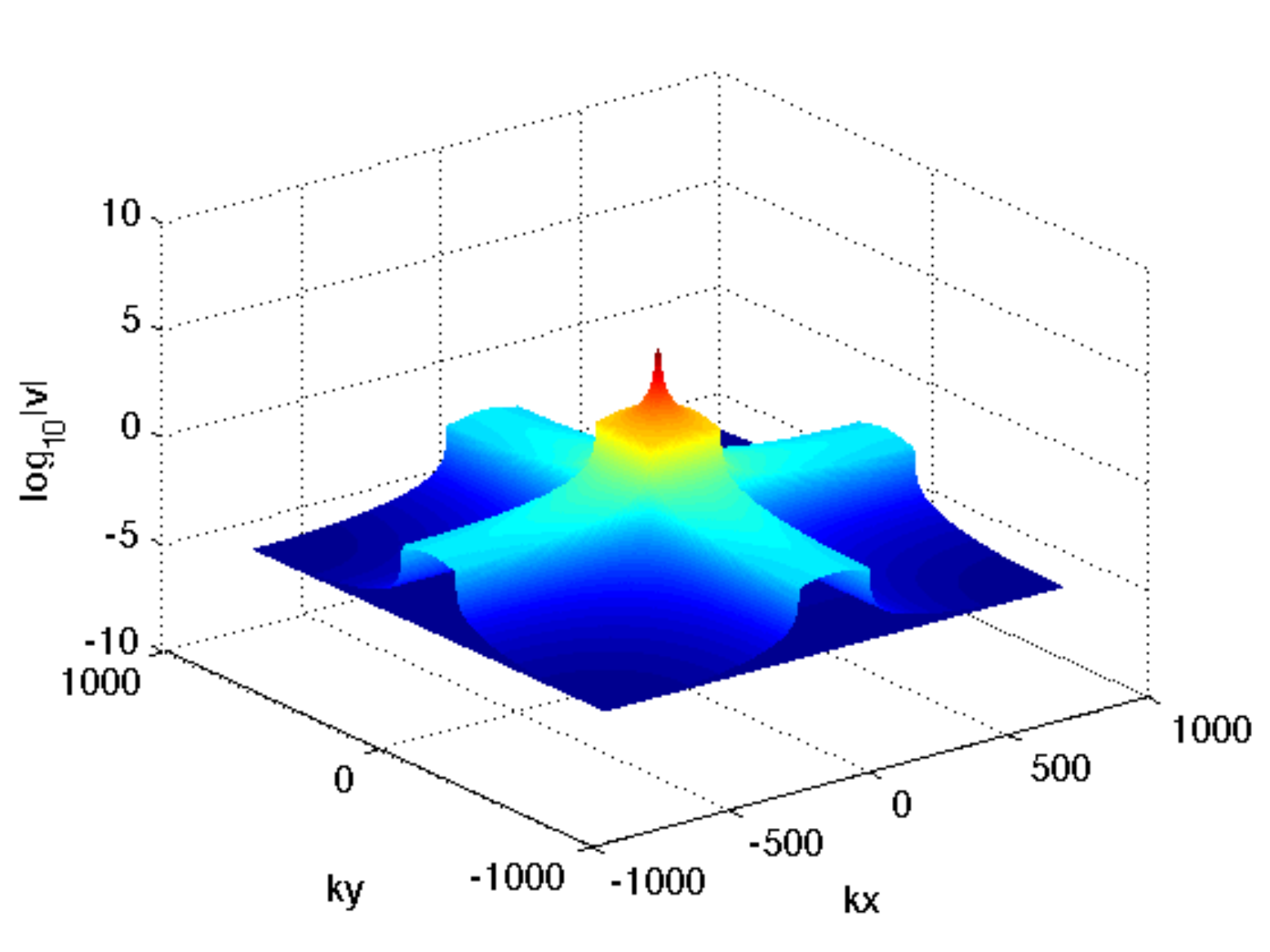}
\includegraphics[width=0.45\textwidth]{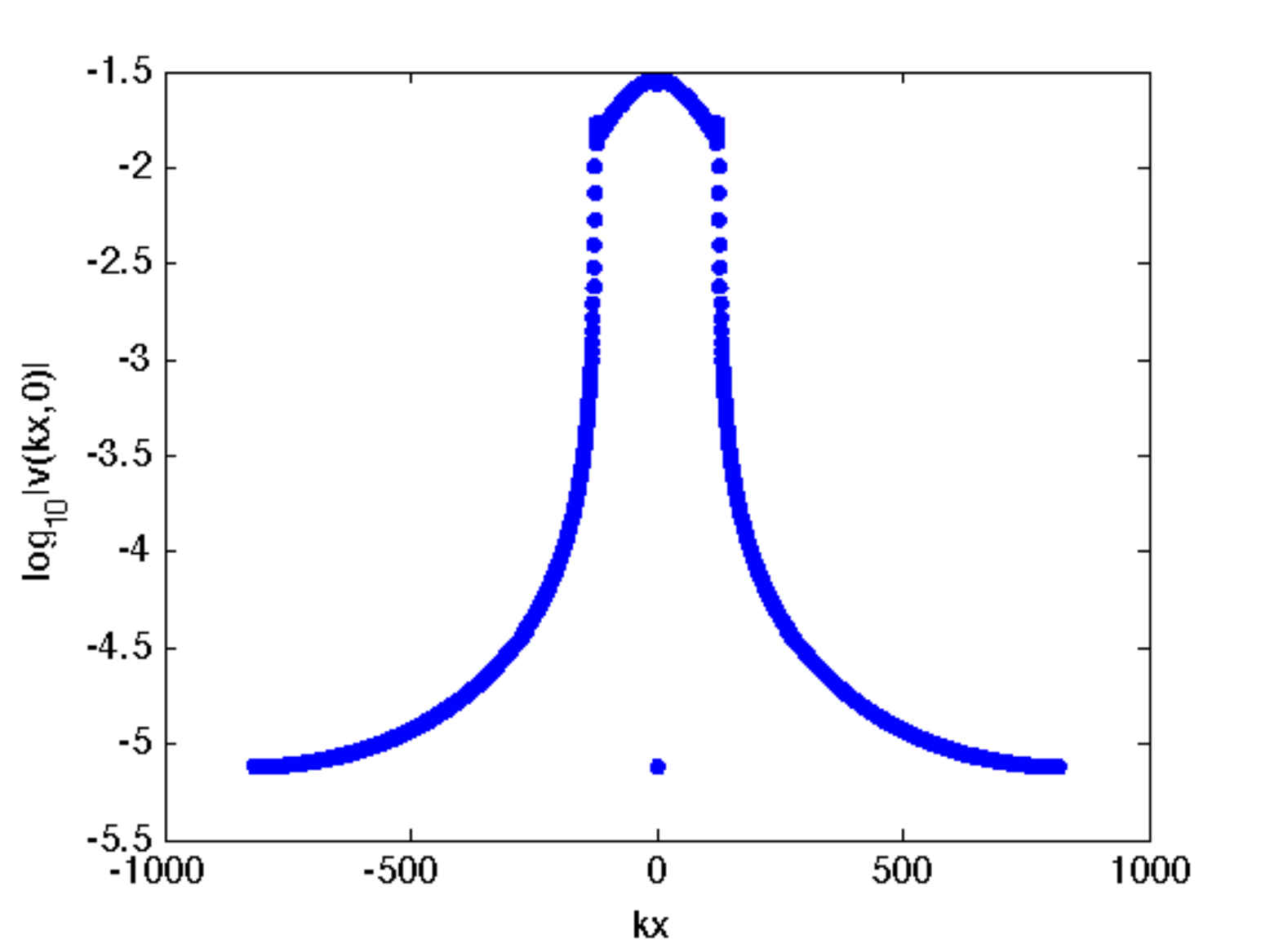}
\caption{Fourier coefficients of $u$ at $t=0.15$
for an initial condition of the form (\ref{ozawaini}).} 
\label{ozexcf}
\end{figure}

The time evolution of $\underset{x,y}{\max} |u(x,y,t)|^2$ and the 
difference between the numerical and the exact solution can be seen 
in Fig. \ref{ampluoz} (the critical time $t_{c}$ is not on the shown 
grid, thus the solution is always finite on the grid points).
The code continues to run after the critical time, but the numerical solution 
obviously no longer represents the Ozawa solution. 
The numerically computed energy jumps at the blow up time as can be 
seen in Fig.~\ref{uozenergy}.
The Fourier coefficients at $t=0.15$ are shown in 
Fig.~\ref{ozexcf}. Despite the Gibbs phenomenon the Fourier 
coefficients for the initial data decrease to $10^{-8}$. Spatial 
resolution is still satisfactory at half the blowup time. 

\begin{figure}[htb!]
\centering
\includegraphics[width=0.45\textwidth]{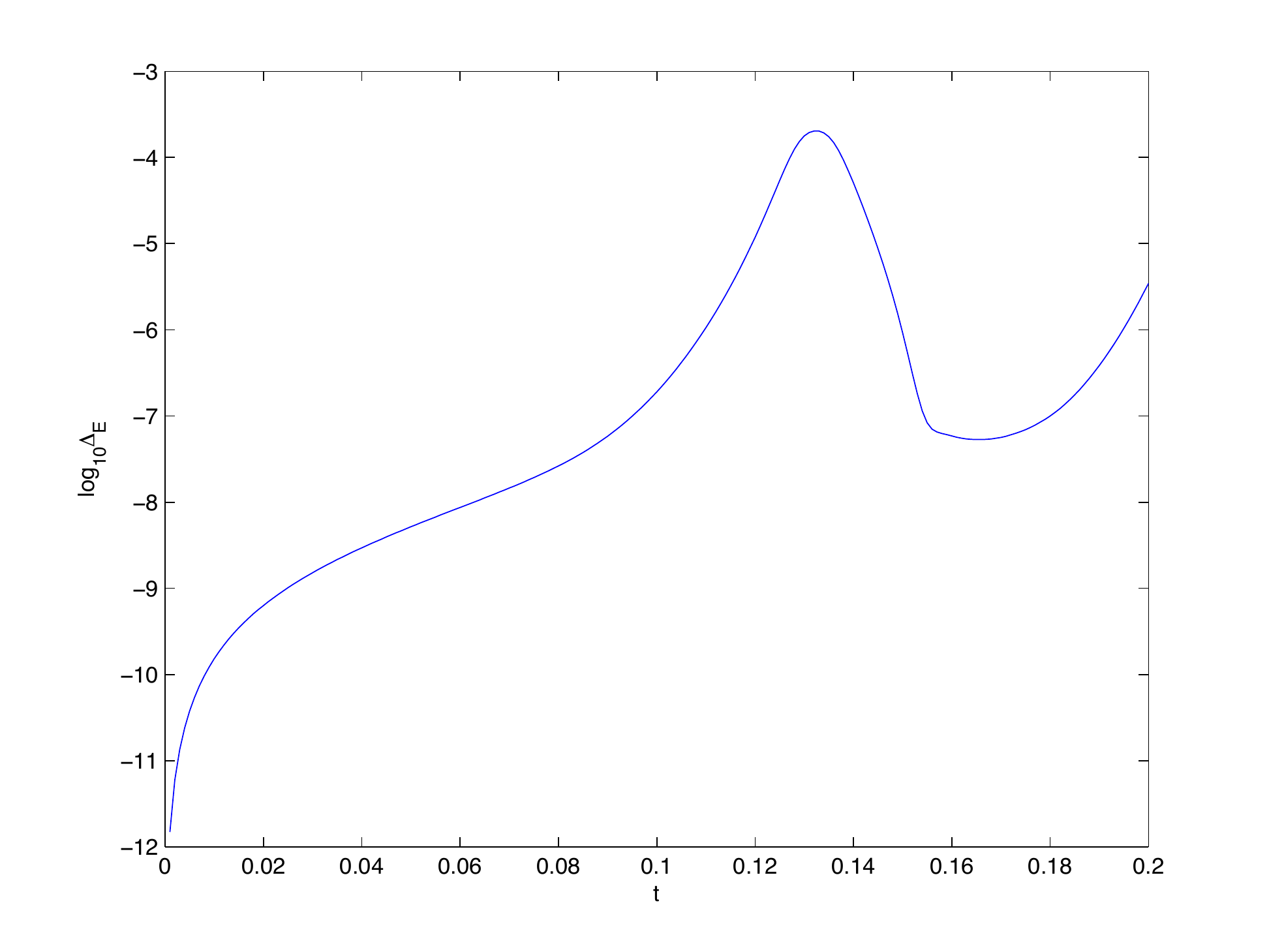}
\includegraphics[width=0.45\textwidth]{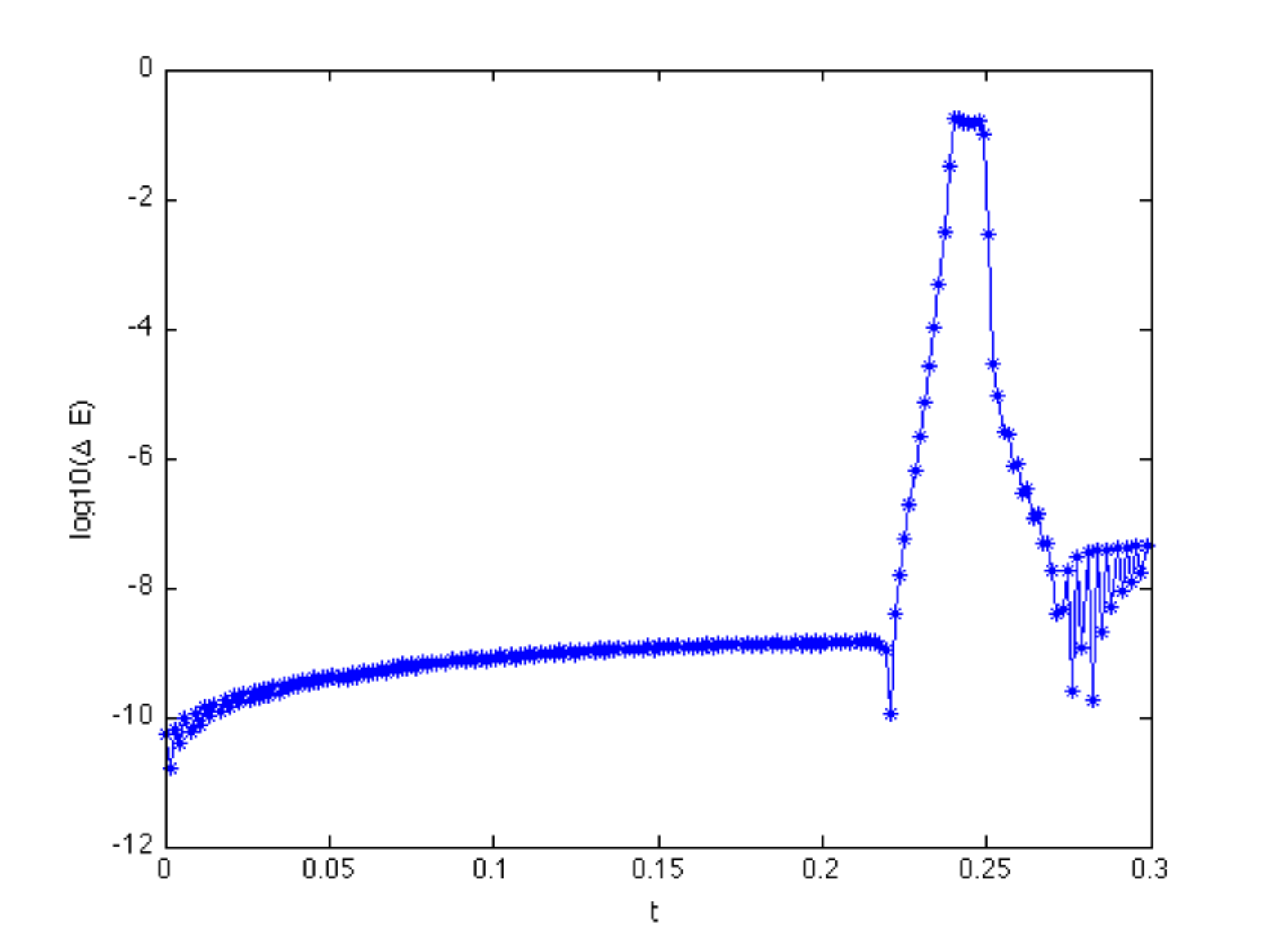}
\caption{Computed numerical energy for quintic NLS in 
Fig.~\ref{nlsquintE} with $N=2^{8}$ and for the Ozawa solution in 
Fig.~\ref{uozenergy} with $N_{x}=N_{y}=2^{12}$.} 
\label{nlsozres}
\end{figure}

\begin{remark}
    The jump of the computed energy at blowup is dependent on sufficient 
    spatial resolution as can be seen in Fig.~\ref{nlsozres} for the 
    example of the quintic NLS of Fig.~\ref{nlsquintE} and the Ozawa 
    solution in Fig.~\ref{uozenergy}. For low resolution blow-up can 
    be still clearly recognized from the computed energy, but the 
    energy does not stay on the level at blow-up.  
\end{remark}

\section{Perturbations of the lump solution}

In this section we consider perturbations of the lump solution 
(\ref{lump}). First we propagate initial data obtained from the lump 
 after multiplication with some scalar factor. Then we 
consider a perturbation  with a Gaussian and 
a deformed lump.

\subsection{Perturbation of the lump by rescaled initial data}
\begin{figure}[htb!]
\centering
\includegraphics[width=0.45\textwidth]{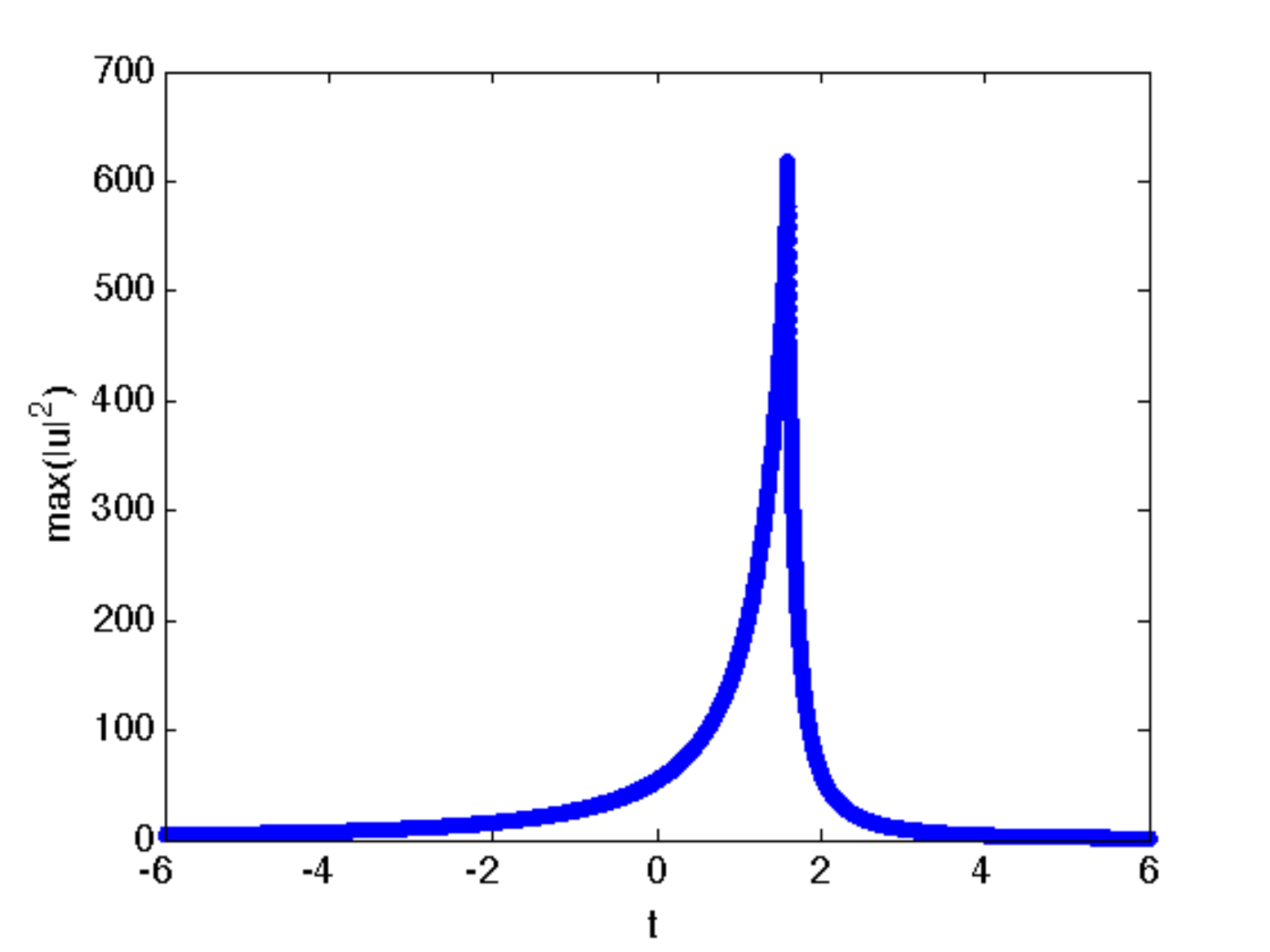}
\includegraphics[width=0.45\textwidth]{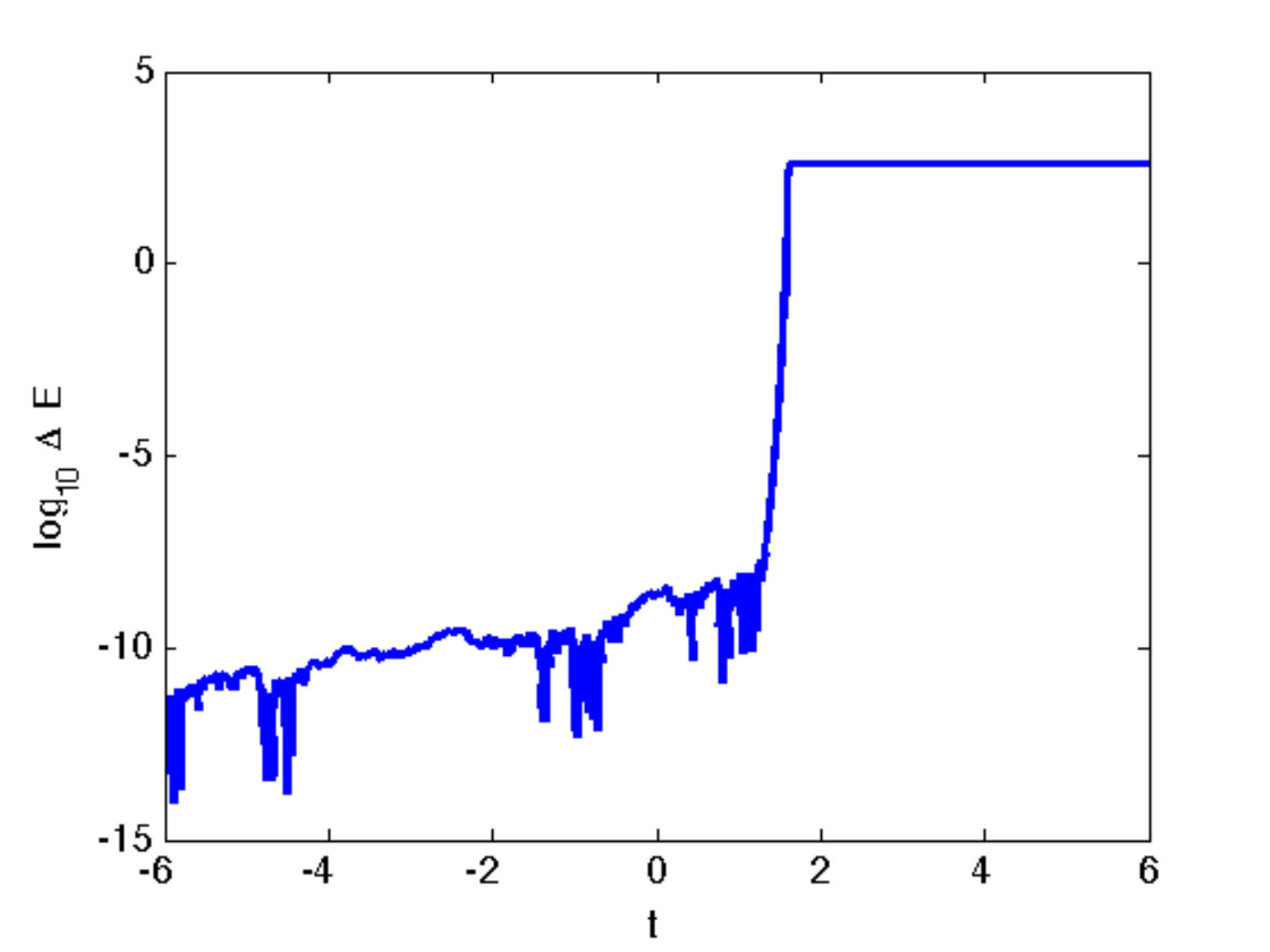}
\caption{Evolution of $max(|u|^{2})$ and the numerically computed 
energy in dependence of time for a solution to the 
focusing DS II equation (\ref{DSII})
for an initial condition of the form $u(x,y,-6)=1.1 u_{l}$.}
\label{amplul11}
\end{figure}

\begin{figure}[htb!]
\centering
\includegraphics[width=0.45\textwidth]{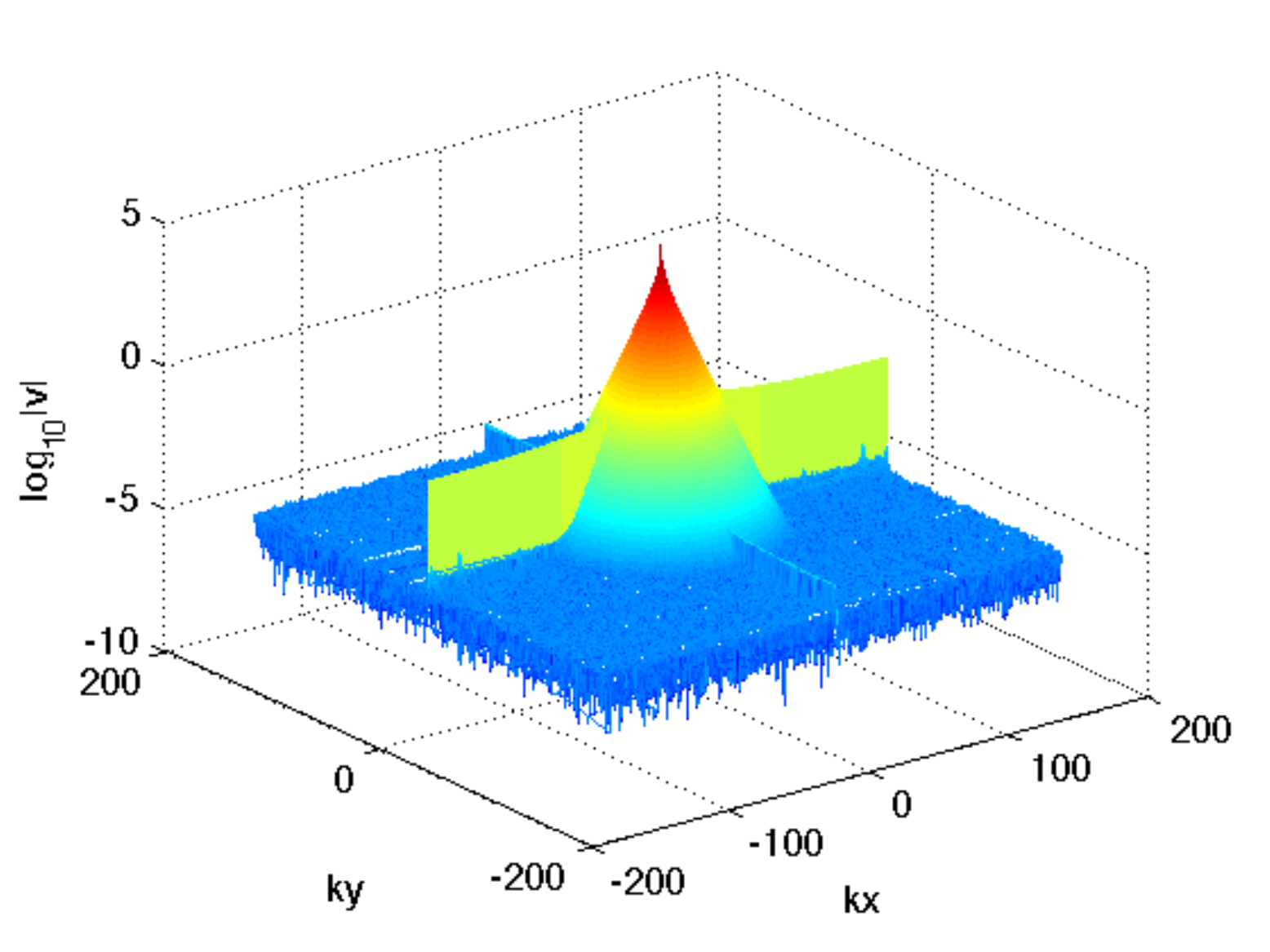}
\includegraphics[width=0.45\textwidth]{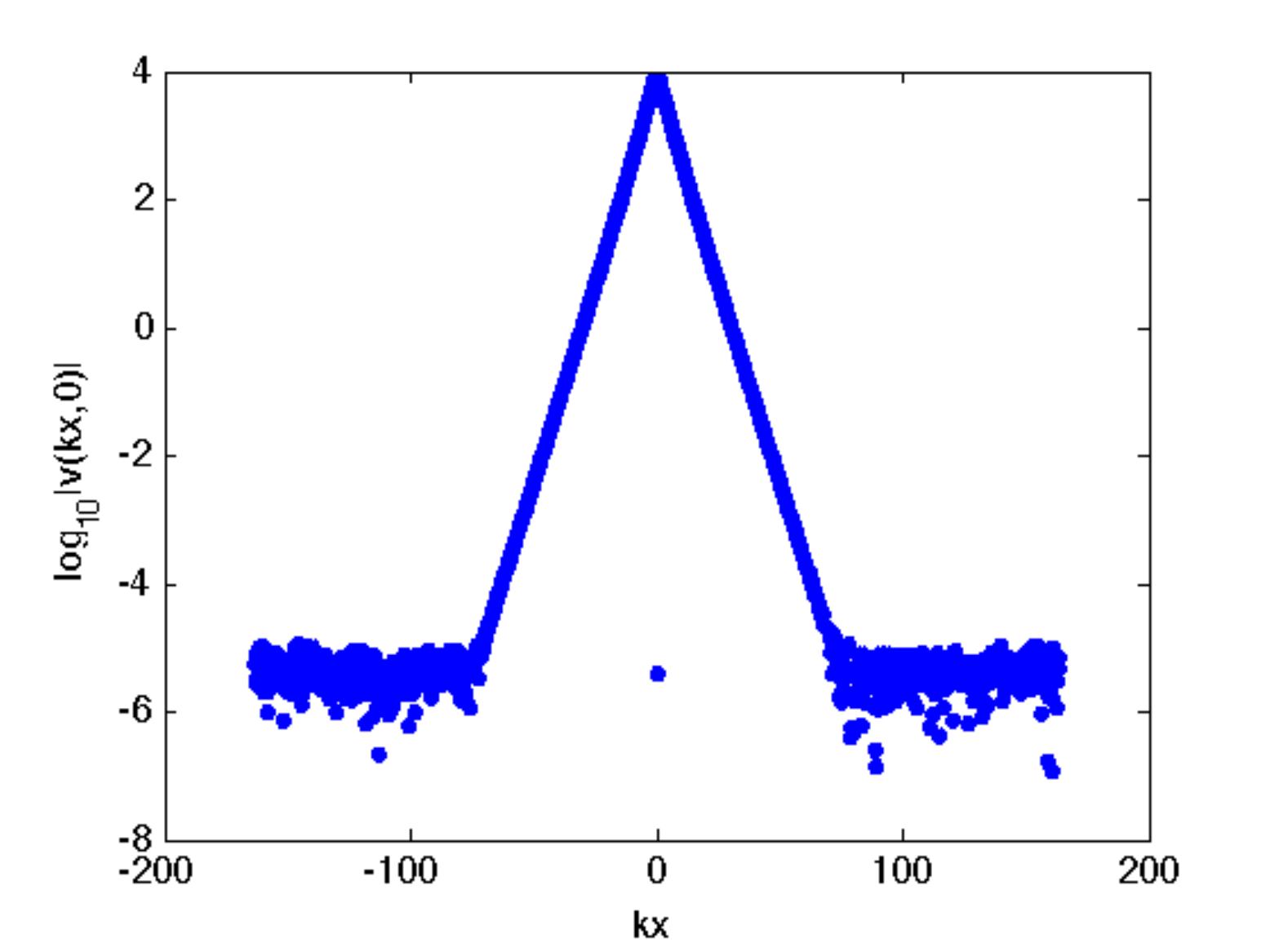}
\caption{Fourier coefficients  at $t=0$ for a solution to the 
focusing DS II equation (\ref{DSII})
for an initial condition of the form $u(x,y,-6)=1.1 u_{l}$.} 
\label{l11cf}
\end{figure}

We first consider  rescaled initial data from 
the lump (\ref{lump}) denoted by $u_{l}$
$$u(x,y,-6)=Au_{l},$$
where $A\in\mathbb{R}$ is a scaling factor.
The computations are carried out with $N_{x}=N_{y}=2^{14}$ points for 
$x\times y \in [-50\pi, 50\pi] \times [-50\pi, 50\pi]$ and $t\in[-6,6]$.
\\
\\
For $A=1.1$, and 
$N_t = 1000$, we observe a blowup of the solution at $t_c\sim1.6$.
The time evolution of $\underset{x,y}{\max} |u(x,y,t)|^2$ and of the energy 
is shown in Fig.~\ref{amplul11}.
The maximum of  $|u|^2$ in Fig.~\ref{amplul11} is 
clearly smaller than in the case of the Ozawa solution. This is due 
to the lower resolution in time which is used for this computation.
Nevertheless, 
the jump in the energy is obviously present.
The  Fourier coefficients at $t=0$ can be seen in Fig.~\ref{l11cf}. 
They again decrease by almost 6 orders of magnitude. 

\begin{figure}[htb!]
\centering
\includegraphics[width=0.45\textwidth]{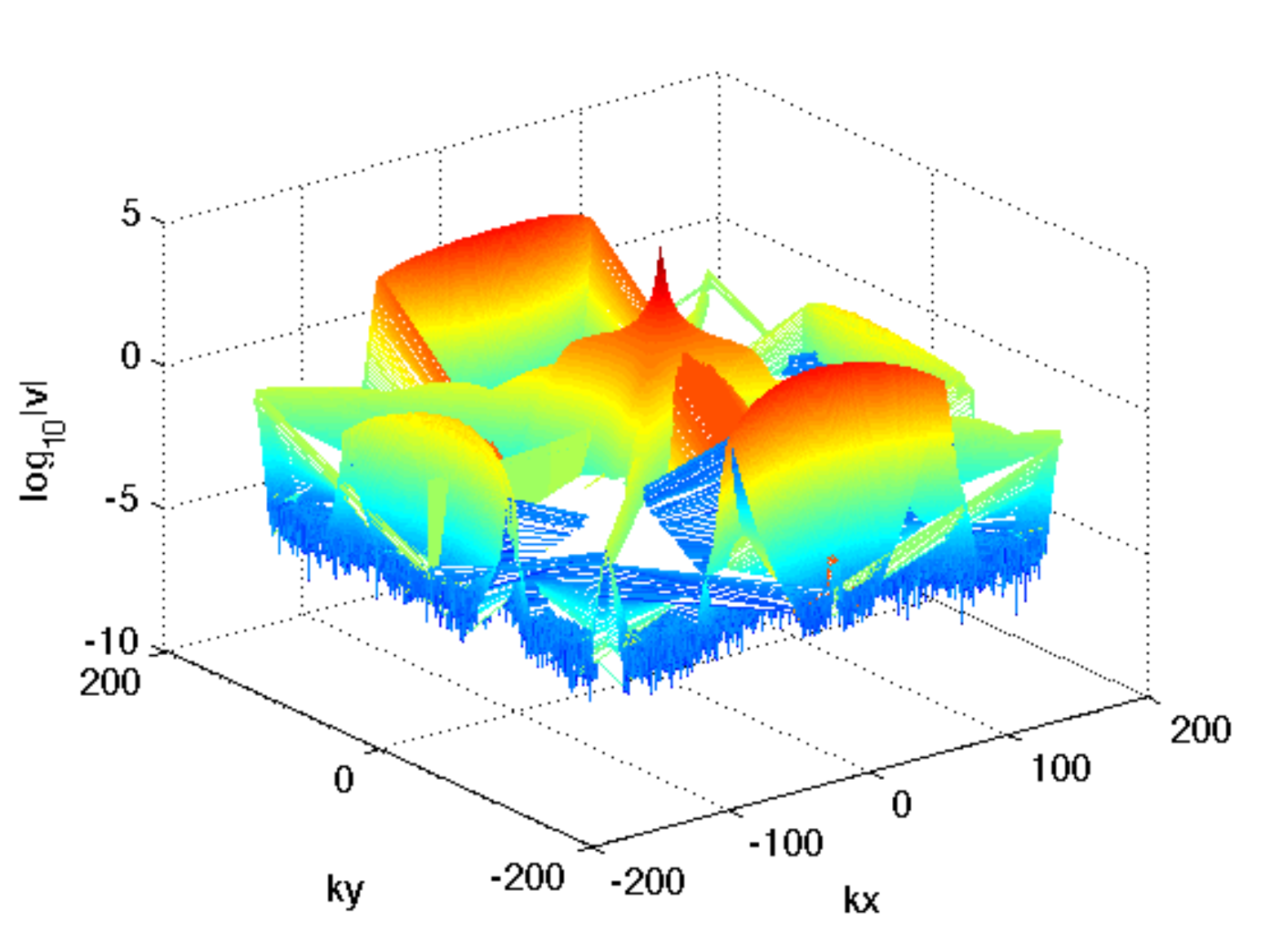}
\includegraphics[width=0.45\textwidth]{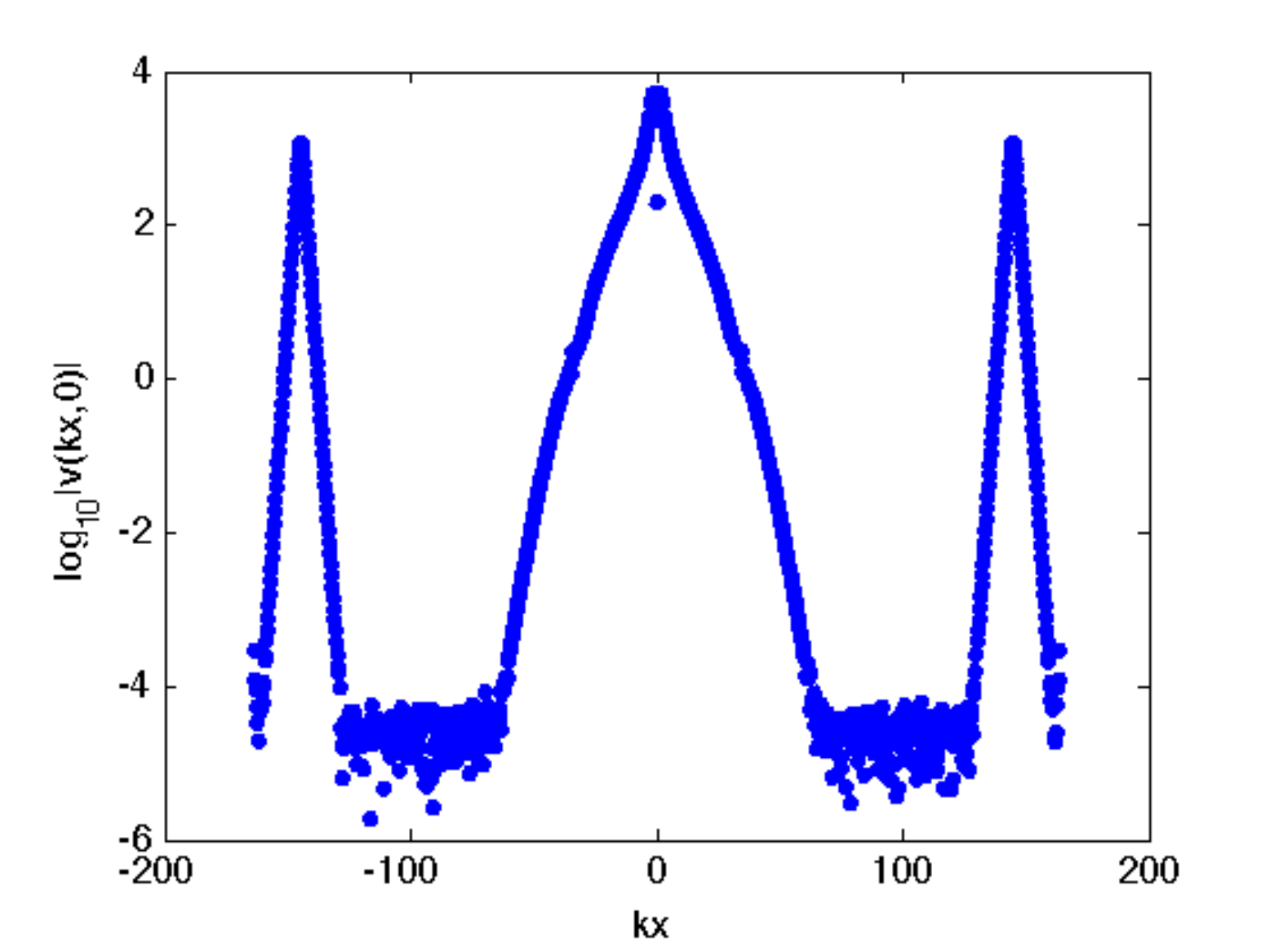}
\caption{Fourier coefficients  at $t=6$ for a solution to the 
focusing DS II equation (\ref{DSII})
for an initial condition of the form $u(x,y,-6)=1.1 u_{l}$.} 
\label{l11cff}
\end{figure}

\begin{figure}[htb!]
\centering
\includegraphics[width=0.45\textwidth]{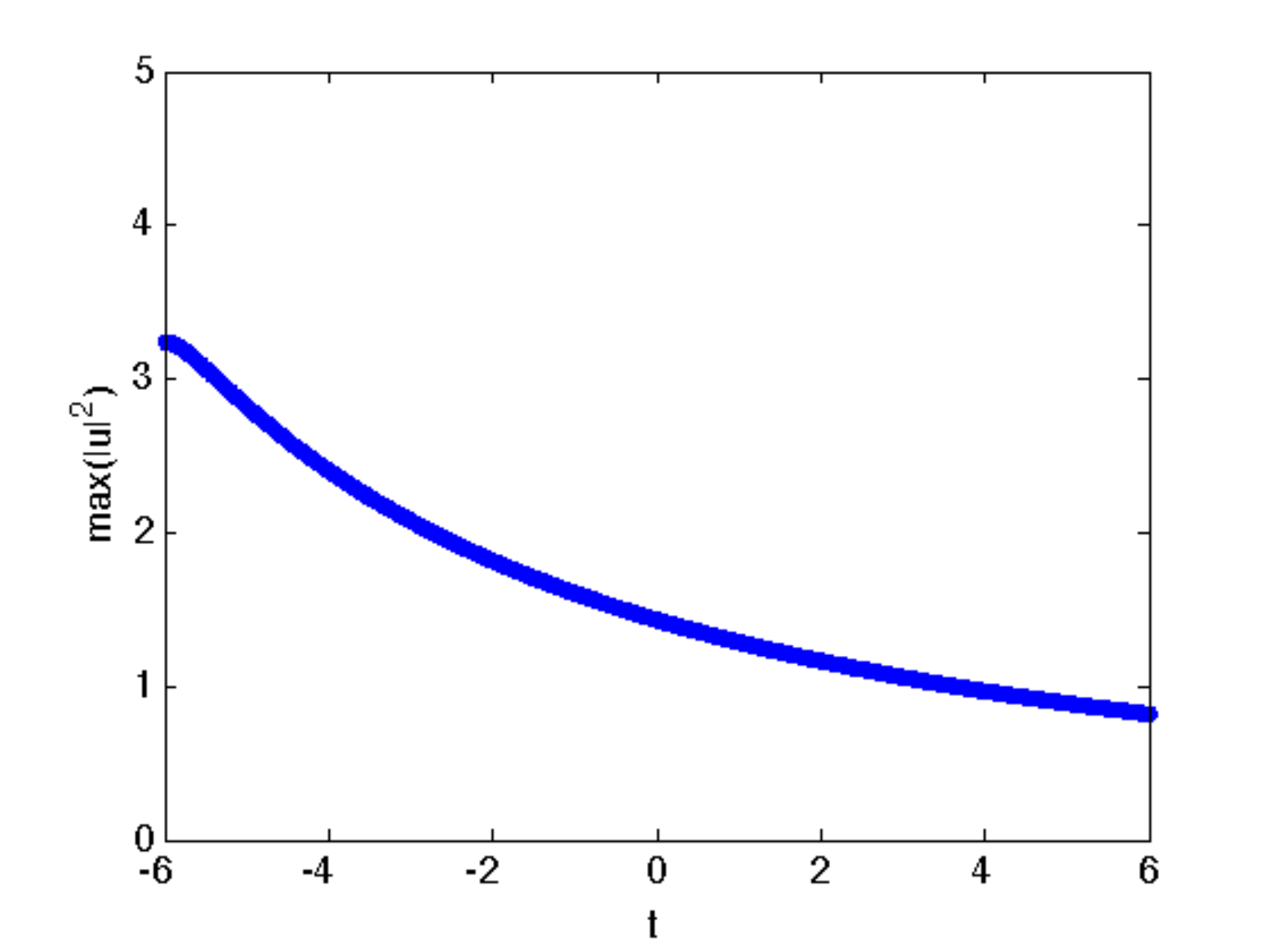}
\includegraphics[width=0.45\textwidth]{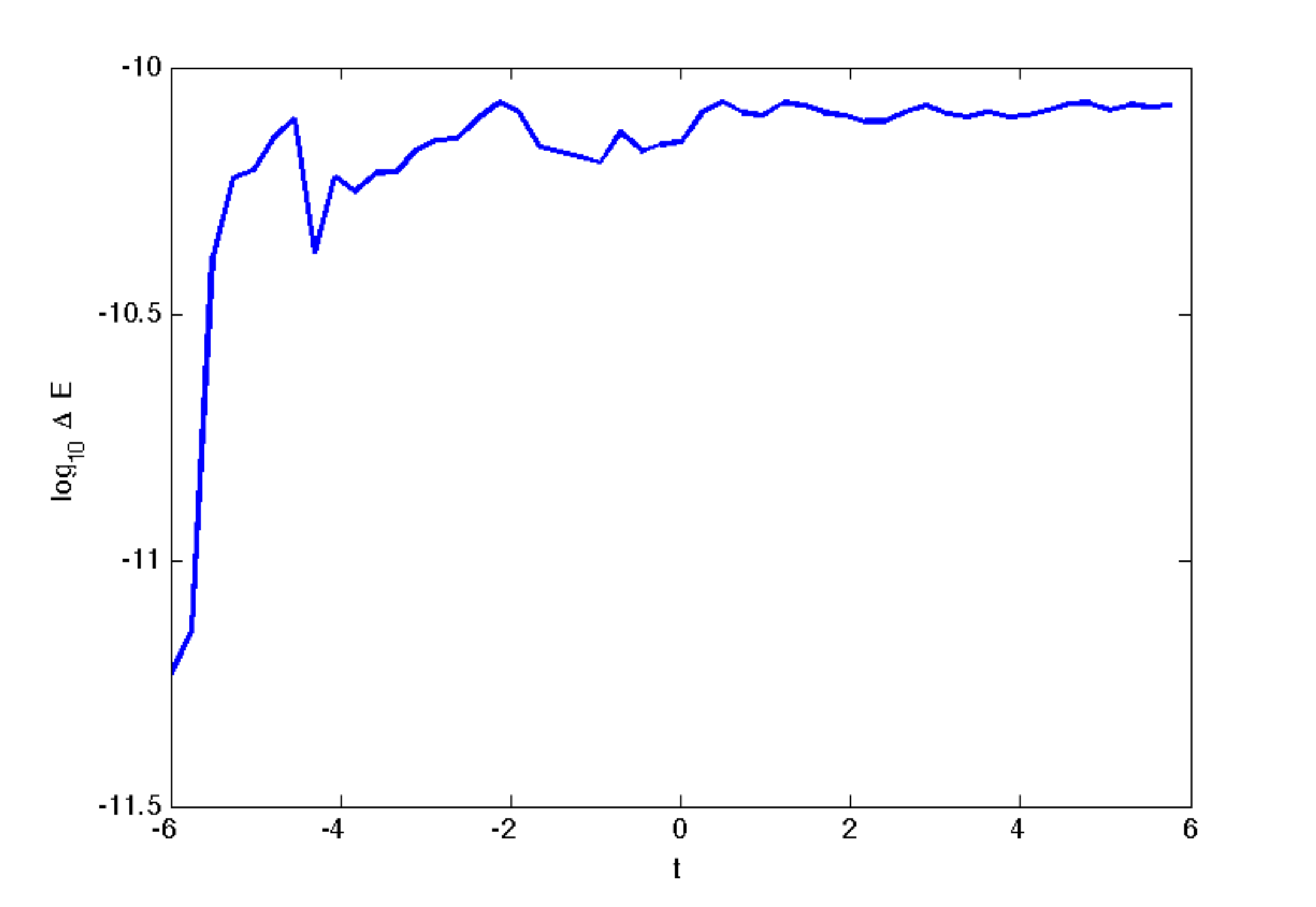}
\caption{Evolution of $max(|u|^{2})$ and the numerically computed energy 
in dependence of time  for a solution to the 
focusing DS II equation (\ref{DSII}) for an initial condition 
of the form $u(x,y,-6)=0.9 u_l$.}
\label{ampllump09}
\end{figure}

\begin{figure}[htb!]
\centering
\includegraphics[width=0.45\textwidth]{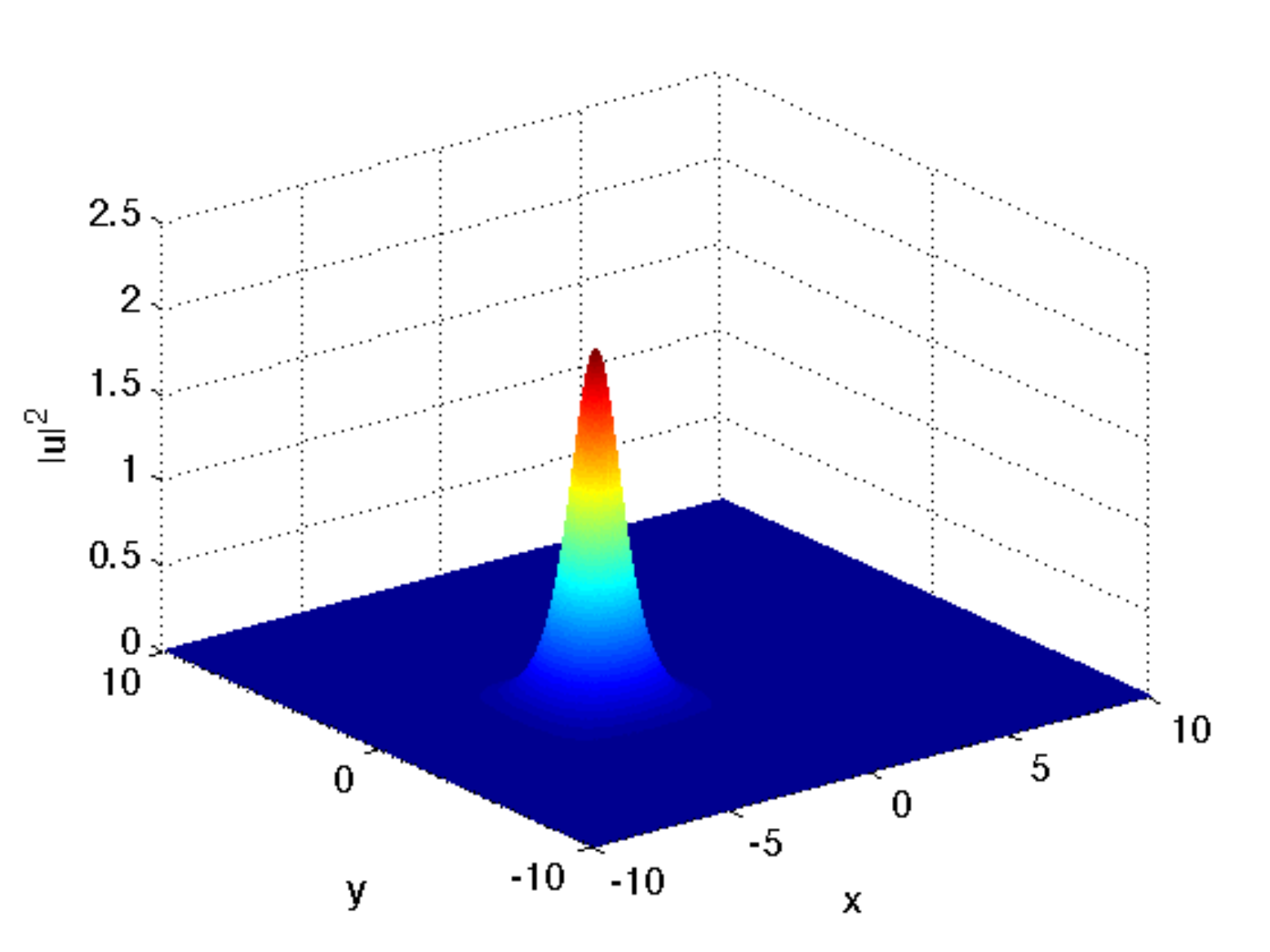}
\includegraphics[width=0.45\textwidth]{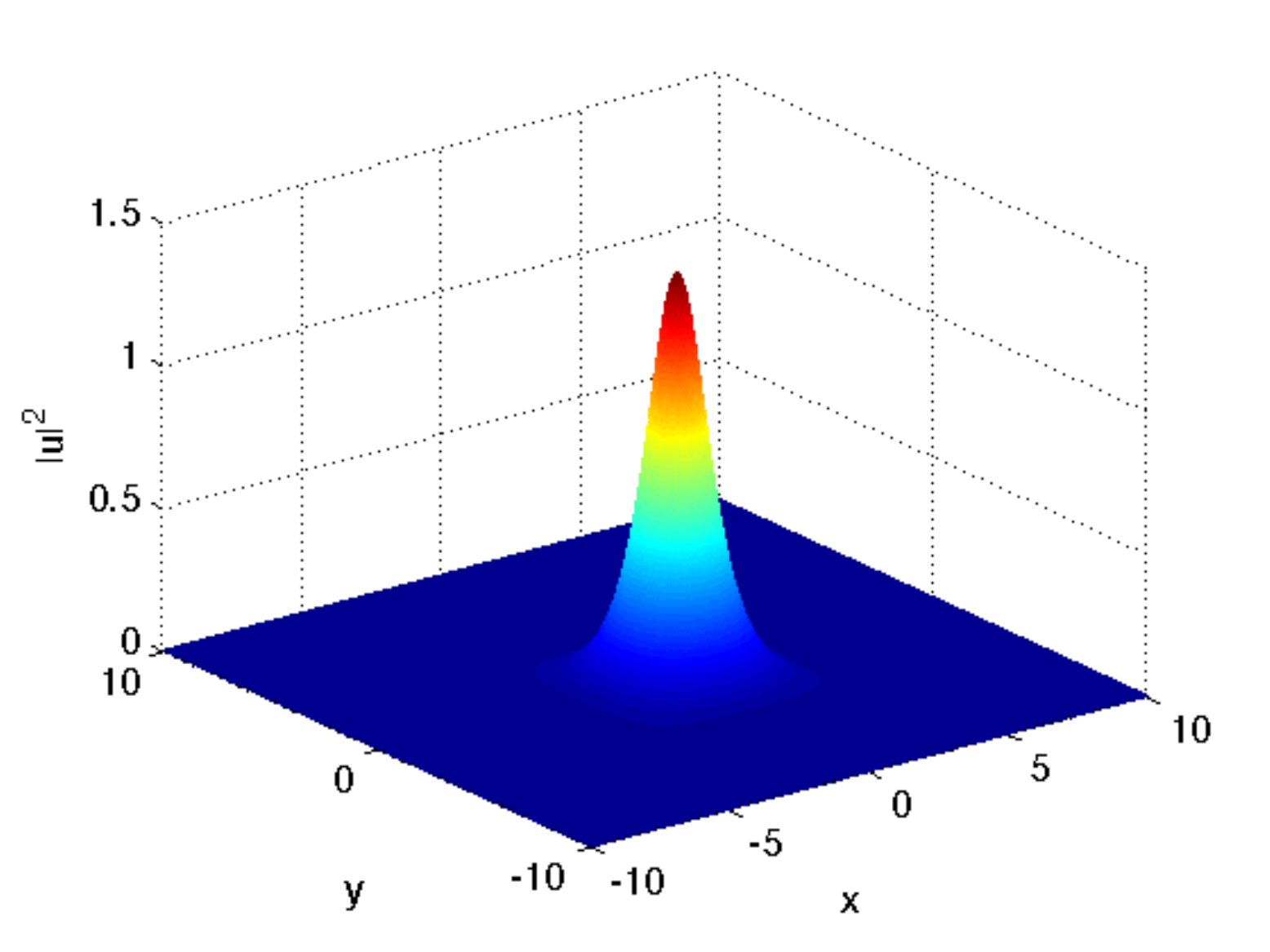}
\includegraphics[width=0.45\textwidth]{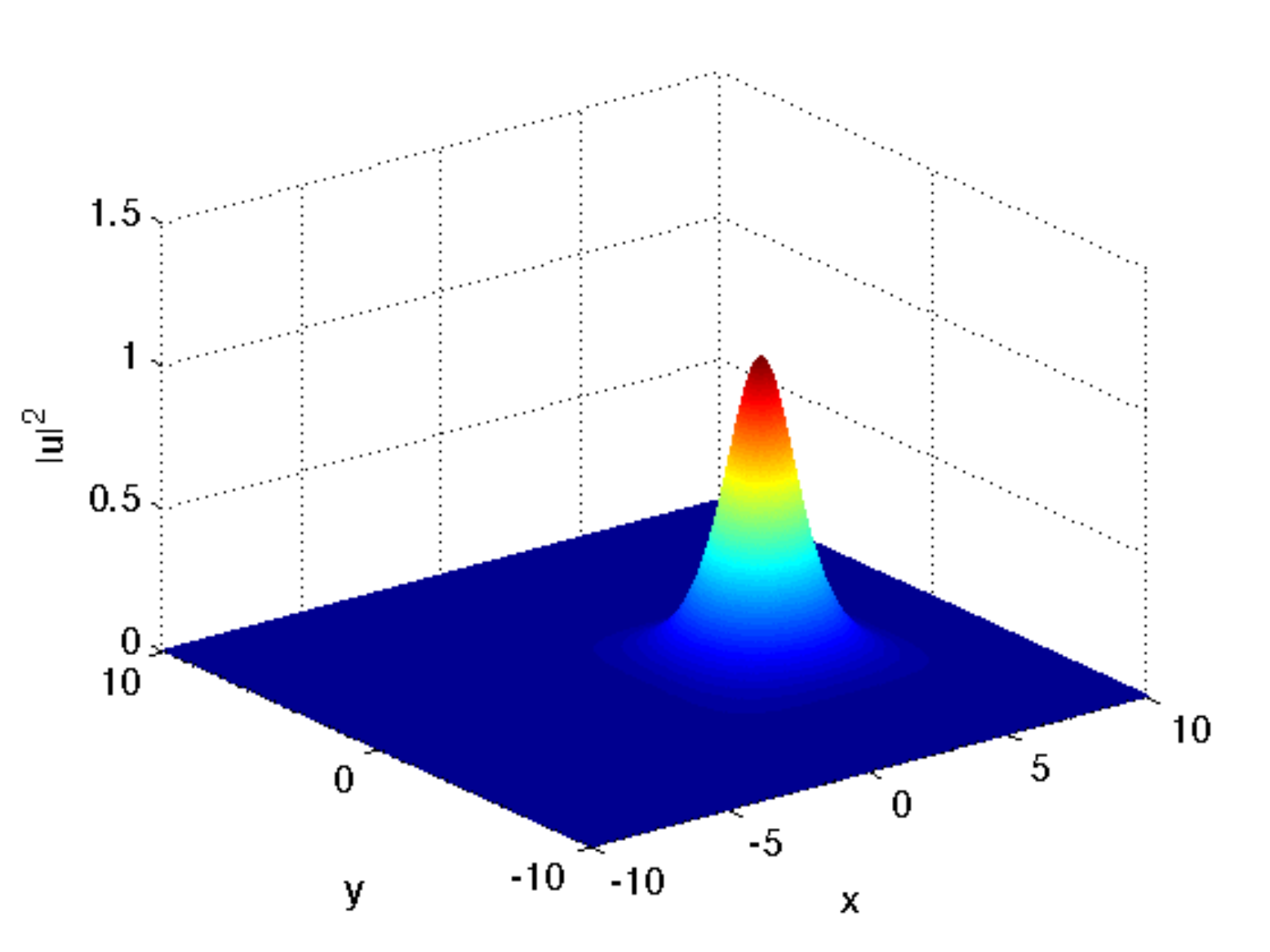}
\includegraphics[width=0.45\textwidth]{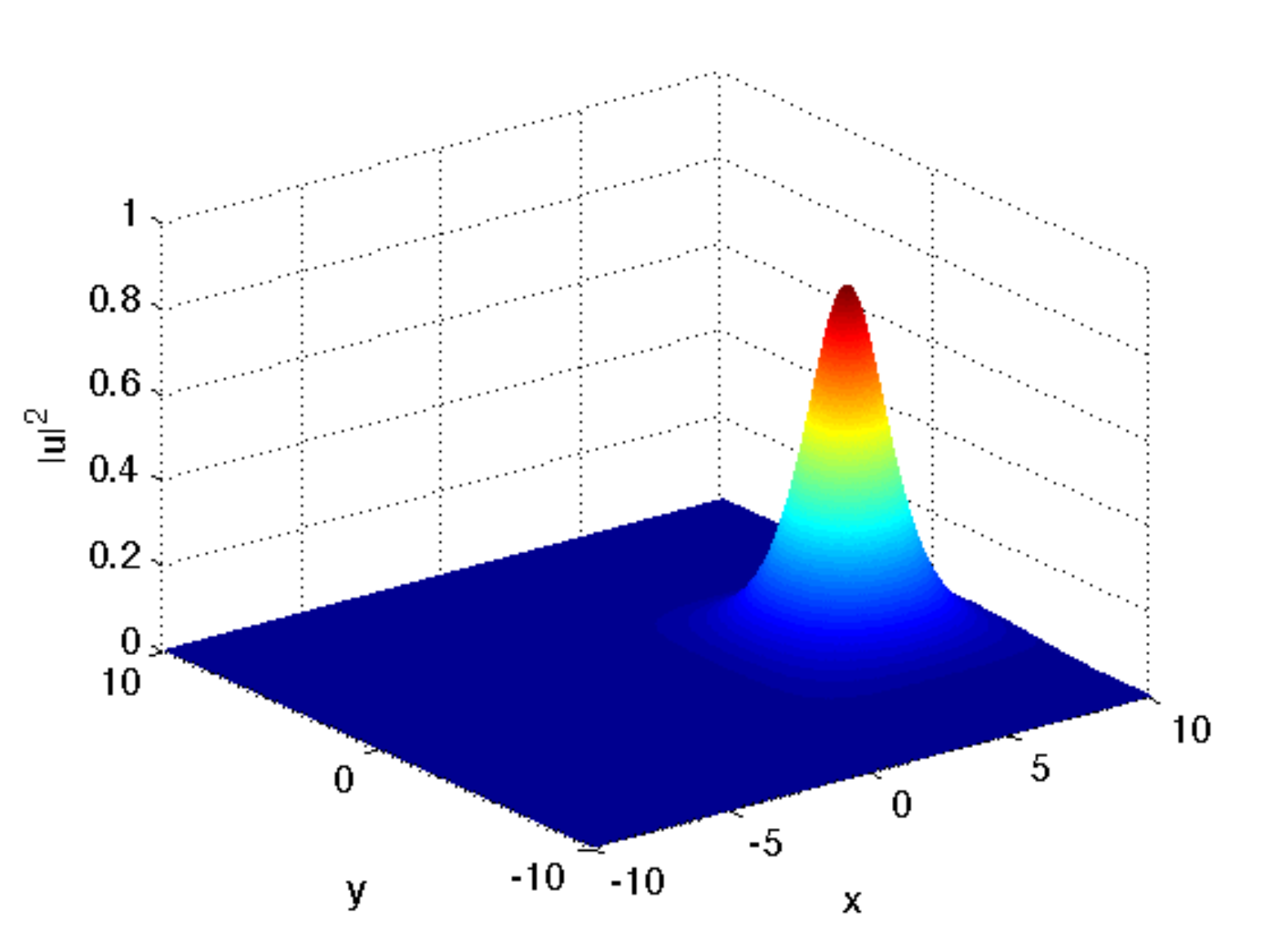}
\caption{Solution to the focusing DS II equation (\ref{DSII}) for an initial condition 
of the form $u(x,y,-6)=0.9 u_l$ for 
$t=-3$ and $t=0$ in the first row and $t=3$ and $t=6$ 
below.}
\label{lump09}
\end{figure}

\begin{figure}[htb!]
\centering
\includegraphics[width=0.45\textwidth]{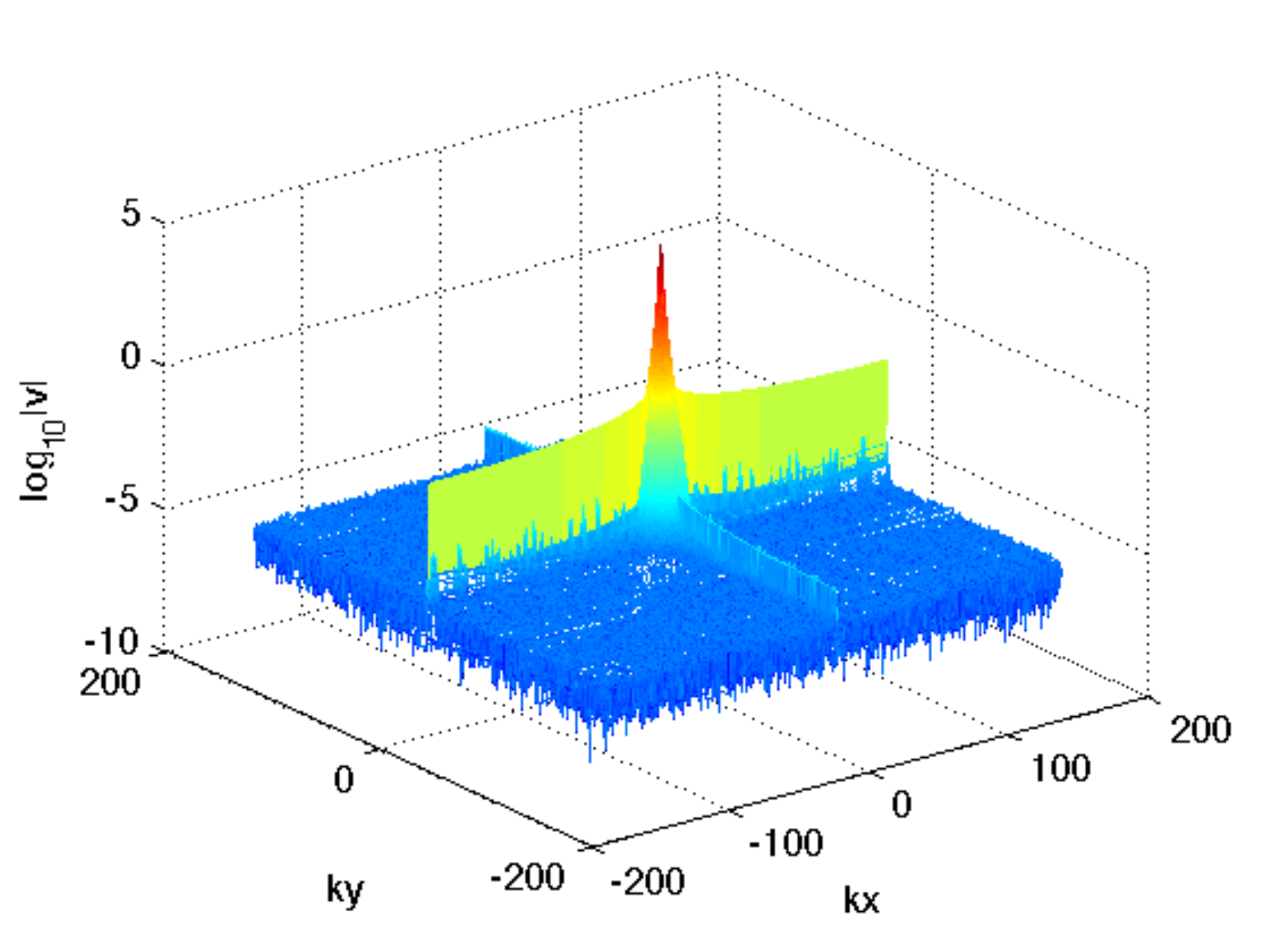}
\includegraphics[width=0.45\textwidth]{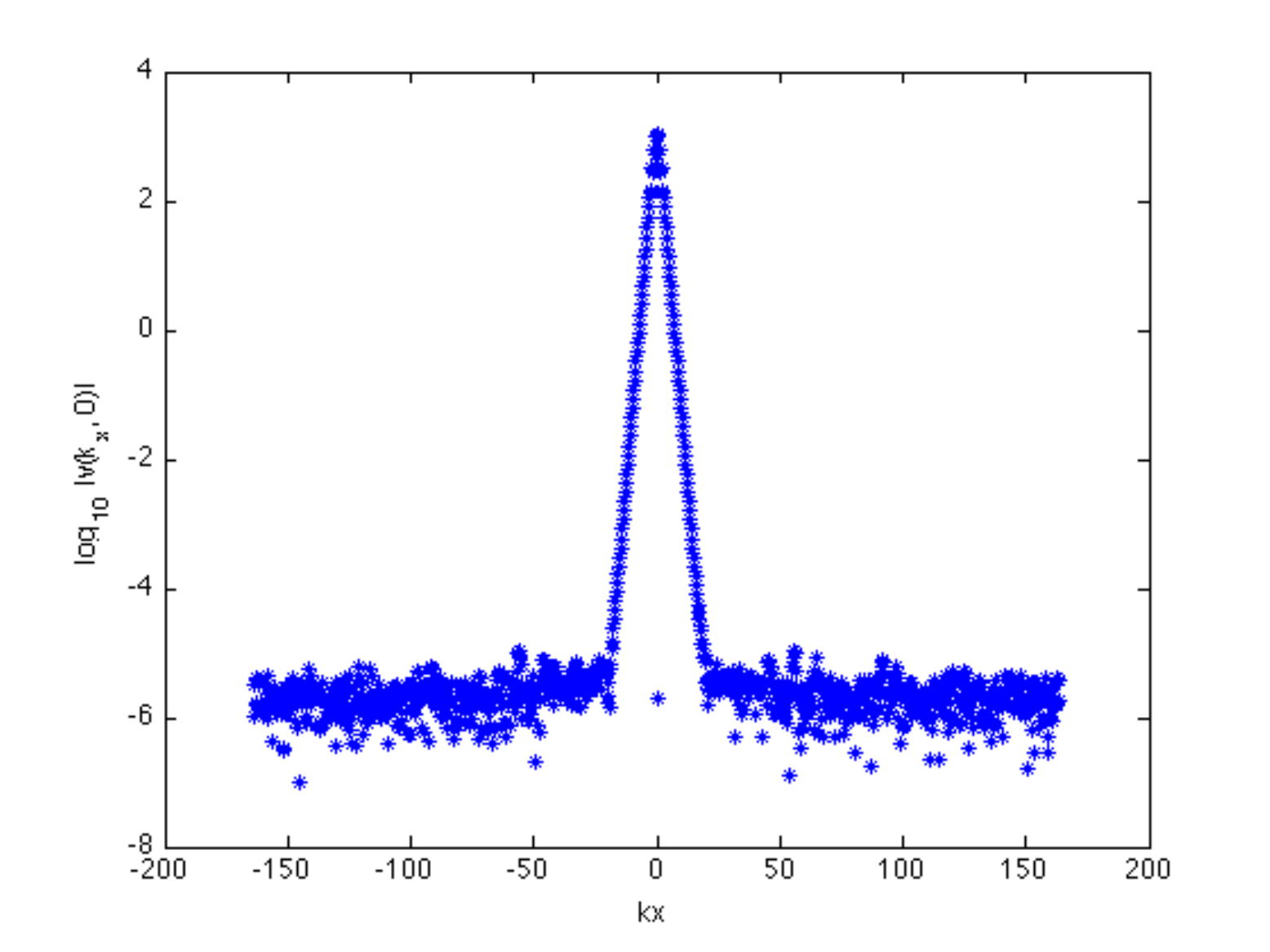}
\caption{The Fourier coefficients at $t=0$ of the solution to the 
focusing DS II equation (\ref{DSII})
for an initial condition of the form $u(x,y,-6)=0.9 u_l$.} 
\label{l09cf}
\end{figure}

To illustrate  the modulational instability at a concrete example, we 
show the Fourier coefficients after the critical time in 
Fig.~\ref{l11cff}. It can be seen that the modulus 
of the coefficients of the high wavenumbers increases instead of decreasing as to be expected for smooth functions. 
This indicates once more that the computed solution after the blowup 
time has to be taken with a grain of salt. 
\\
\\
For $A=0.9$, the initial pulse travels in the same direction as 
the exact solution, but loses speed and height and is broadened, see 
Fig.~\ref{ampllump09}. It appears that this modified lump just 
disperses asymptotically.
The solution can be seen in 
Fig.~\ref{lump09}. Its Fourier coefficients  in 
Fig.~\ref{l09cf} show that the resolution of the initial data is 
almost maintained. 

\subsection{Perturbation of the lump with a Gaussian}

\begin{figure}[htb!]
\centering
\includegraphics[width=0.45\textwidth]{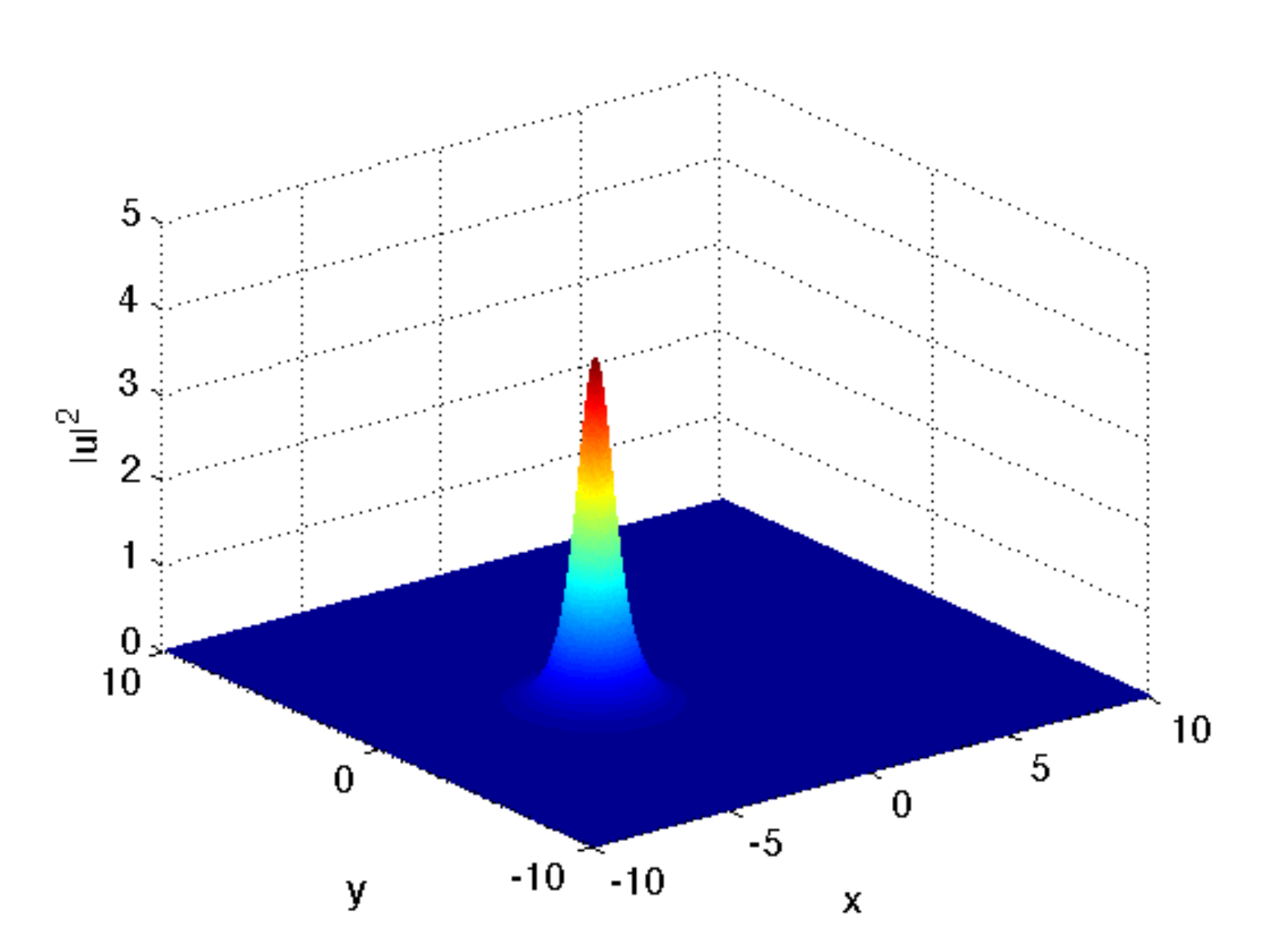}
\includegraphics[width=0.45\textwidth]{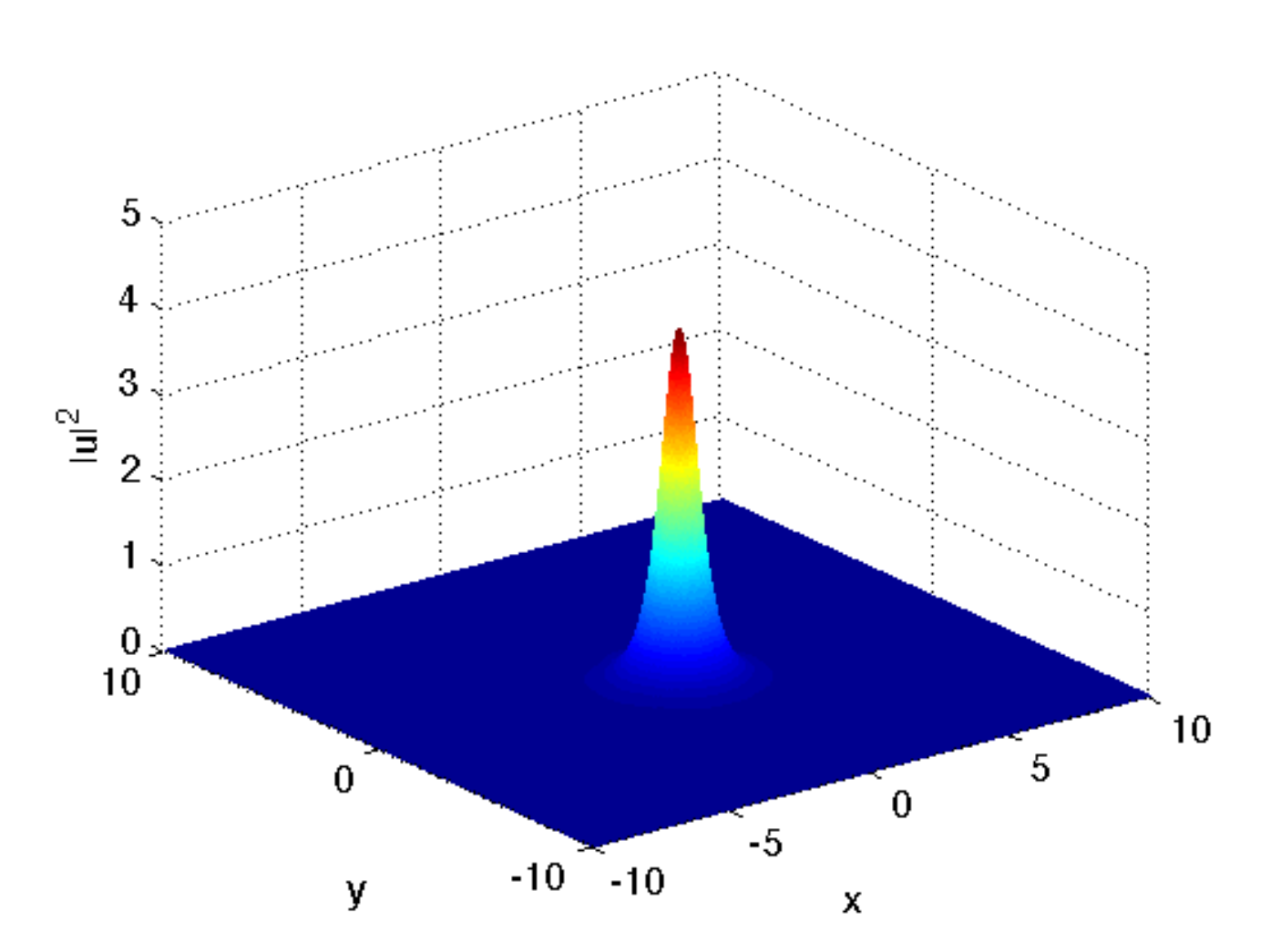}
\includegraphics[width=0.45\textwidth]{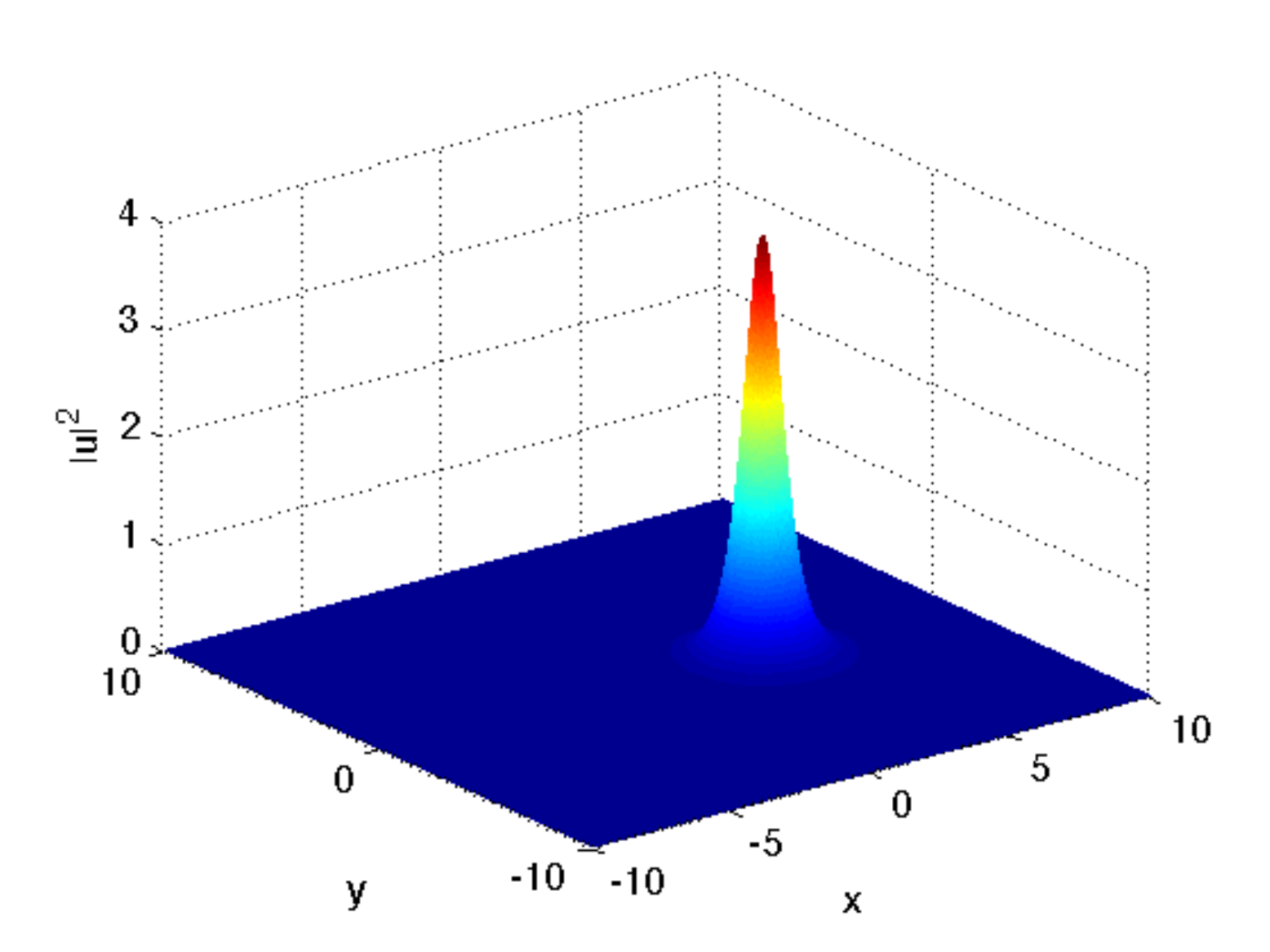}
\includegraphics[width=0.45\textwidth]{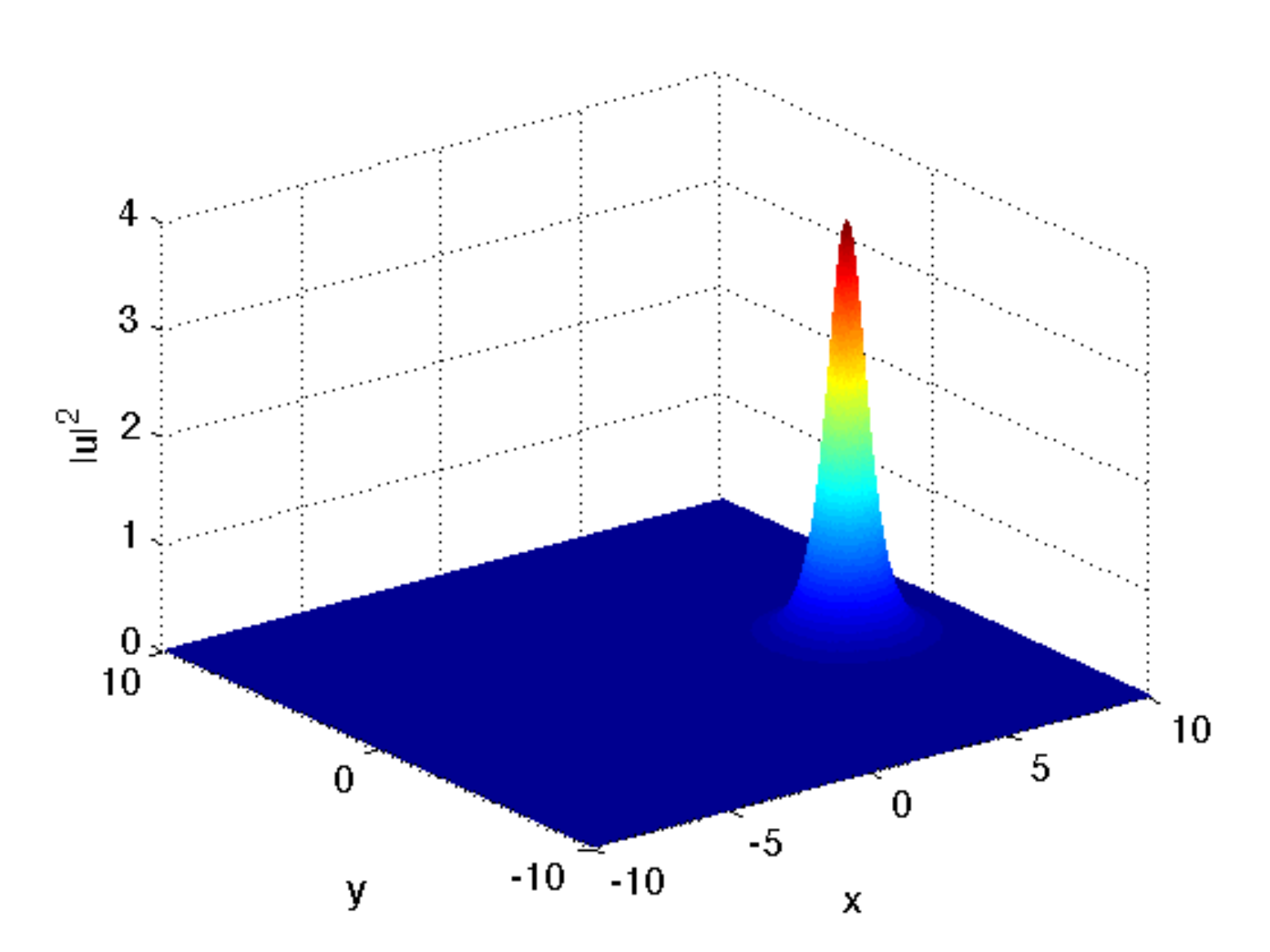}
\caption{Solution to the focusing DS II equation (\ref{DSII}) 
for an initial condition of the form (\ref{lumpgauss}) with $B=0.1$ for 
$t=-3$ and $t=0$ in the first row and $t=3$ and $t=6$ 
below.} 
\label{ulupg1}
\end{figure}

\begin{figure}[htb!]
\centering
\includegraphics[width=0.45\textwidth]{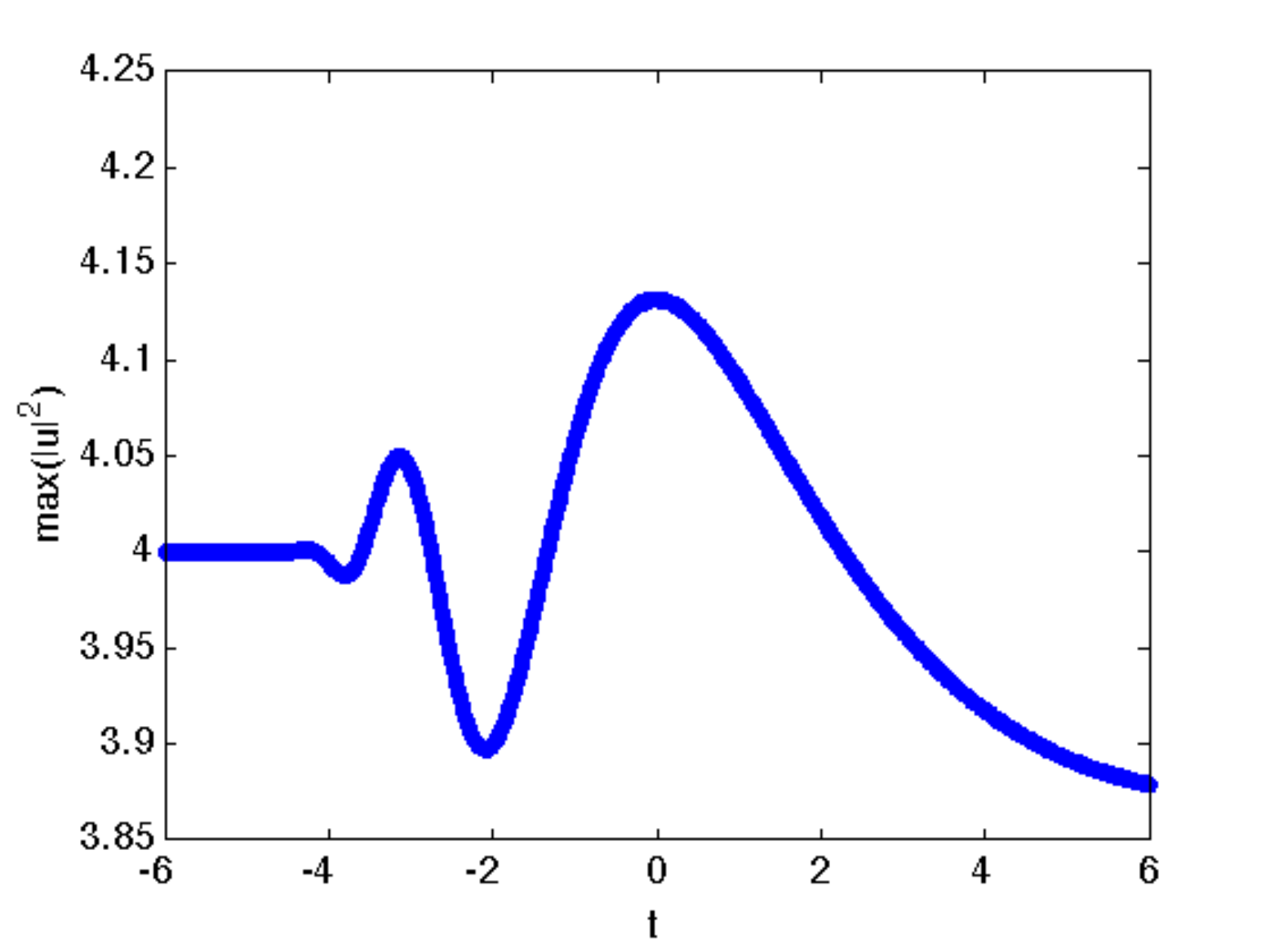}
\includegraphics[width=0.45\textwidth]{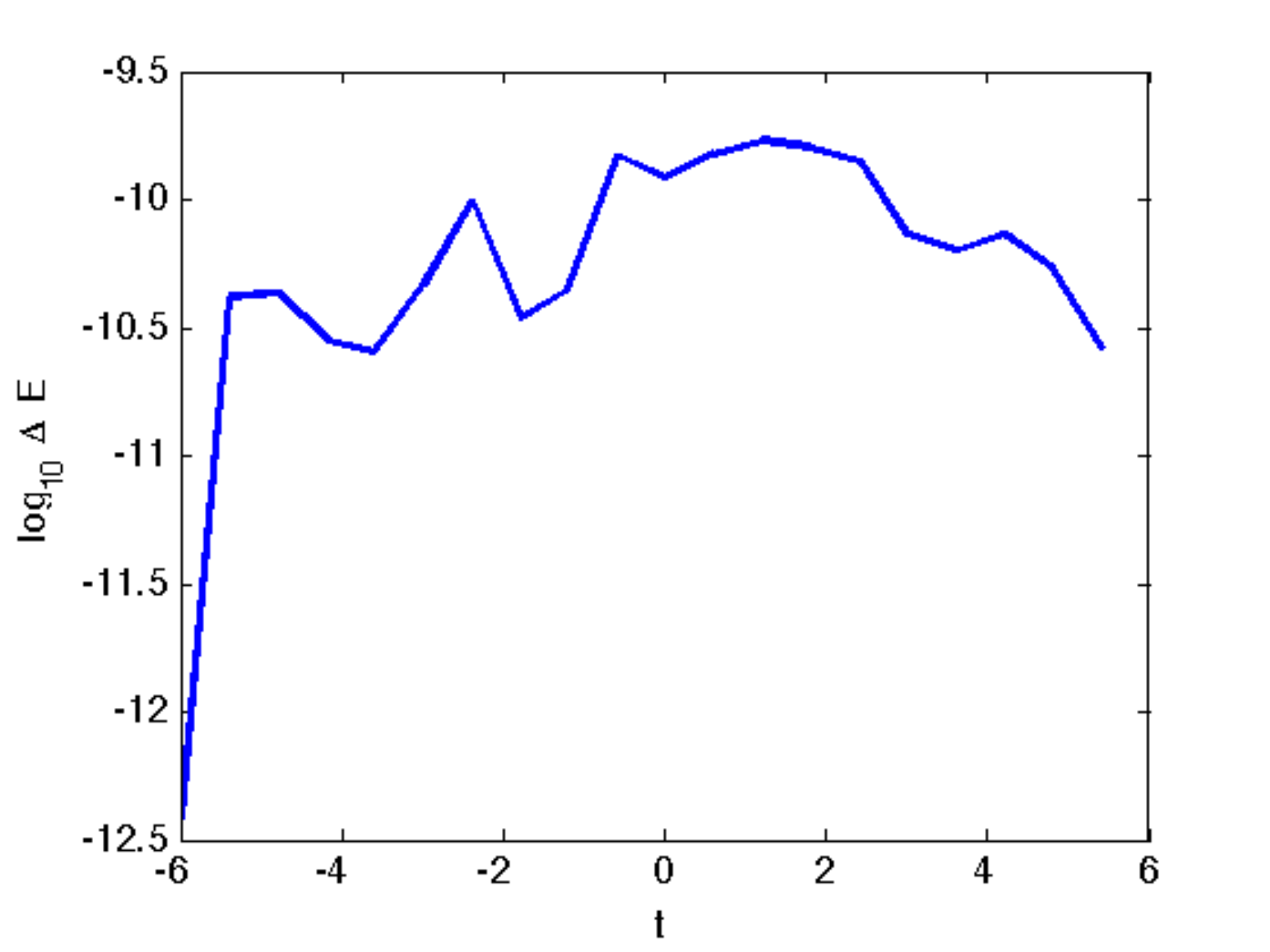}
\caption{Evolution of $max(|u|^{2})$ and of the energy in dependence of time 
for an initial condition of the form (\ref{lumpgauss}) with $B=0.1$.}
\label{amplulupg1}
\end{figure}

\begin{figure}[htb!]
\centering
\includegraphics[width=0.45\textwidth]{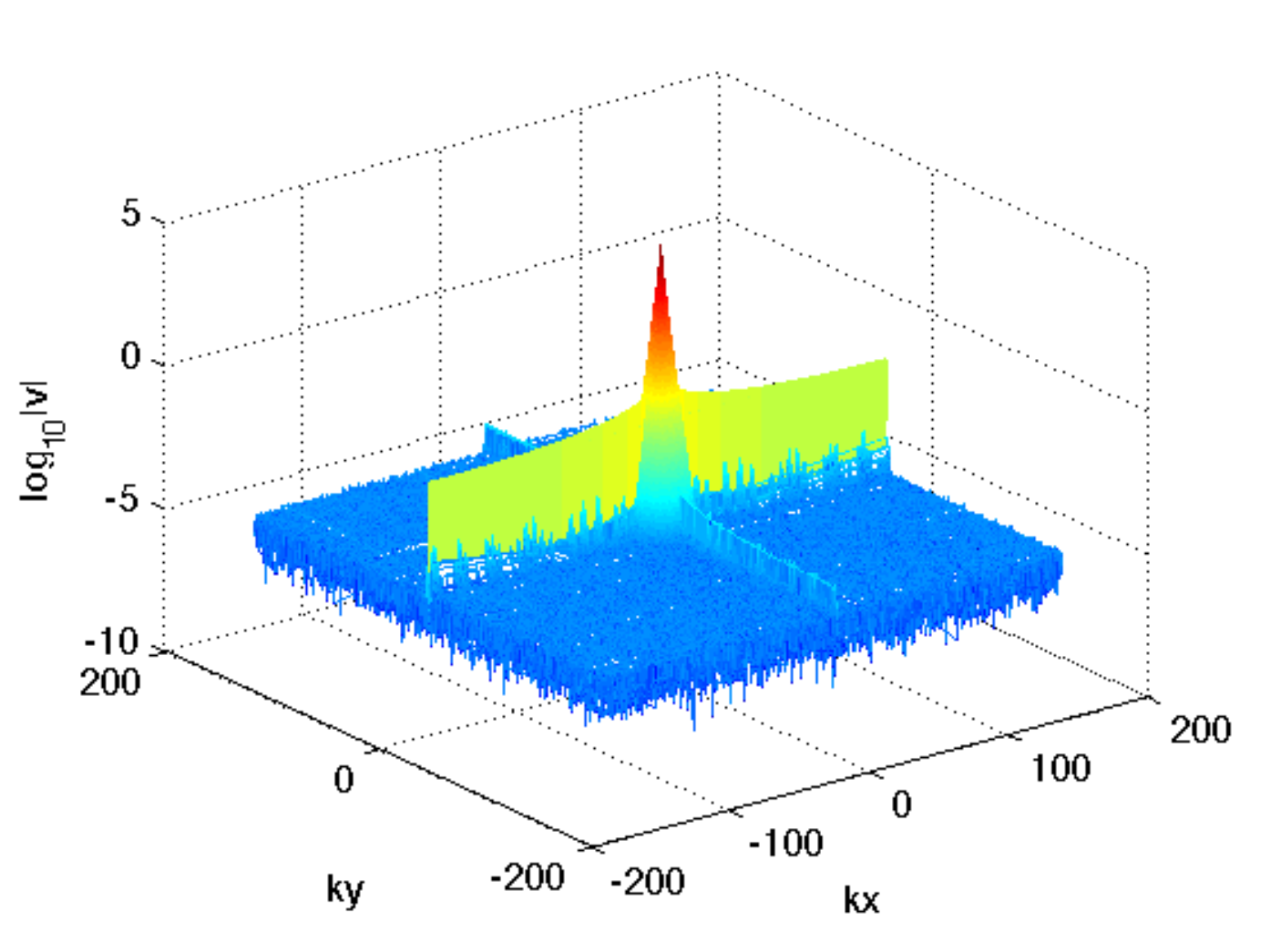}
\includegraphics[width=0.45\textwidth]{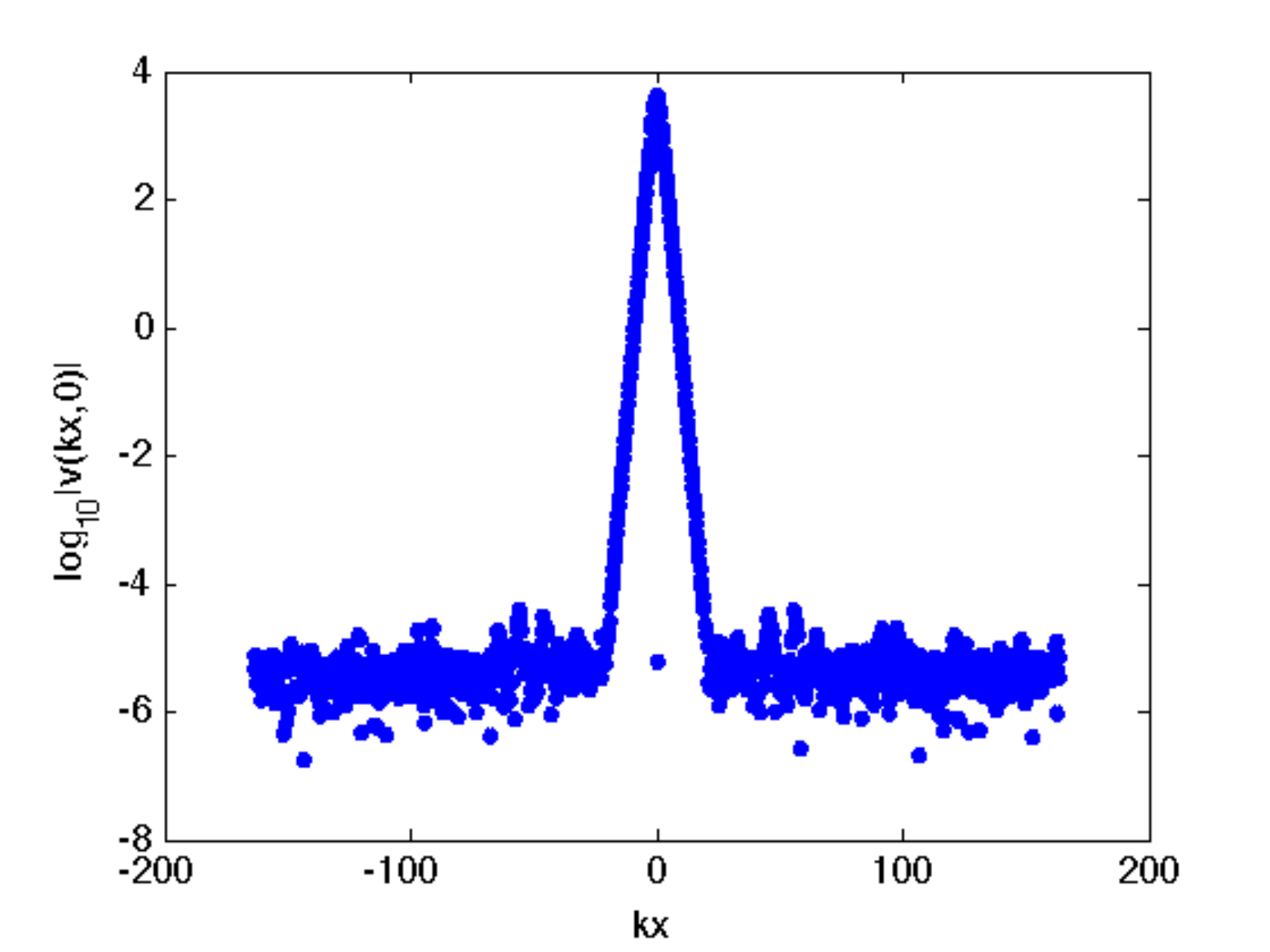}
\caption{Fourier coefficients of $u$ at $t=6$
for an initial condition of the form (\ref{lumpgauss}) with $B=0.1$.} 
\label{lu01cf}
\end{figure}

We consider an initial condition of the form 
\begin{equation}
    u(x,y,-6) = u_l + B \exp(-(x^2+y^2)), \,\, B \in \mathbb{R}.
    \label{lumpgauss}
\end{equation}
For $B=0.1$ and
$N_t = 1000$, 
we show the solution at different times in Fig.~\ref{ulupg1}.
The solution travels at the same speed as before, but its amplitude varies, 
growing and decreasing successively, see Fig.~\ref{amplulupg1}.
The time evolution of 
 the energy can be seen in Fig.~\ref{amplulupg1}. There is no 
 indication of blowup in this example. The solution appears to 
 disperse for $t\to\infty$.
The Fourier coefficients at $t=6$ in Fig.~\ref{lu01cf} show the 
wanted spatial resolution. 
\\
\\
A similar behavior is observed if a larger value for the amplitude of 
the perturbation is chosen, e.g., $B=0.5$.

\subsection{Deformation of the Lump}

\begin{figure}[htb!]
\centering
\includegraphics[width=0.5\textwidth]{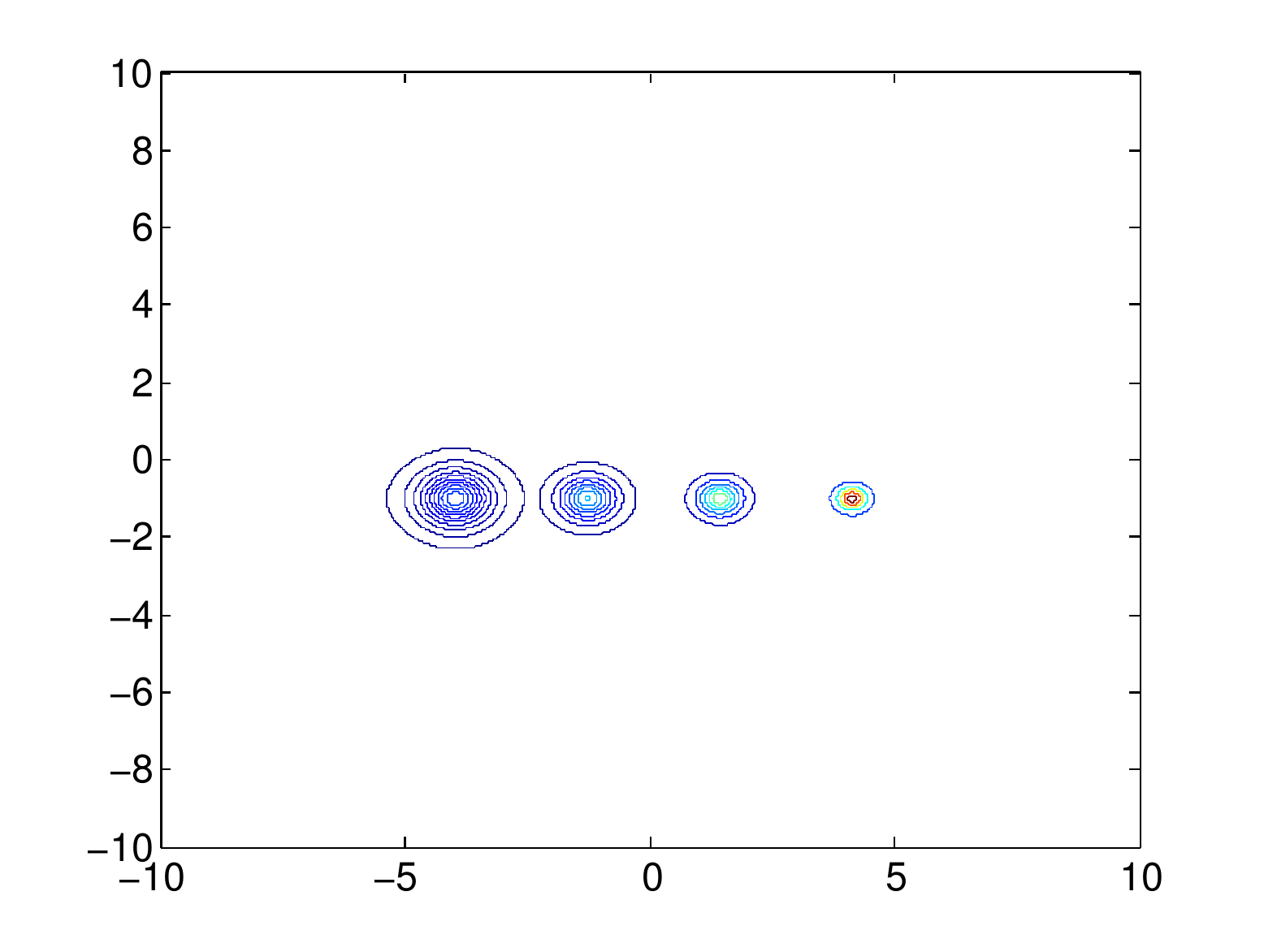}
\caption{Contour plot for a solution to the focusing DS II equation 
(\ref{DSII}) 
for an initial condition of the form (\ref{lumpyd}) with $\kappa=0.9$ 
for different times.}
\label{conluy09}
\end{figure}

\begin{figure}[htb!]
\centering
\includegraphics[width=0.45\textwidth]{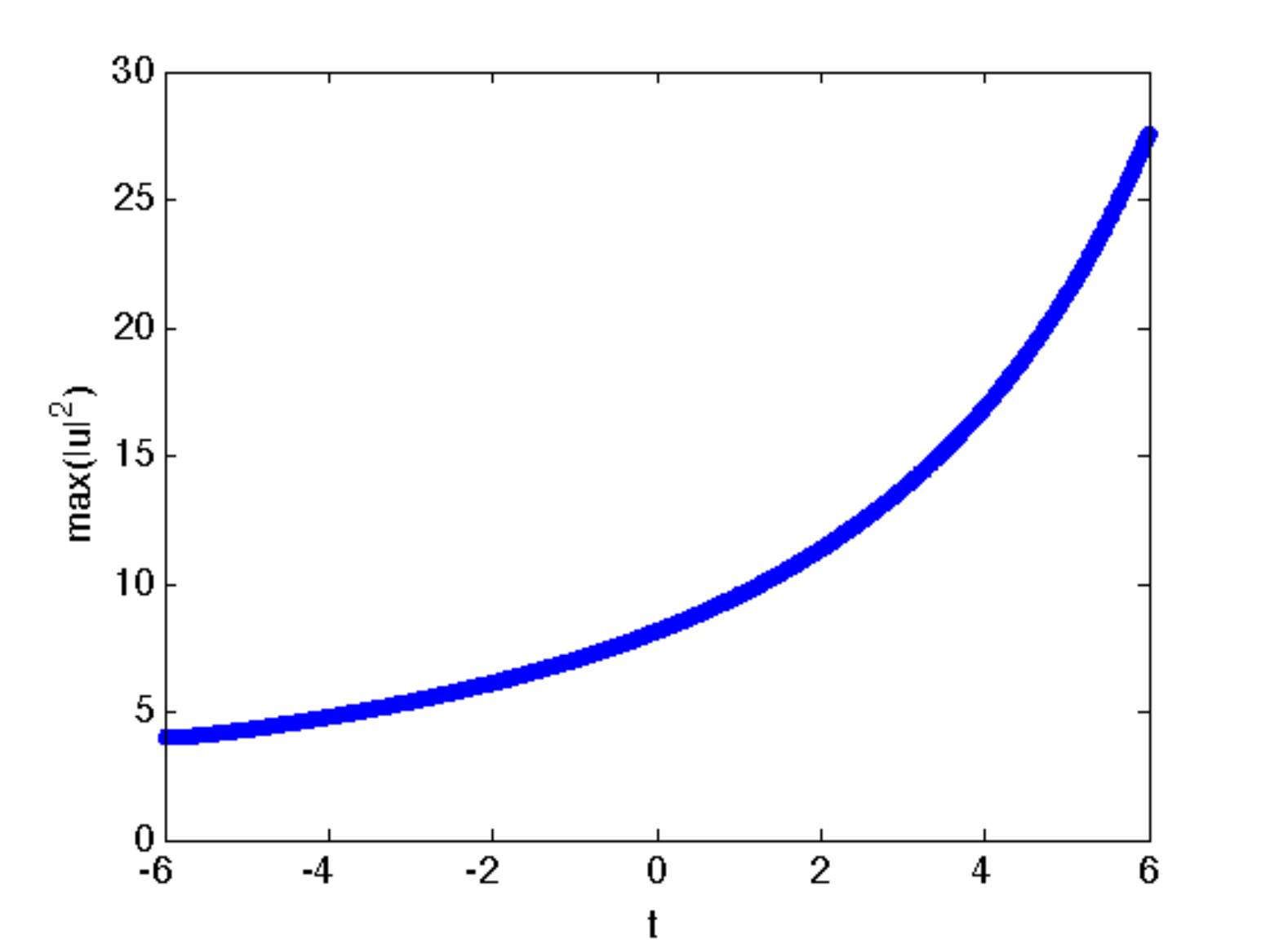}
\includegraphics[width=0.45\textwidth]{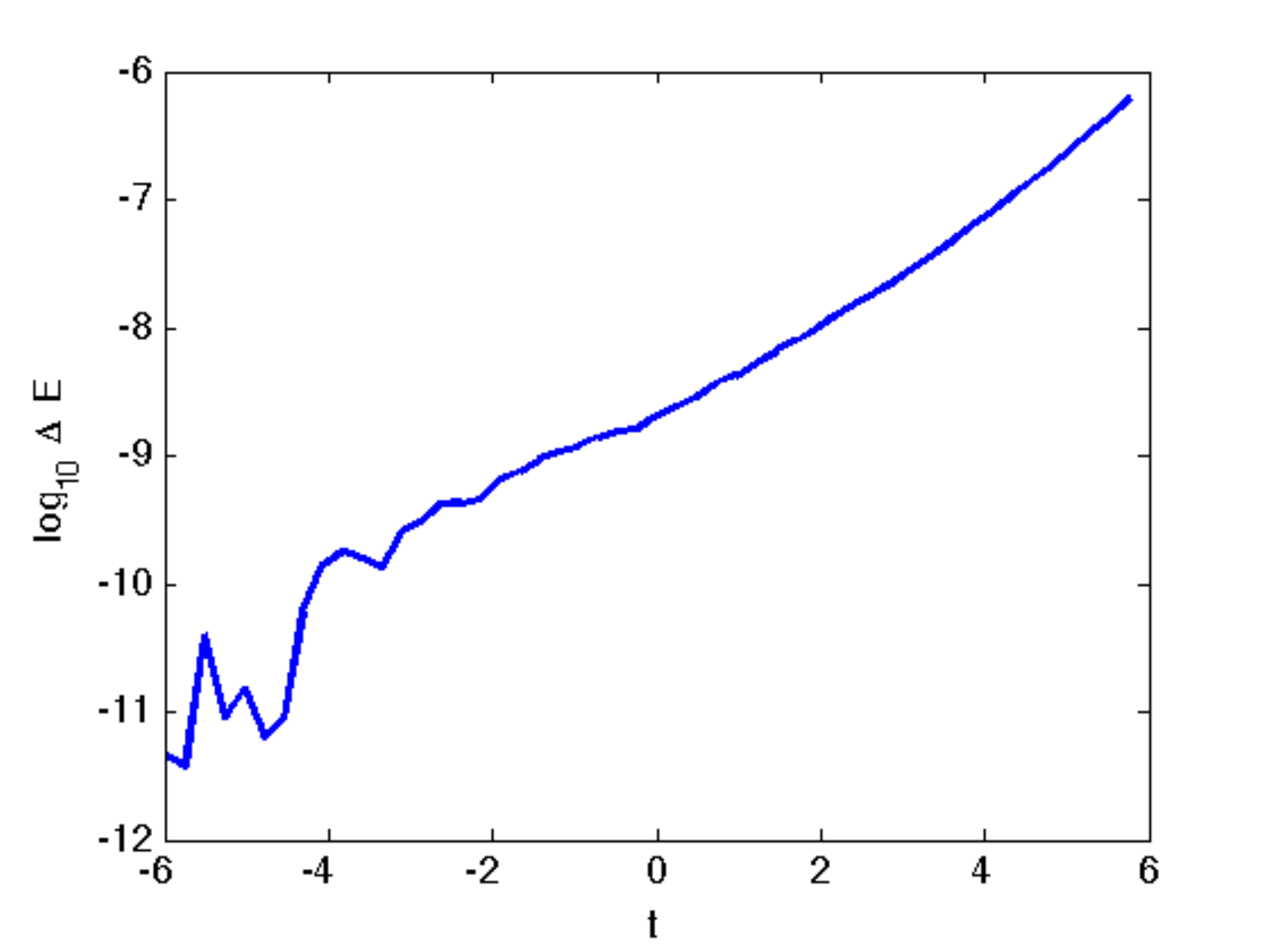}
\caption{Evolution of  $max(|u|^{2})$ and the numerically computed 
energy in dependence of time for the 
focusing DS II equation 
(\ref{DSII}) for an initial condition of the form (\ref{lumpyd}) with 
$\kappa=0.9$.}
\label{ampluy09}
\end{figure}

\begin{figure}[htb!]
\centering
\includegraphics[width=0.45\textwidth]{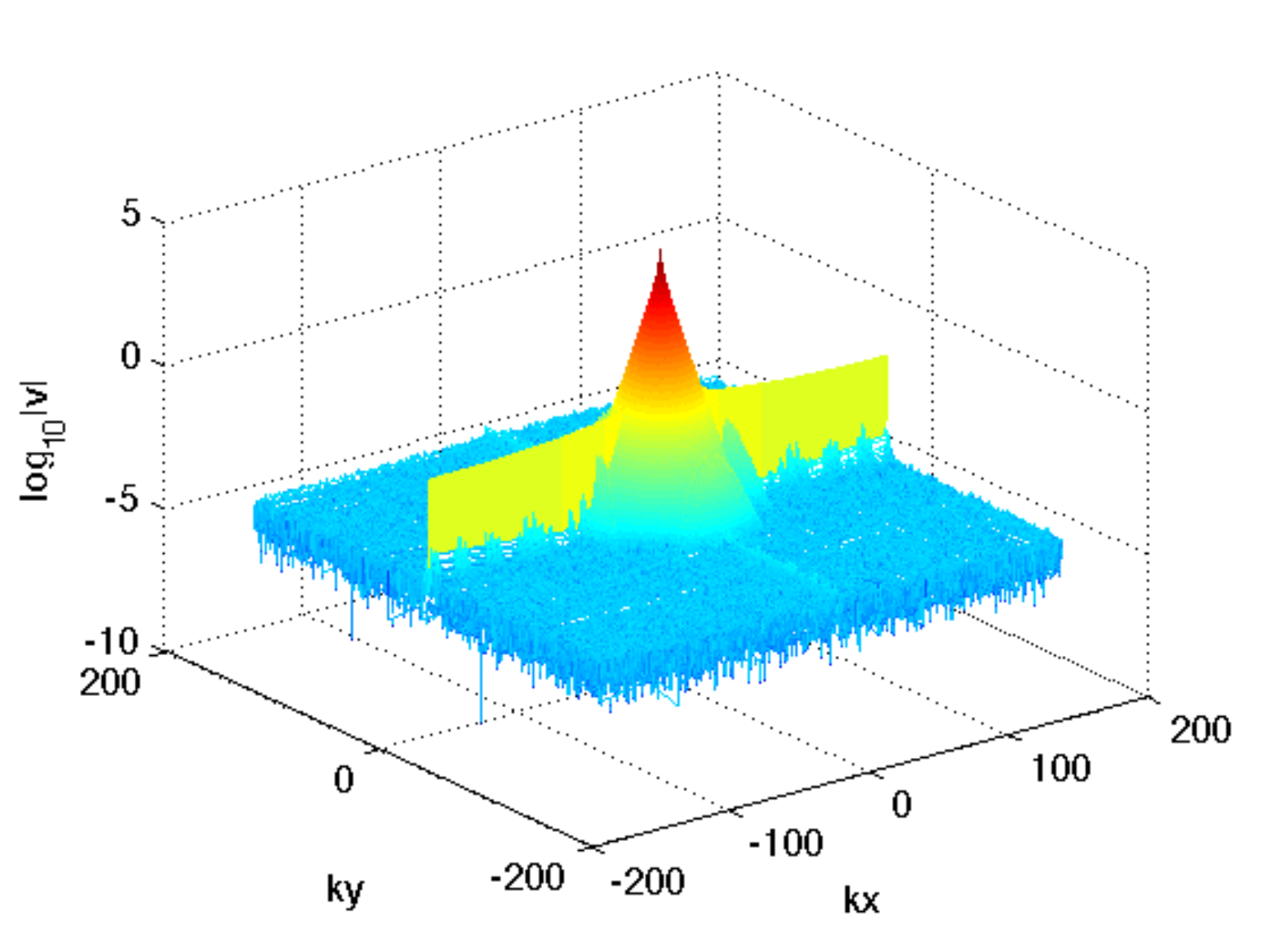}
\includegraphics[width=0.45\textwidth]{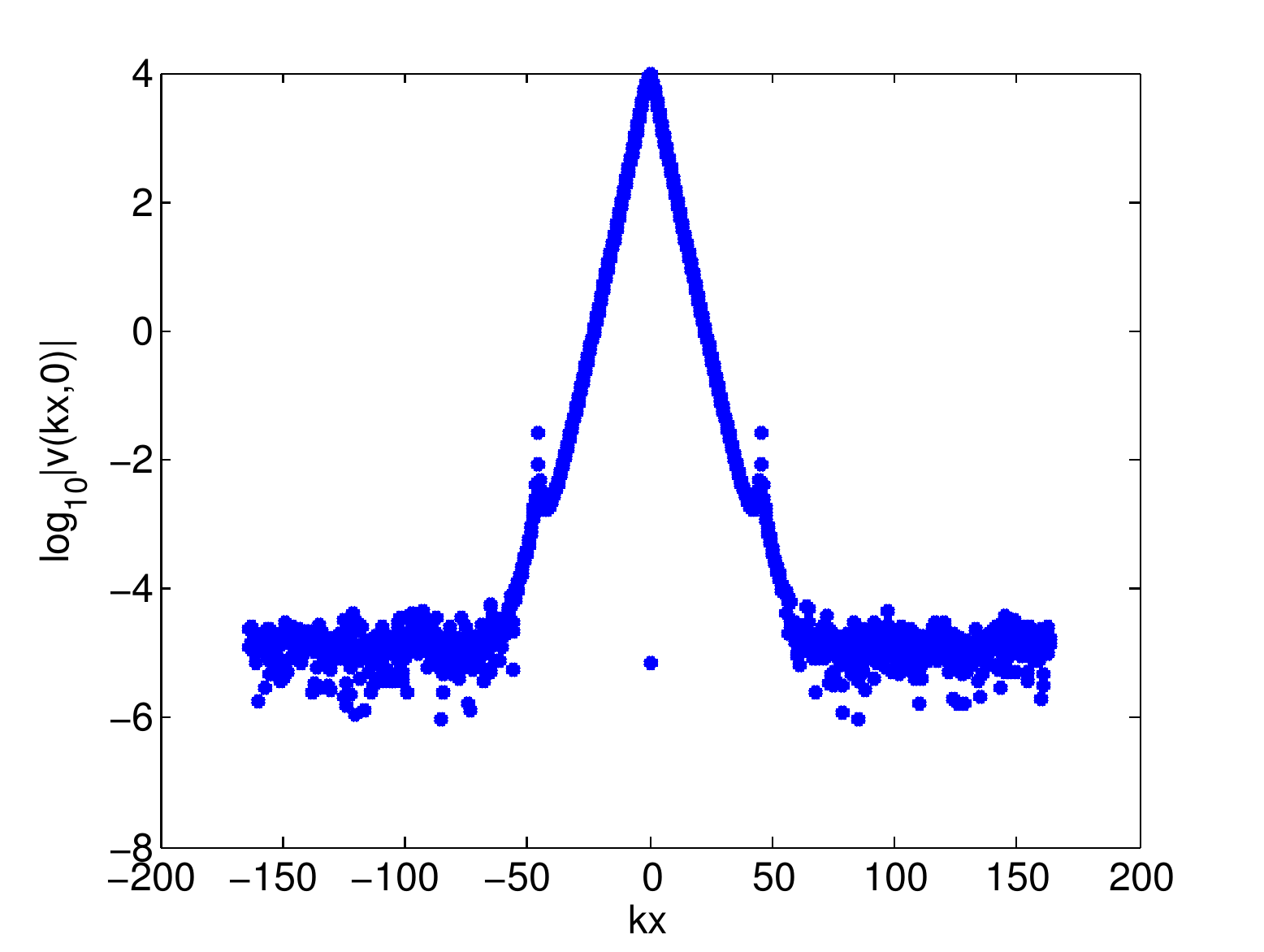}
\caption{Fourier coefficients of the solution to the focusing DS II equation 
(\ref{DSII}) 
for an initial condition of the form (\ref{lumpyd}) with $\kappa=0.9$  at $t=0$.}  
\label{ly09cf}
\end{figure}

\begin{figure}[htb!]
\centering
\includegraphics[width=0.5\textwidth]{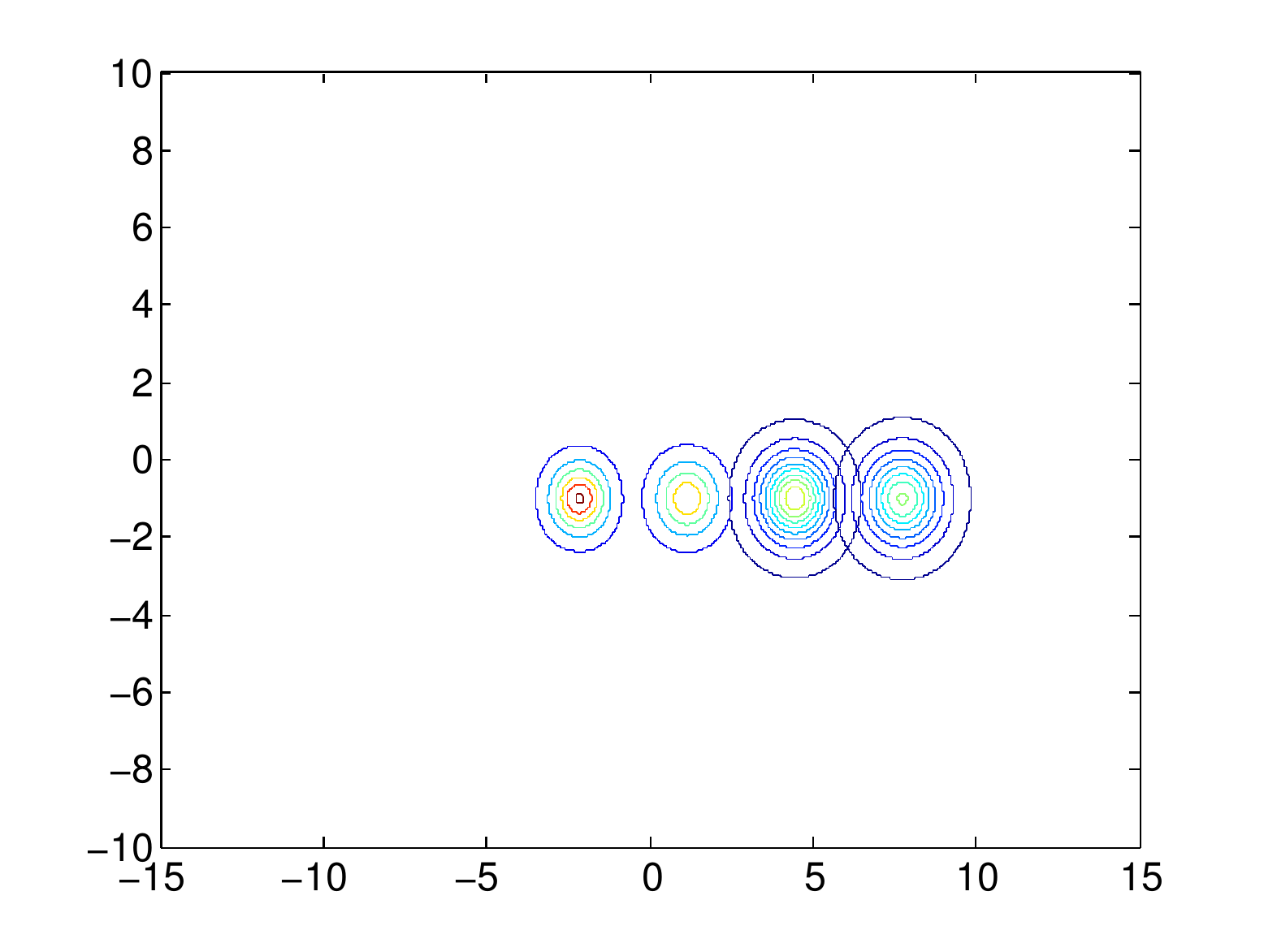}
\caption{Contour plot for a solution to the focusing DS II equation 
(\ref{DSII}) 
for an initial condition of the form (\ref{lumpyd}) with $\kappa=1.1$ 
for different times.}
\label{conluy11}
\end{figure}

\begin{figure}[htb!]
\centering
\includegraphics[width=0.45\textwidth]{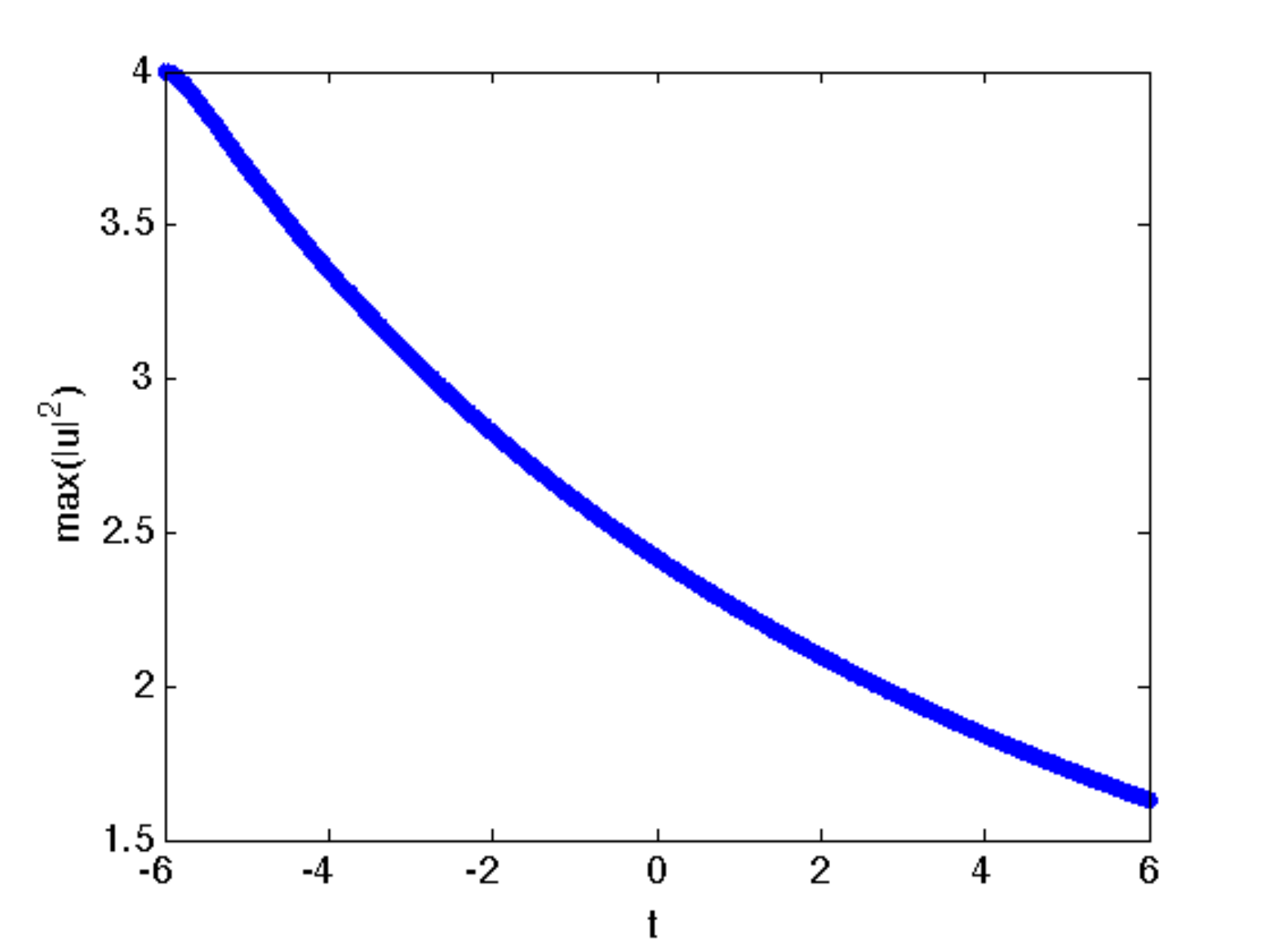}
\includegraphics[width=0.45\textwidth]{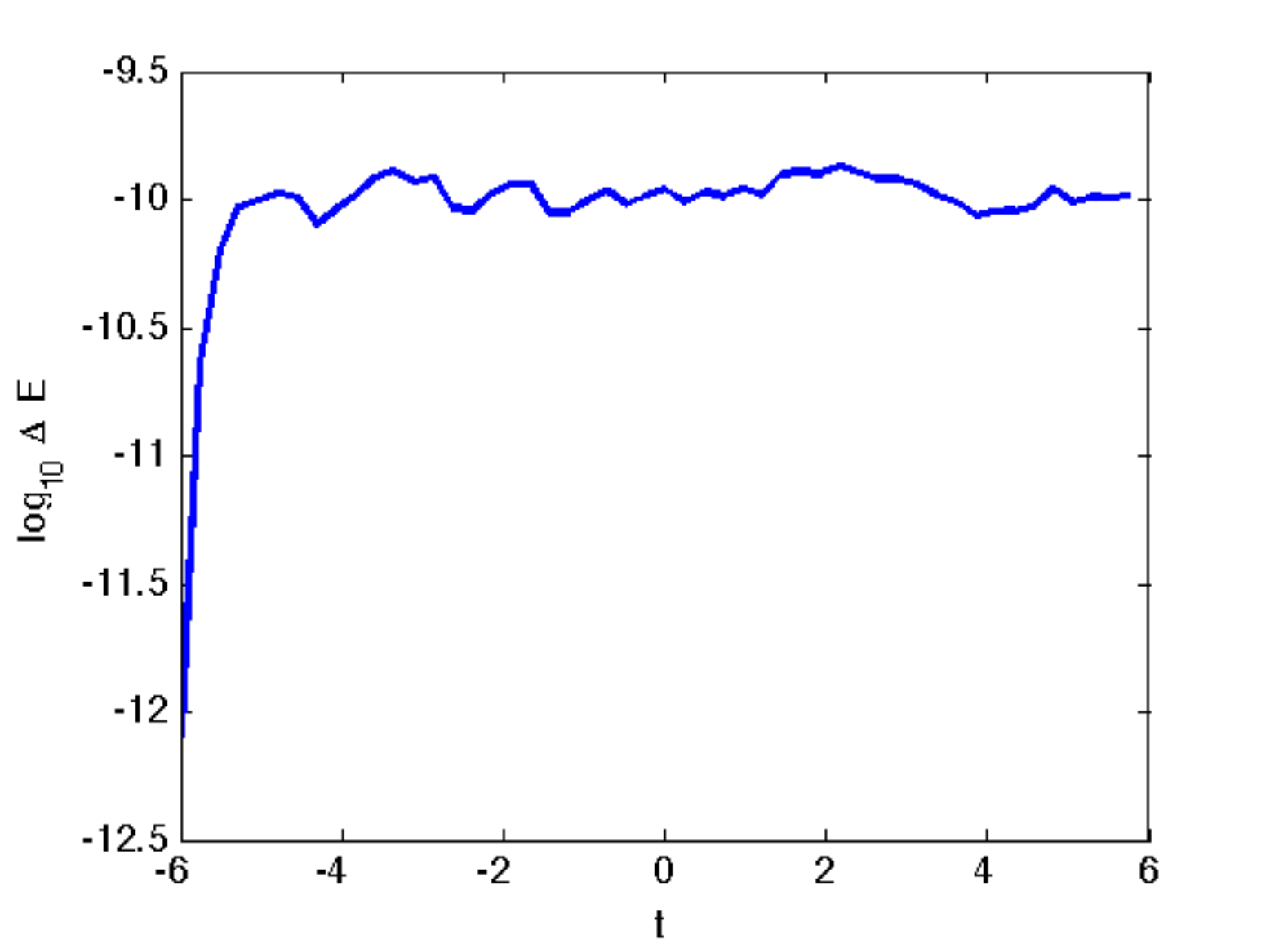}
\caption{Evolution of  $max(|u|^{2})$ and the numerically computed 
energy $E$for a solution to the focusing DS II equation 
(\ref{DSII}) 
for an initial condition of the form (\ref{lumpyd}) with $\kappa=1.1$.}
\label{ampluy11}
\end{figure}

\begin{figure}[htb!]
\centering
\includegraphics[width=0.45\textwidth]{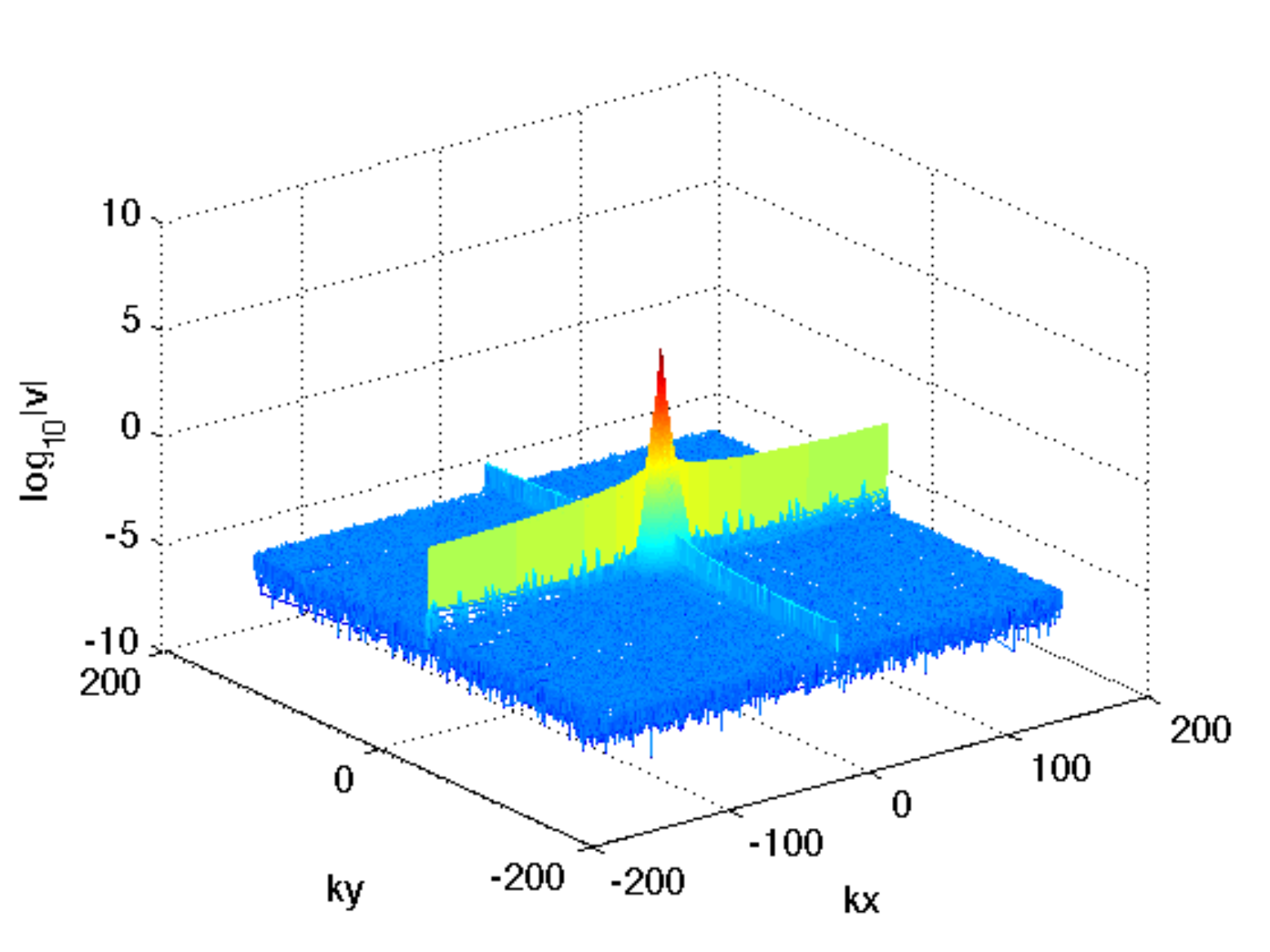}
\includegraphics[width=0.45\textwidth]{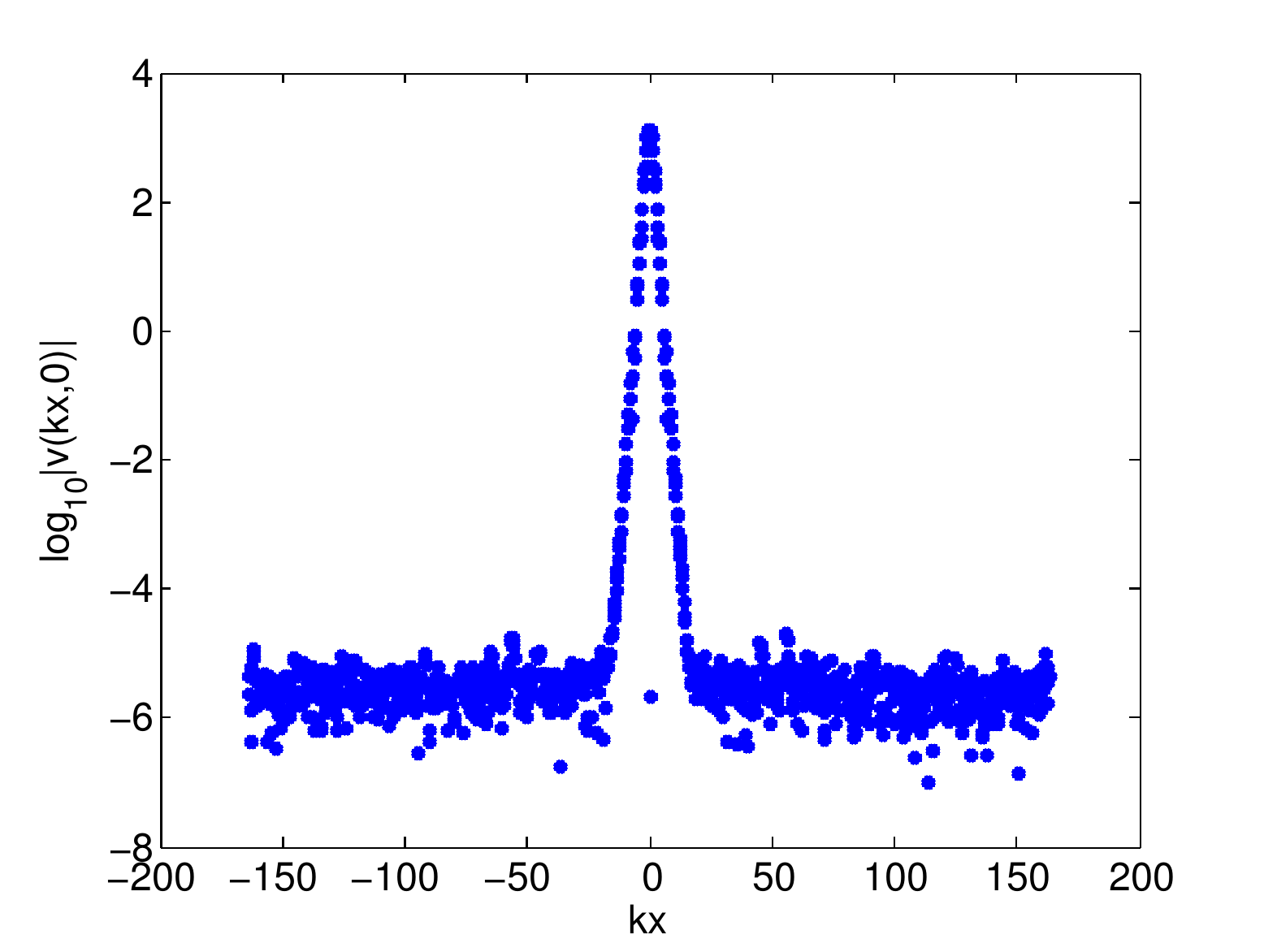}
\caption{Fourier coefficients of the solution to the focusing DS II equation 
(\ref{DSII}) 
for an initial condition of the form (\ref{lumpyd}) with $\kappa=1.1$  at $t=0$.}  
\label{ly11cf}
\end{figure}

We consider initial data of the form
\begin{equation}\label{lumpyd}
 u(x,y,-6) = u_{l}(x,\kappa y,-6),
\end{equation}
i.e., a deformed (in $y$-direction) initial lump in this 
subsection.
The computations are carried out with $N_{x}=N_{y}=2^{14}$ points for 
$x\times y \in [-50\pi, 50\pi] \times [-50\pi, 50\pi]$ and $t\in[-6,6]$.
\\
\\
For $\kappa=0.9$, the resulting solution loses speed and width as can 
be seen in Fig.~\ref{conluy09}.
Its height and energy grow,  but both stay finite, see Fig.~\ref{ampluy09}. It is 
possible that the solution eventually blows up, but not on the 
time scales studied here.
\\
The  Fourier coefficients at $t=0$  in Fig.~\ref{ly09cf} show  
the wanted spatial resolution.
\\
\\
For $\kappa=1.1$, we observe the opposite behavior in 
Fig.~\ref{conluy11}. The solution 
travels with higher speed than the initial lump and 
is broadened. The energy does not show any sudden change, see Fig.~\ref{ampluy11}.
It seems that the initial pulse will asymptotically disperse.
The  Fourier coefficients at $t=0$ in Fig.~\ref{ly11cf} show the wanted spatial resolution.

\section{Perturbations of the Ozawa solution}

In this section we study as for the lump in the previous section 
various perturbations of initial data for the Ozawa solution to test 
whether blowup is generic for the focusing DS II equation. 

\subsection{Perturbation of the Ozawa solution by multiplication with 
a scalar factor}
We consider initial data of the form  
\begin{equation}
    u(x,y,0) = 2C \frac{ \exp \left( -i (x^2 - y^2)
    \right)}{1+x^2 +y^2},
    \label{ozawafac}
\end{equation}
i.e., initial data of the Ozawa solution multiplied by a scalar factor. 
The computation is carried out with $N_{x}=N_{y}=2^{15}$ points for 
$x\times y \in [-20\pi, 20\pi] \times [-20\pi, 20\pi]$.\\

\begin{figure}[htb!]
\centering
\includegraphics[width=0.45\textwidth]{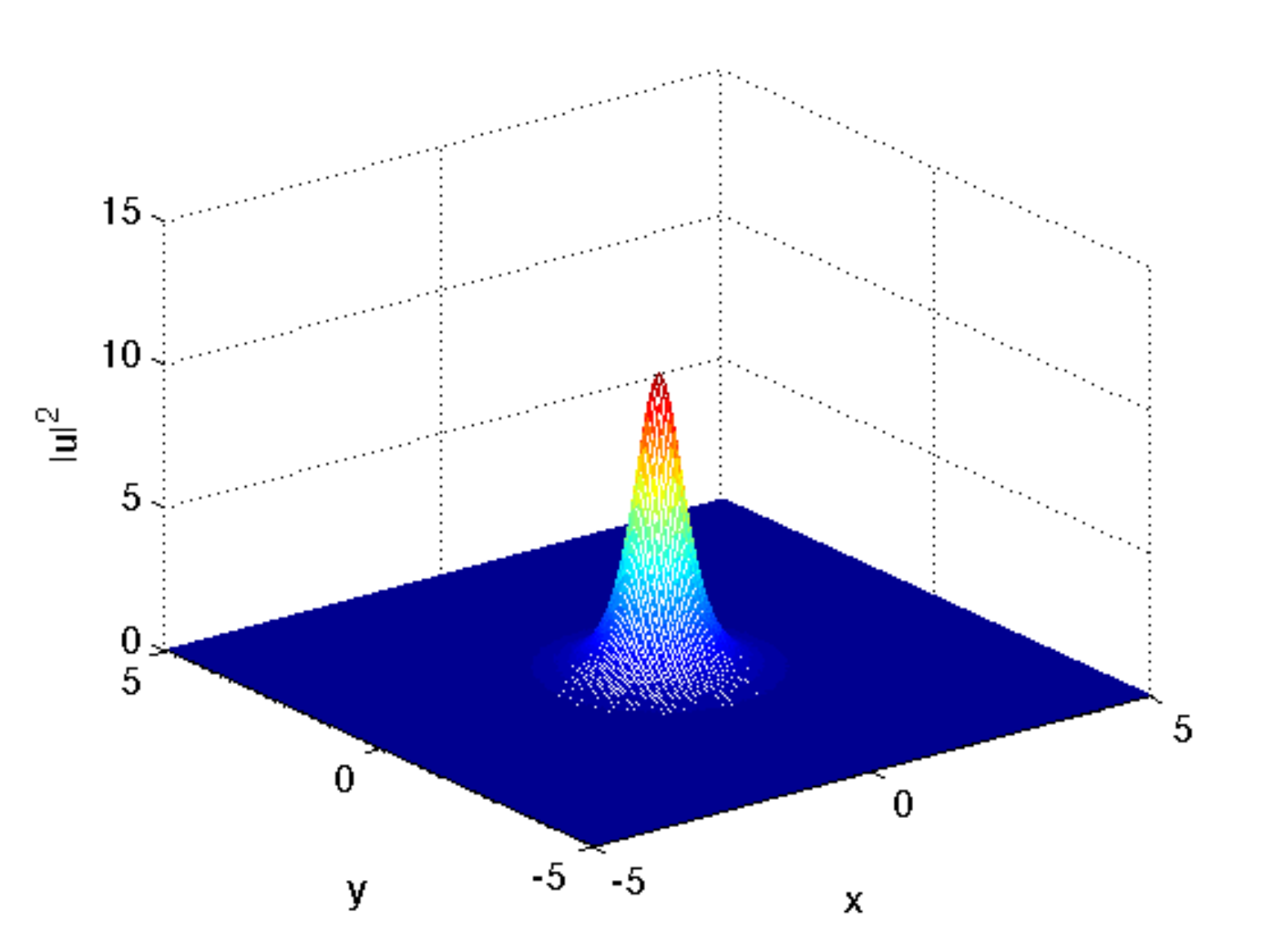}
\includegraphics[width=0.45\textwidth]{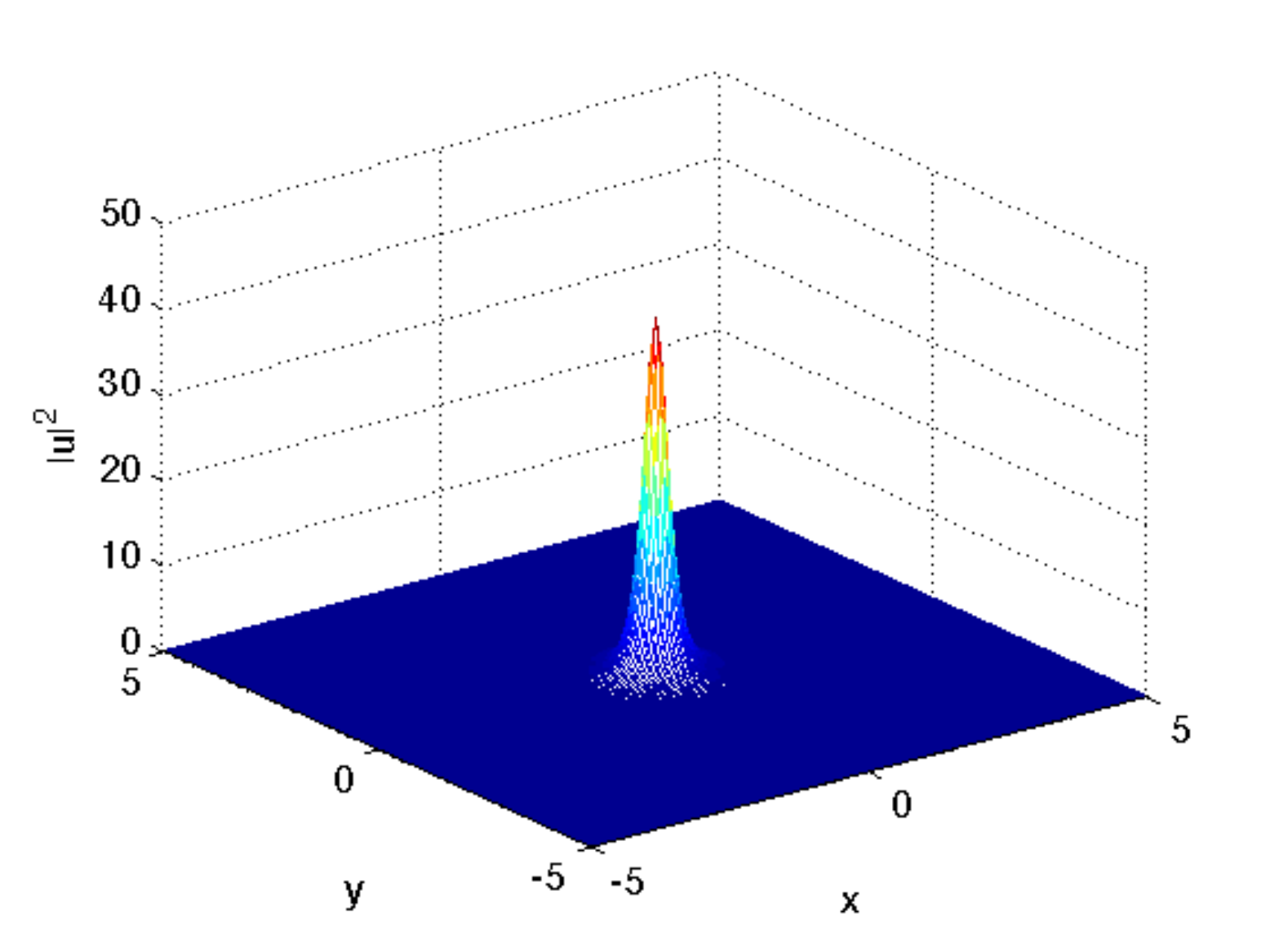}\\
\includegraphics[width=0.45\textwidth]{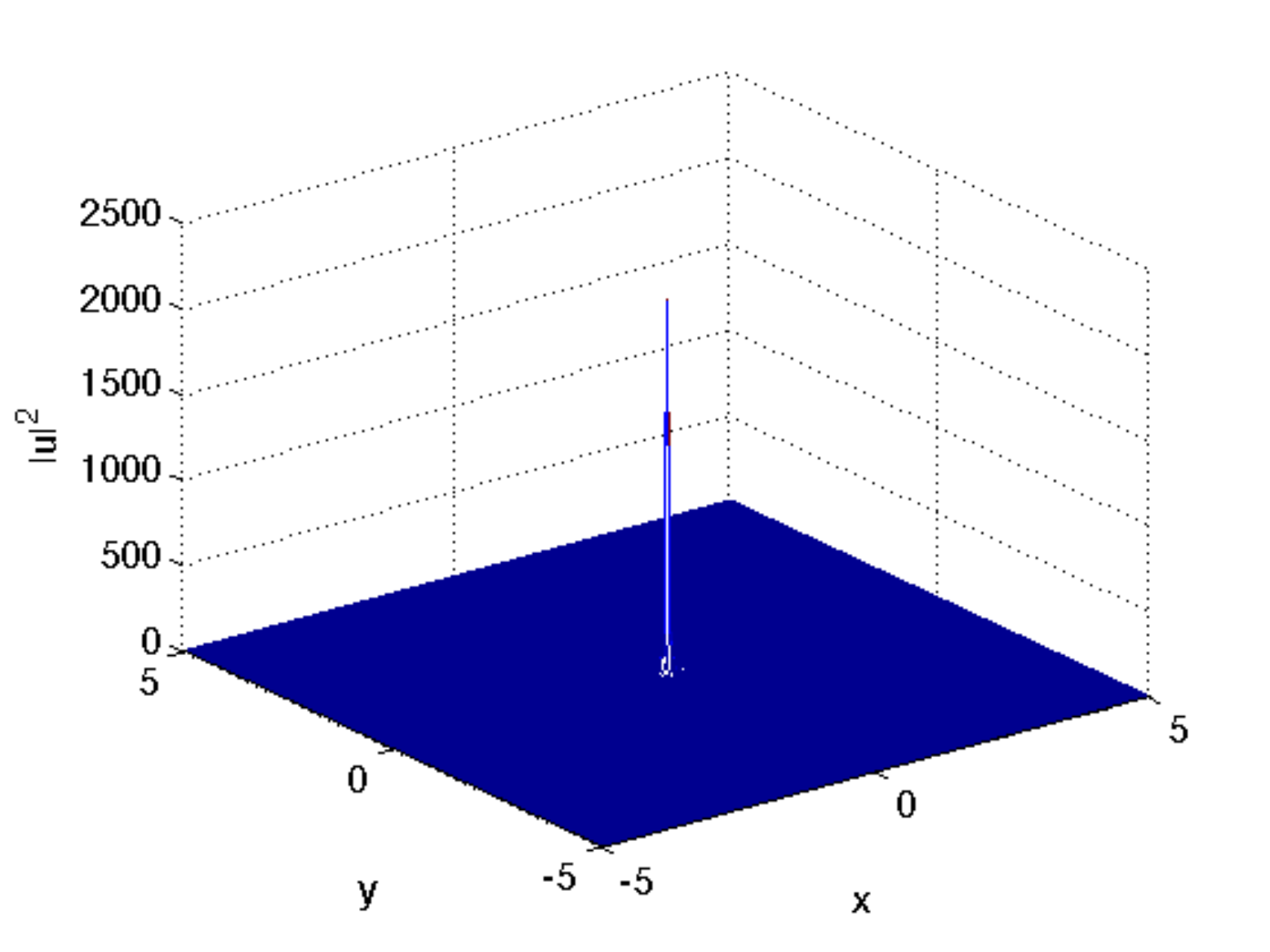}
\includegraphics[width=0.45\textwidth]{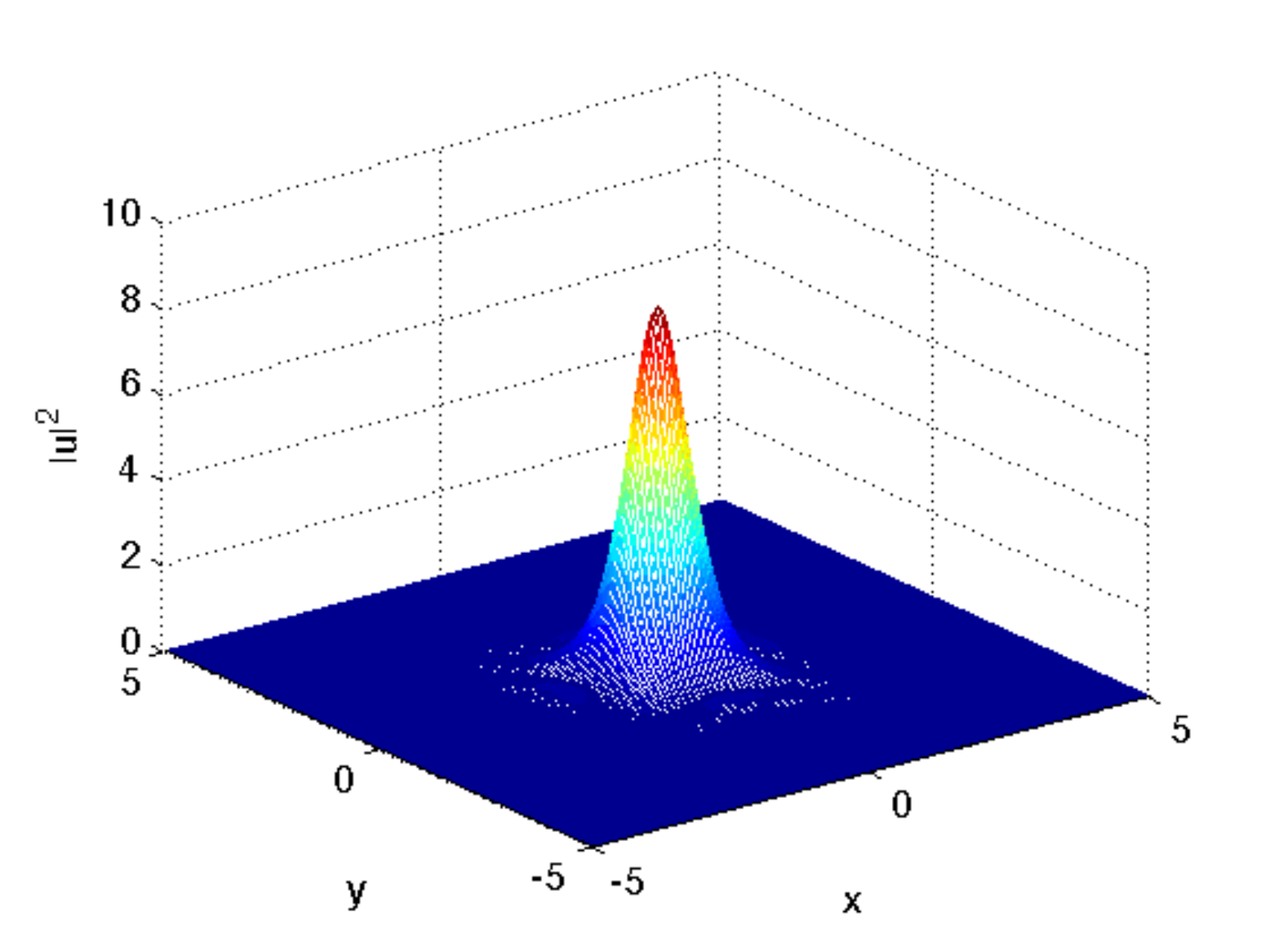}\\
\caption{Solution to the focusing DS II equation (\ref{DSII})
for an initial condition of the form (\ref{ozawafac}) with $C=1.1$ for 
$t=0.075$ and $t=0.15$ in the first row and $t=0.225$ and $t=0.3$ 
below.} 
\label{uoz11}
\end{figure}

\begin{figure}[htb!]
\centering
\includegraphics[width=0.45\textwidth]{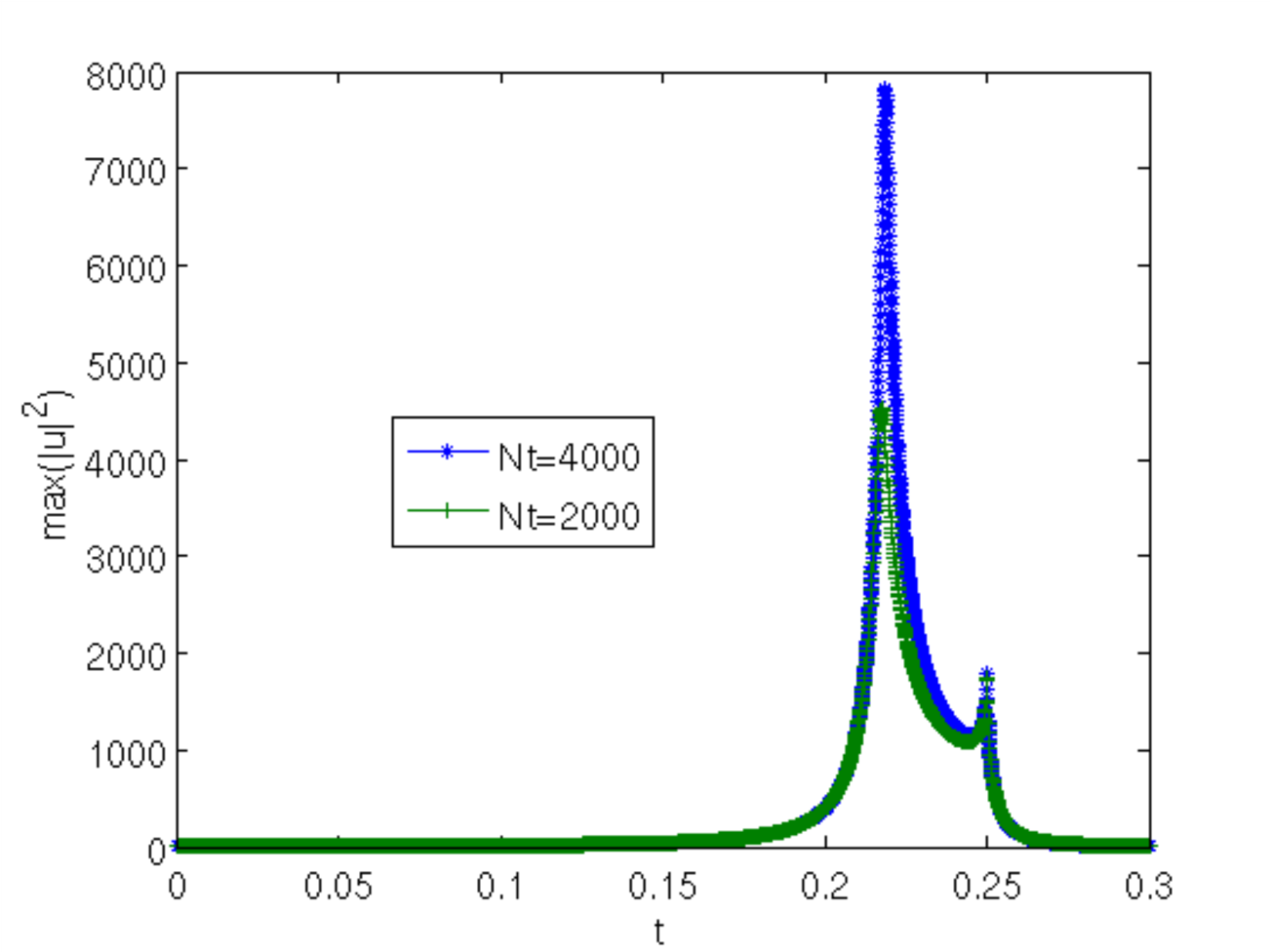}
\includegraphics[width=0.45\textwidth]{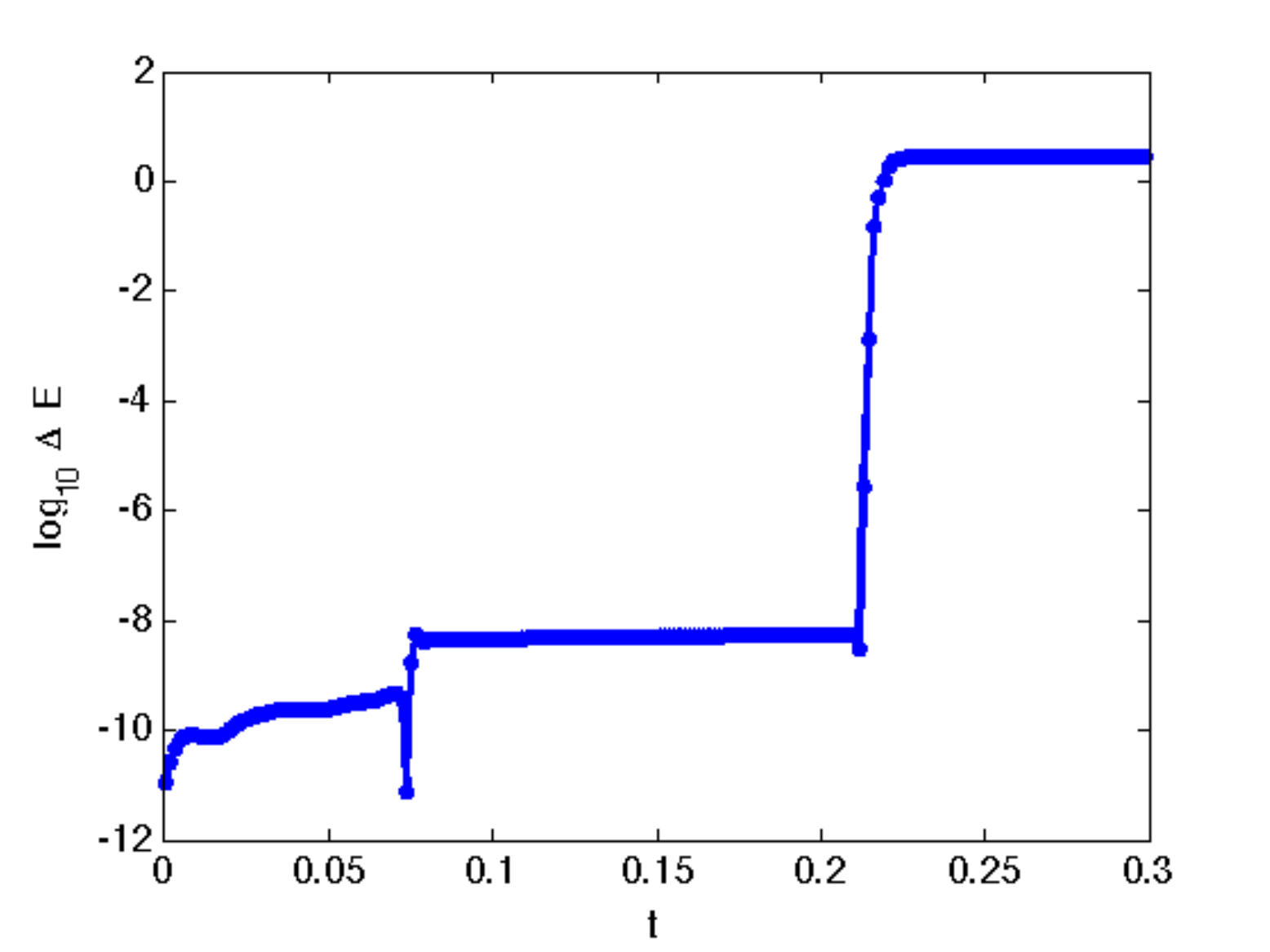}
\caption{Evolution of $max(|u|^{2})$ and the numerically computed energy 
for an initial condition of the form (\ref{ozawafac}) with $C=1.1$.}
\label{ampluoz11}
\end{figure}

\begin{figure}[htb!]
\centering
\includegraphics[width=0.45\textwidth]{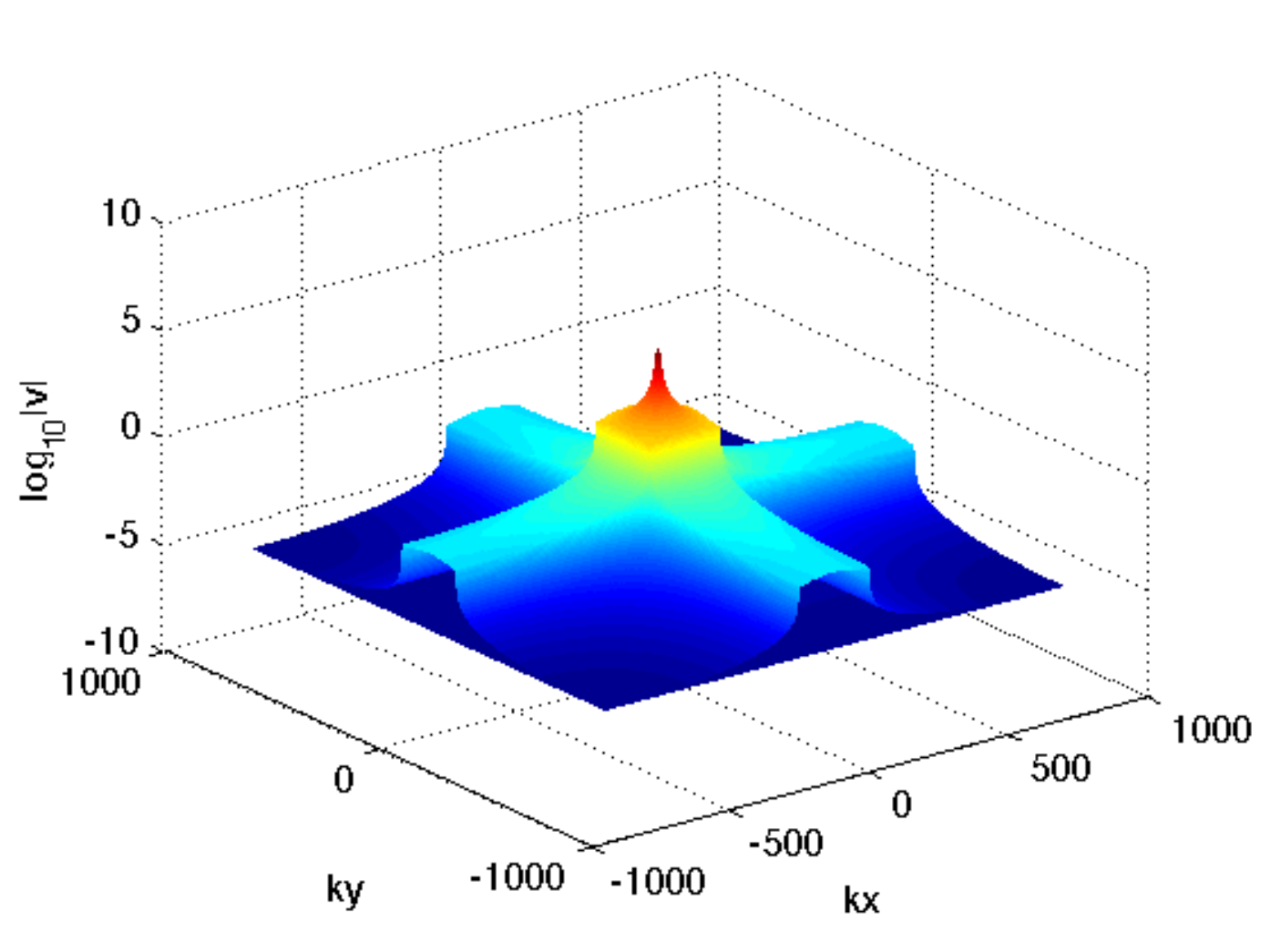}
\includegraphics[width=0.45\textwidth]{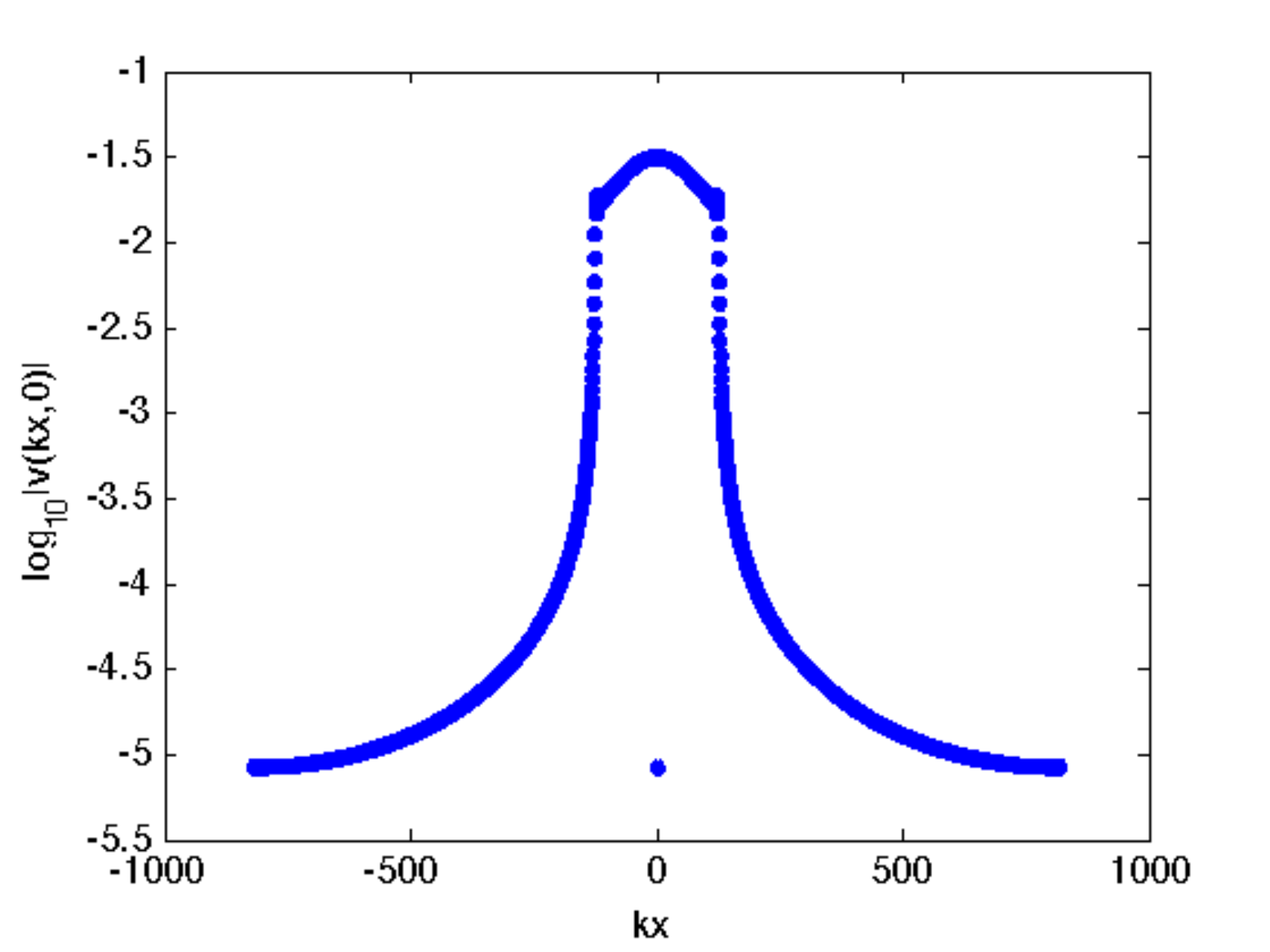}
\caption{Fourier coefficients of solution to the focusing DS II equation (\ref{DSII})
for an initial condition of the form (\ref{ozawafac}) with $C=1.1$ at $t=0.15$.} 
\label{oz11cf}
\end{figure}

For $C=1.1$, and 
$N_t = 2000$, 
we show the behavior of $|u|^2$ at different times in Fig. \ref{uoz11}.
The time evolution of $\underset{x,y}{\max} |u(x,y,t)|^2$ and the 
numerically computed energy are shown in Fig.~\ref{ampluoz11}.
We observe an $L_{\infty}$ blowup at the time $t_c\sim0.2210$.
%
%
The Fourier coefficients at $t=0.15$ (before the blowup) in Fig.~\ref{oz11cf} show the 
wanted spatial resolution.\\

\begin{figure}[htb!]
\centering
\includegraphics[width=0.45\textwidth]{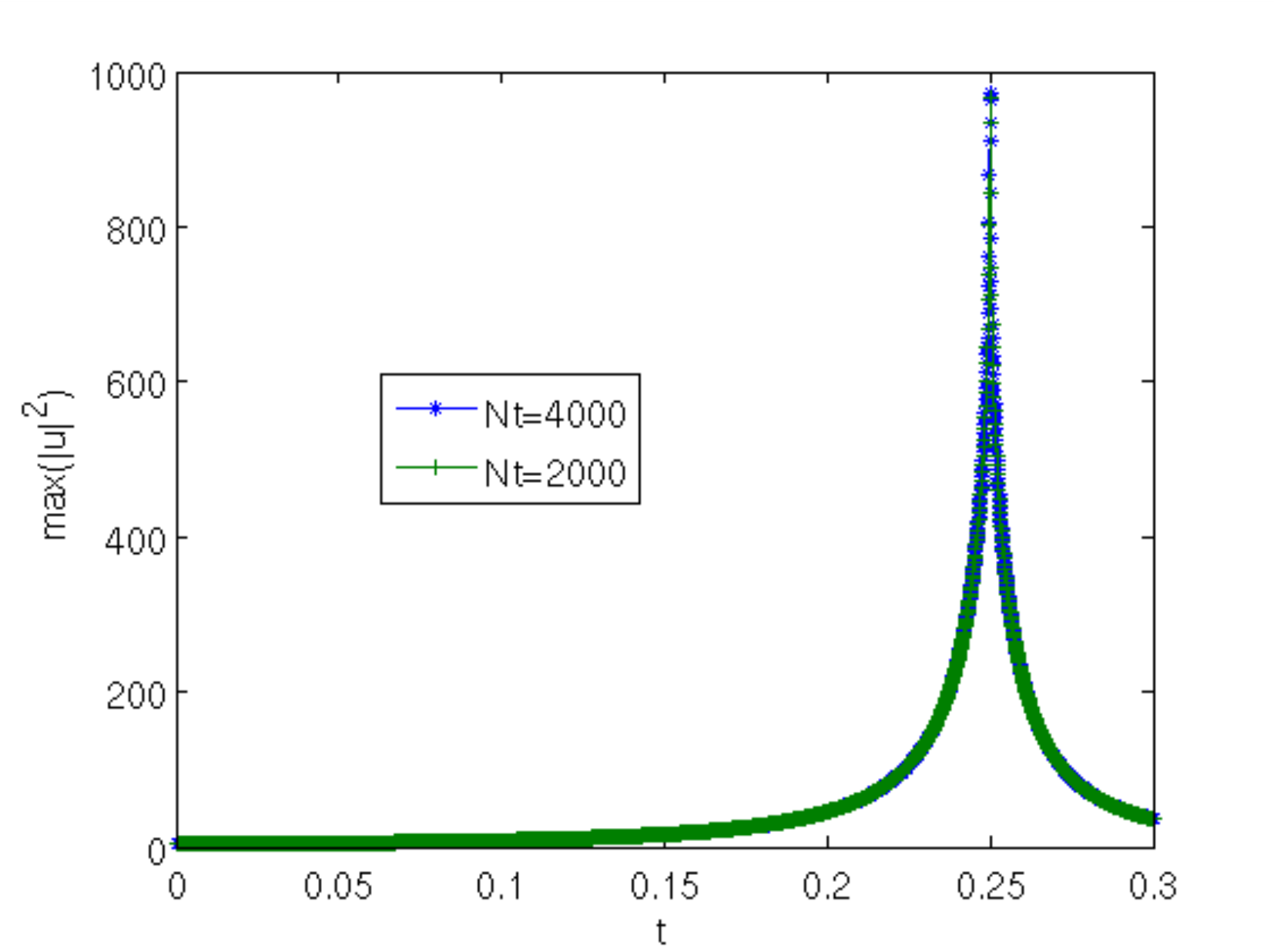}
\includegraphics[width=0.45\textwidth]{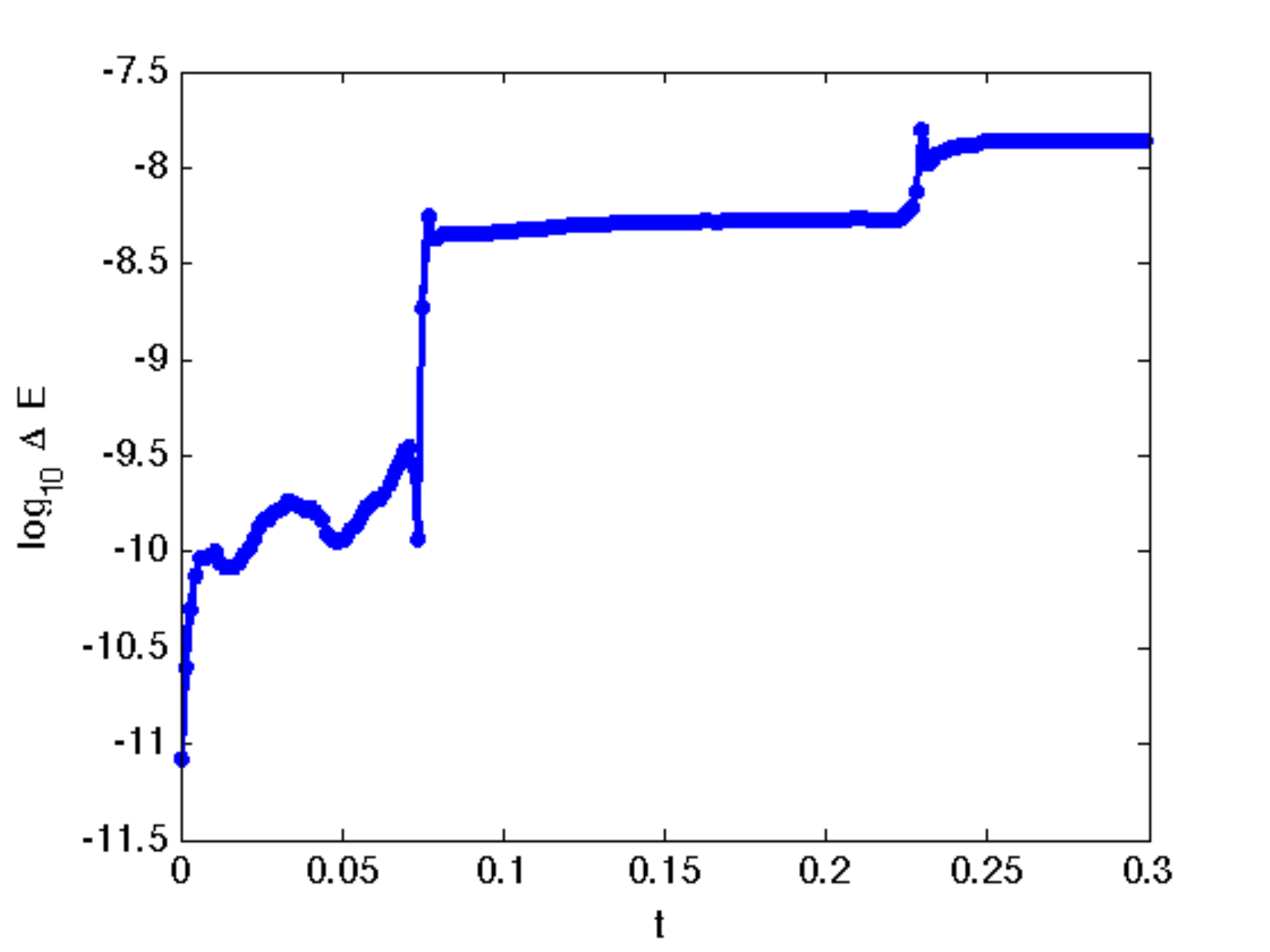}
\caption{Evolution of $max(|u|^{2})$ in dependence of time, for an initial condition 
of the form (\ref{ozawafac}) with $C=0.9$.}
\label{ampluoz09}
\end{figure}

\begin{figure}[htb!]
\centering
\includegraphics[width=0.45\textwidth]{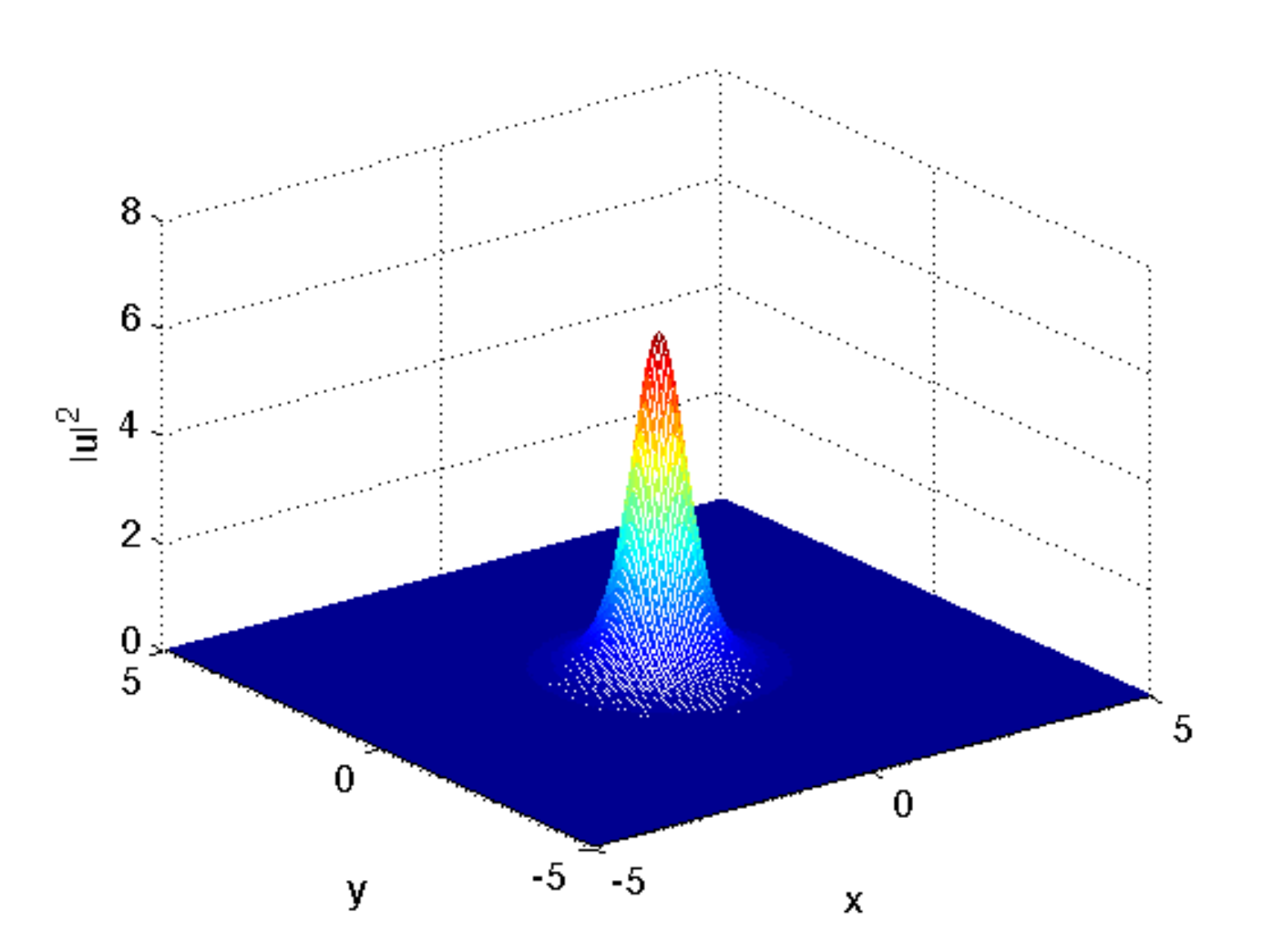}
\includegraphics[width=0.45\textwidth]{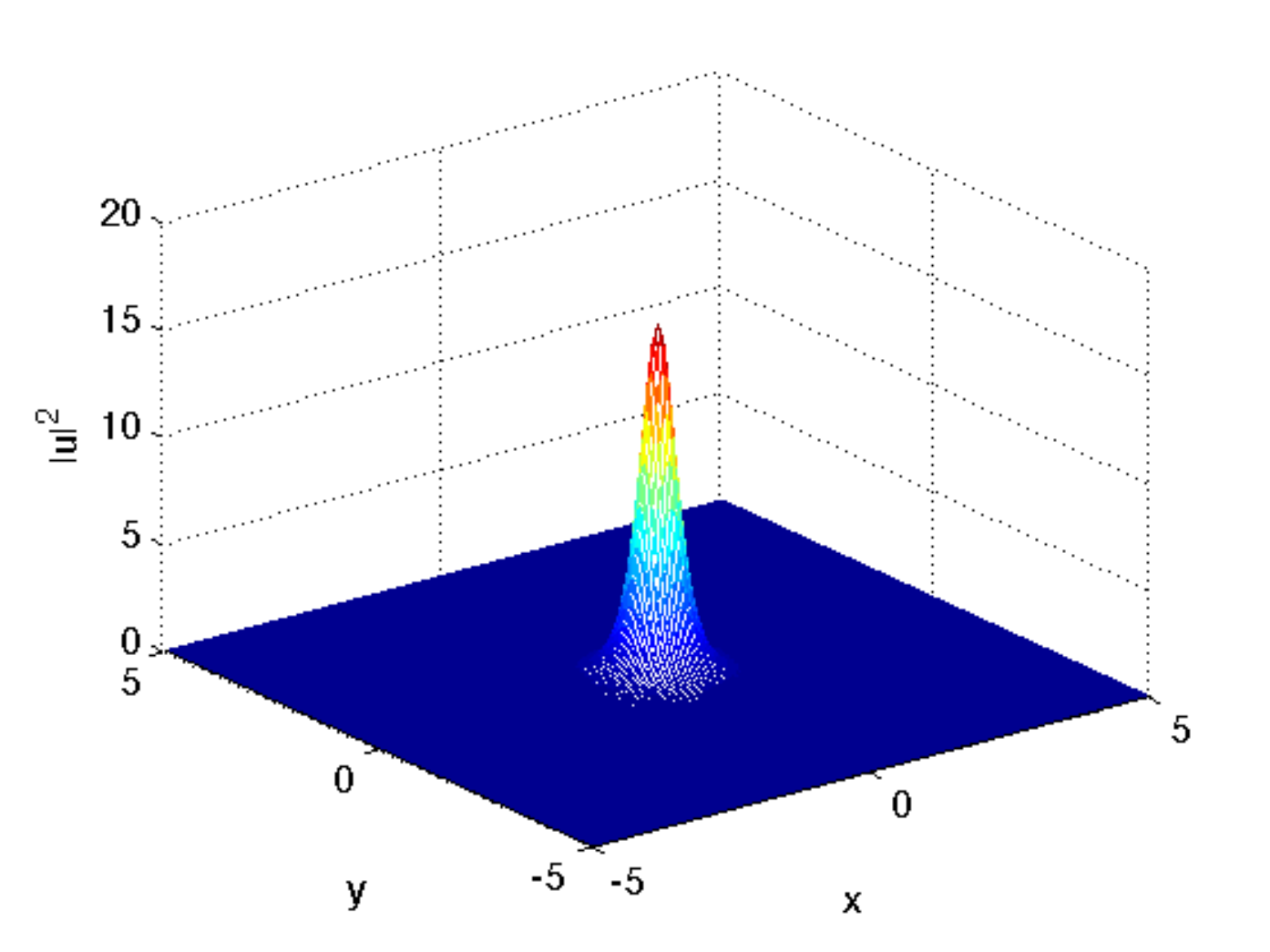}
\includegraphics[width=0.45\textwidth]{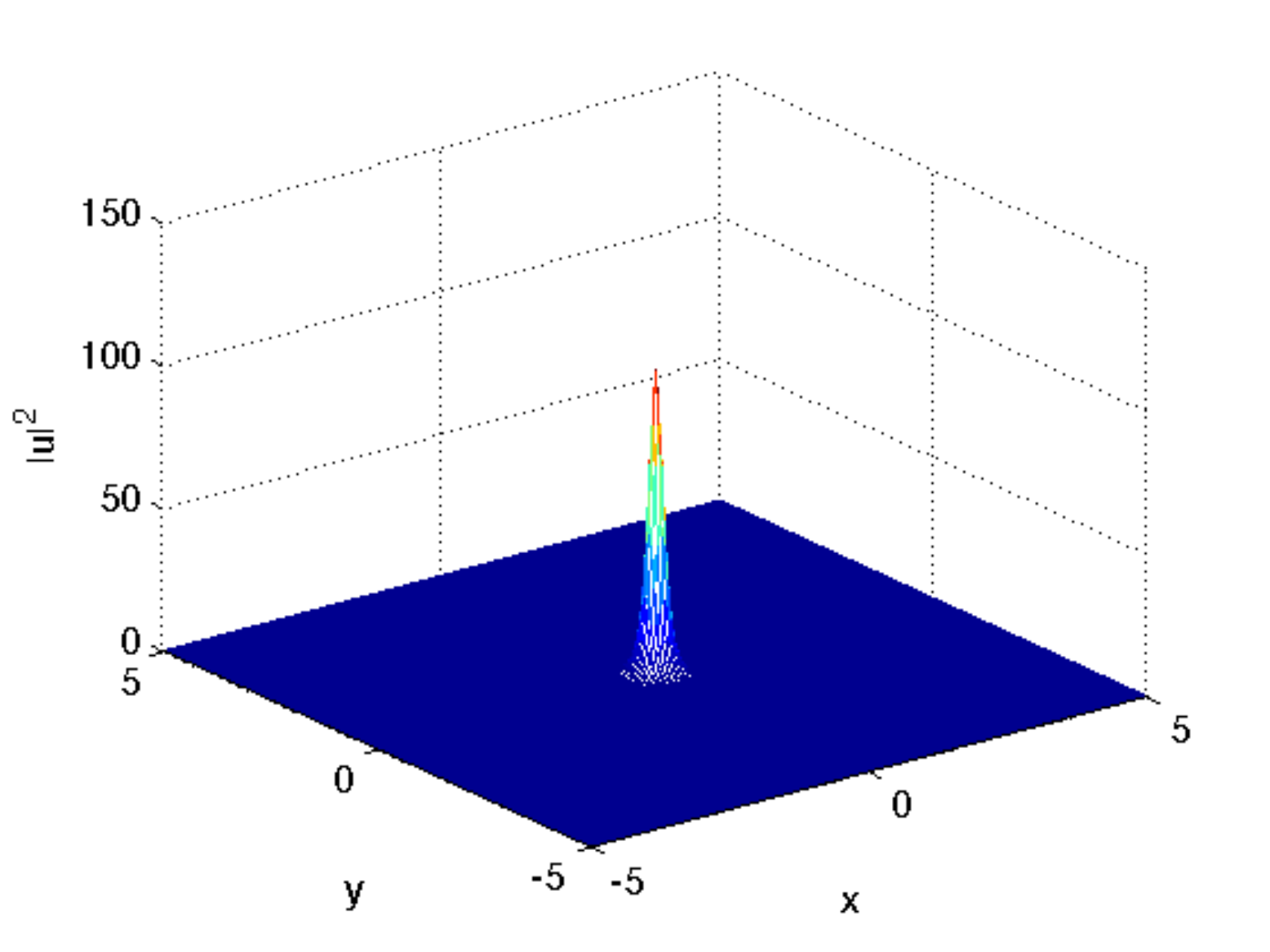}
\includegraphics[width=0.45\textwidth]{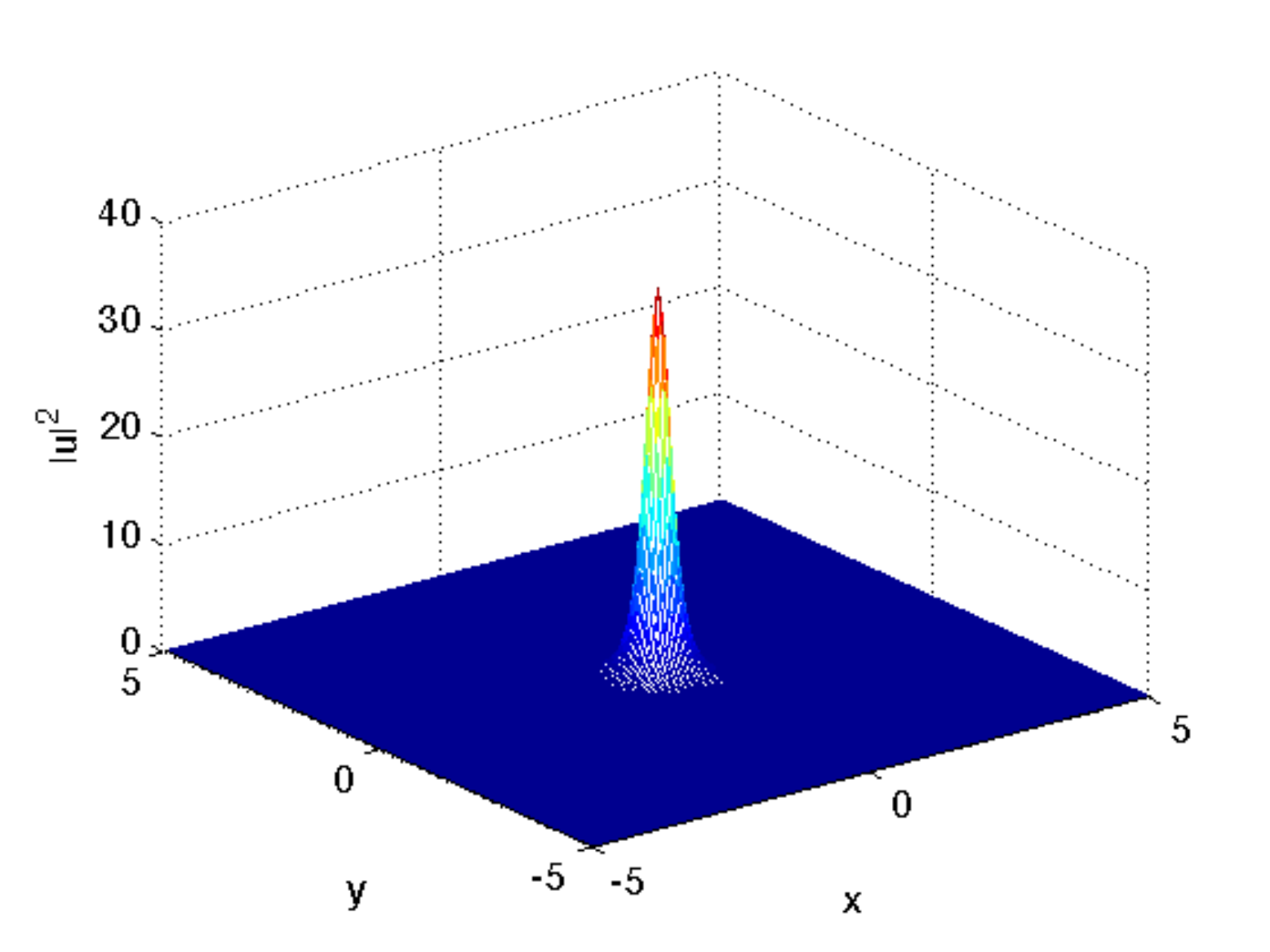}
\caption{Solution to the focusing DS II equation (\ref{DSII}) for an initial condition 
of the form (\ref{ozawafac}) with $C=0.9$, $N_t=2000$ for 
$t=0.075$ and $t=0.15$ in the first row and $t=0.225$ and $t=0.3$ 
below.}
\label{uoz09}
\end{figure}

\begin{figure}[htb!]
\centering
\includegraphics[width=0.45\textwidth]{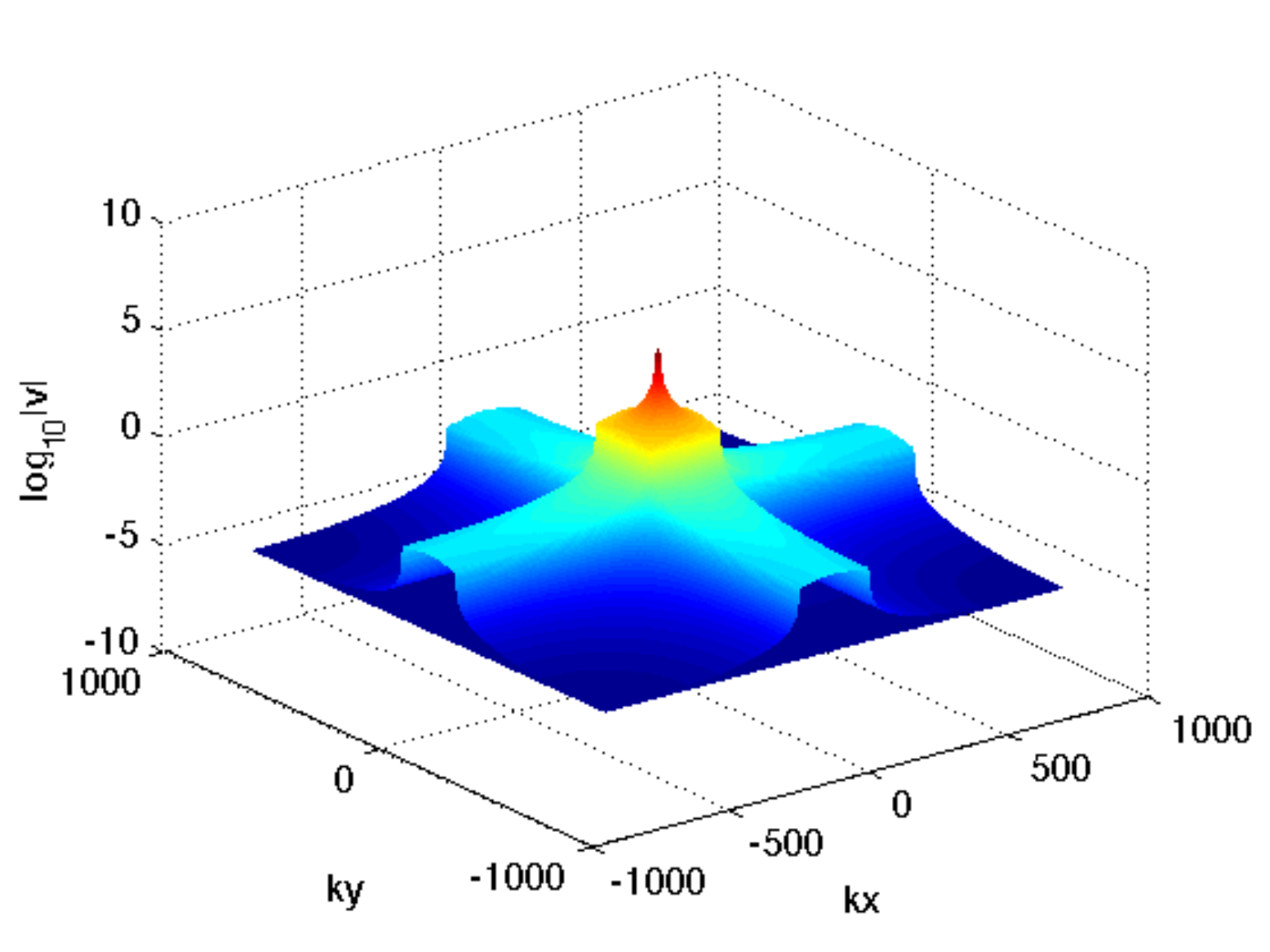}
\includegraphics[width=0.45\textwidth]{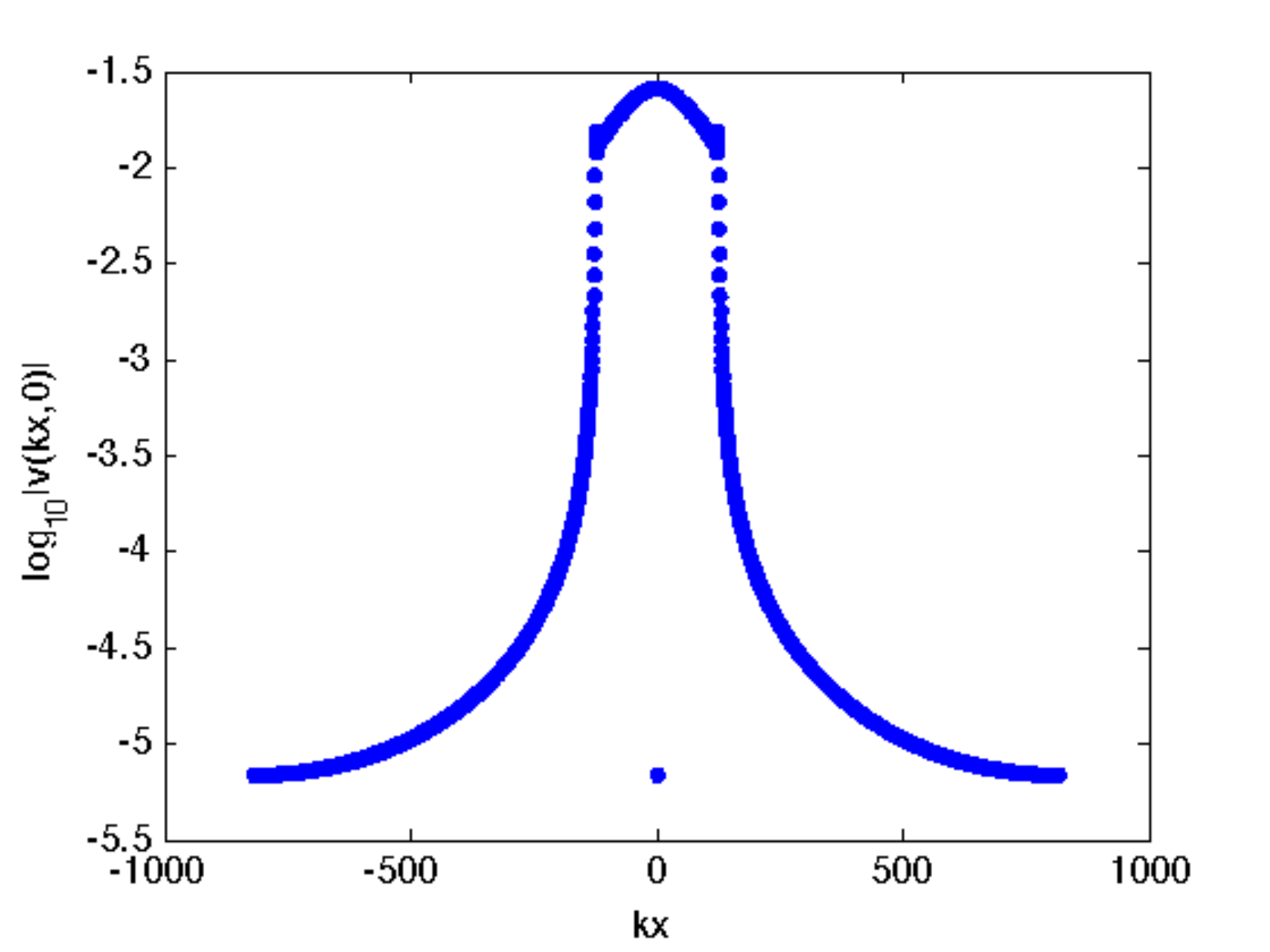}
\caption{Fourier coefficients of the solution to the focusing DS II equation (\ref{DSII}) at $t=0.15$
for an initial condition of the form (\ref{ozawafac}) with $C=0.9$.} 
\label{oz09cf}
\end{figure}

For $C=0.9$, the initial pulse grows until it reaches its maximal height at 
$t=0.2501$, but there is no indication for blowup, see Fig. \ref{ampluoz09}. 
The solution at different times can be seen in Fig.~\ref{uoz09}.
The Fourier coefficients in 
Fig.~\ref{oz09cf} show again that the wanted spatial resolution is 
achieved. 

Thus for initial data given by the Ozawa solution multiplied with a 
factor $C$, we find that for 
$C>1$, blow up seems to occur before the critical time of the 
Ozawa solution, 
and for $C<1$ the solution grows until $t=0.25$ but does not 
blow up. Consequently the Ozawa initial data seem to be critical in this 
sense that data of this form with smaller norm do not blow up.

\subsection{Perturbation of the Ozawa solution with a Gaussian}

\begin{figure}[htb!]
\centering
\includegraphics[width=0.45\textwidth]{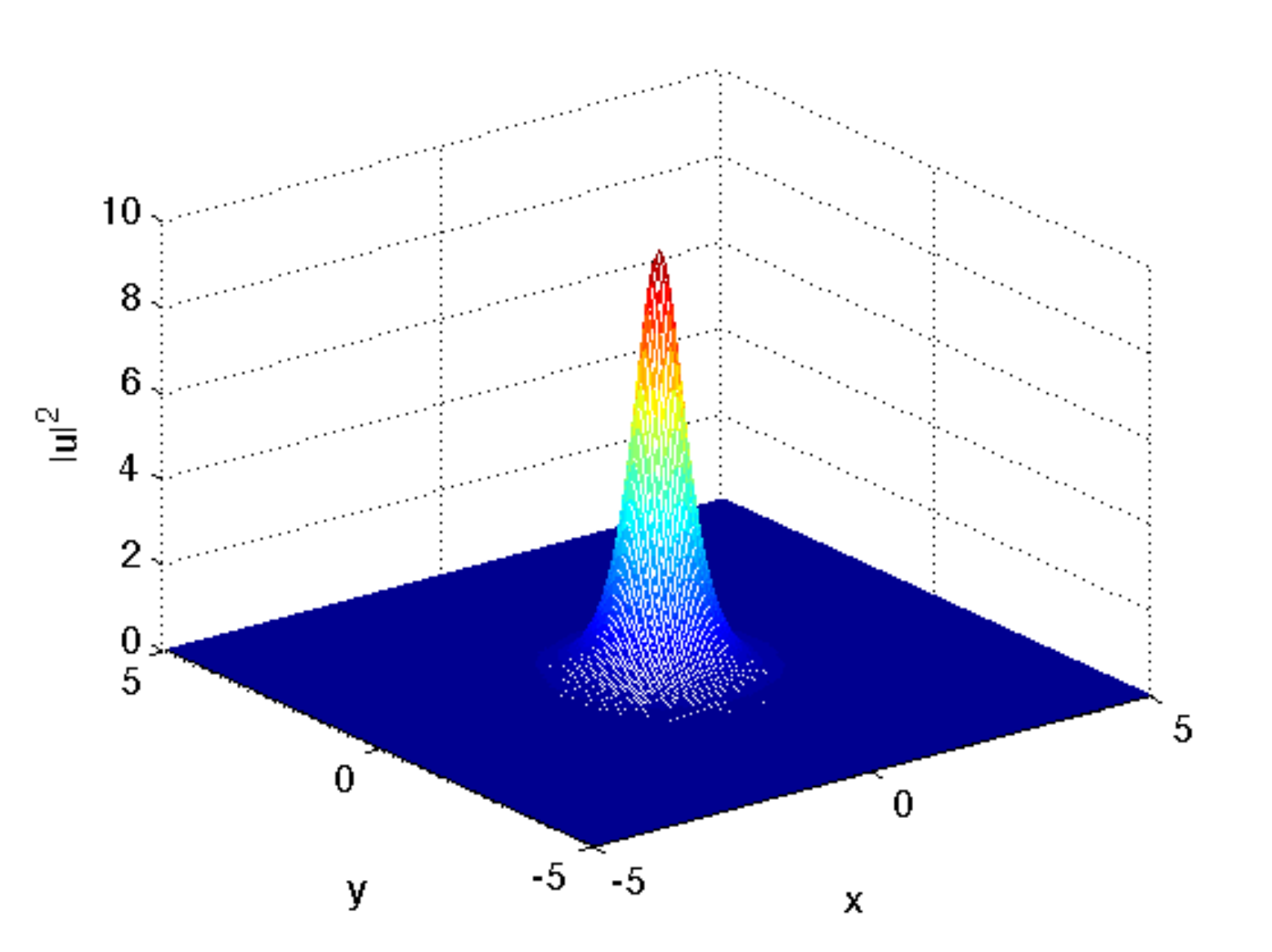}
\includegraphics[width=0.45\textwidth]{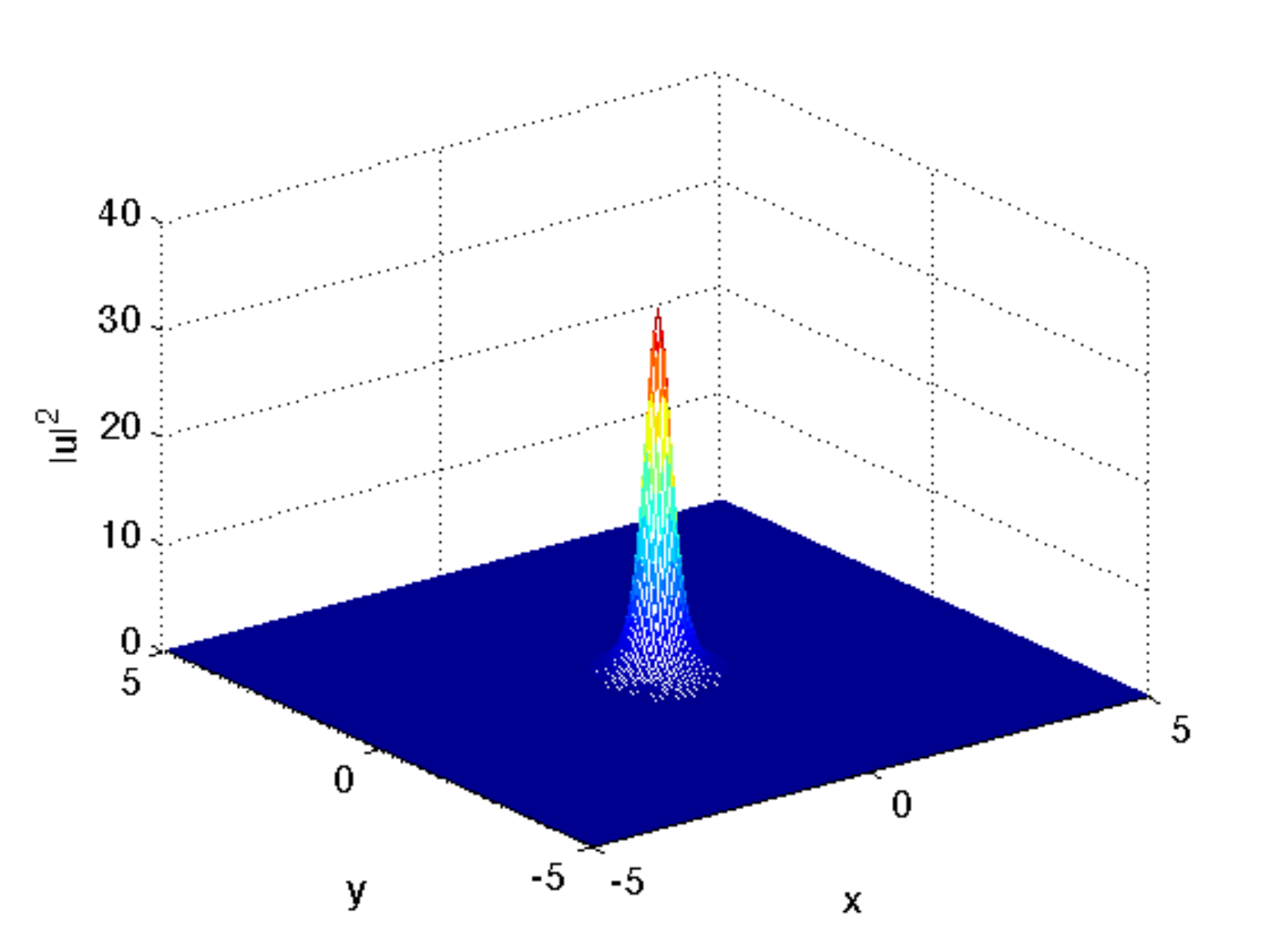}
\includegraphics[width=0.45\textwidth]{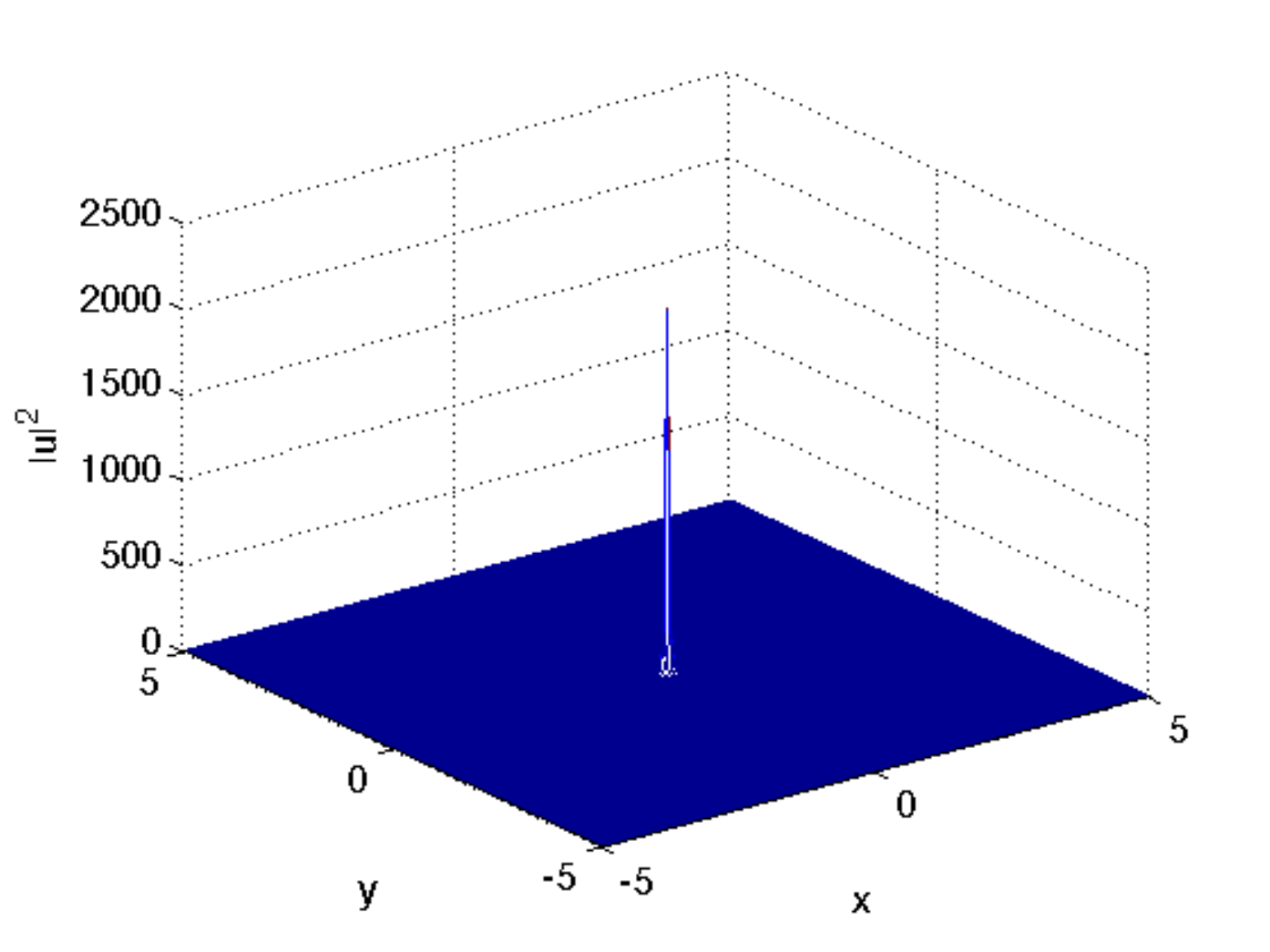} 
\includegraphics[width=0.45\textwidth]{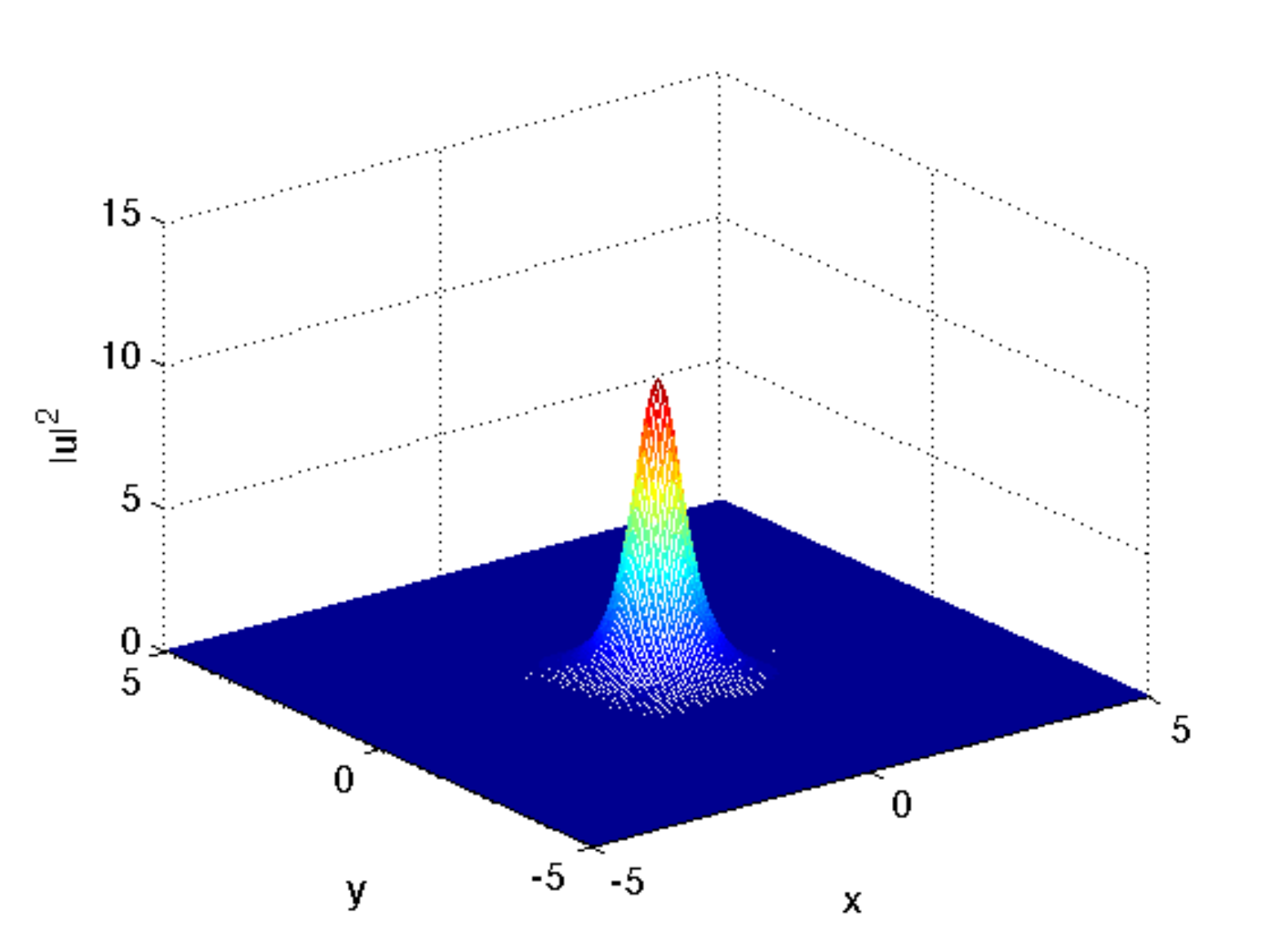}
\caption{Solution to the focusing DS II equation (\ref{DSII})
for an initial condition of the form (\ref{ozawagauss}) with $D=0.1$ for 
$t=0.075$ and $t=0.15$ in the first row and $t=0.225$ and $t=0.3$ 
below .} 
\label{uozpg1}
\end{figure}

\begin{figure}[htb!]
\centering
\includegraphics[width=0.45\textwidth]{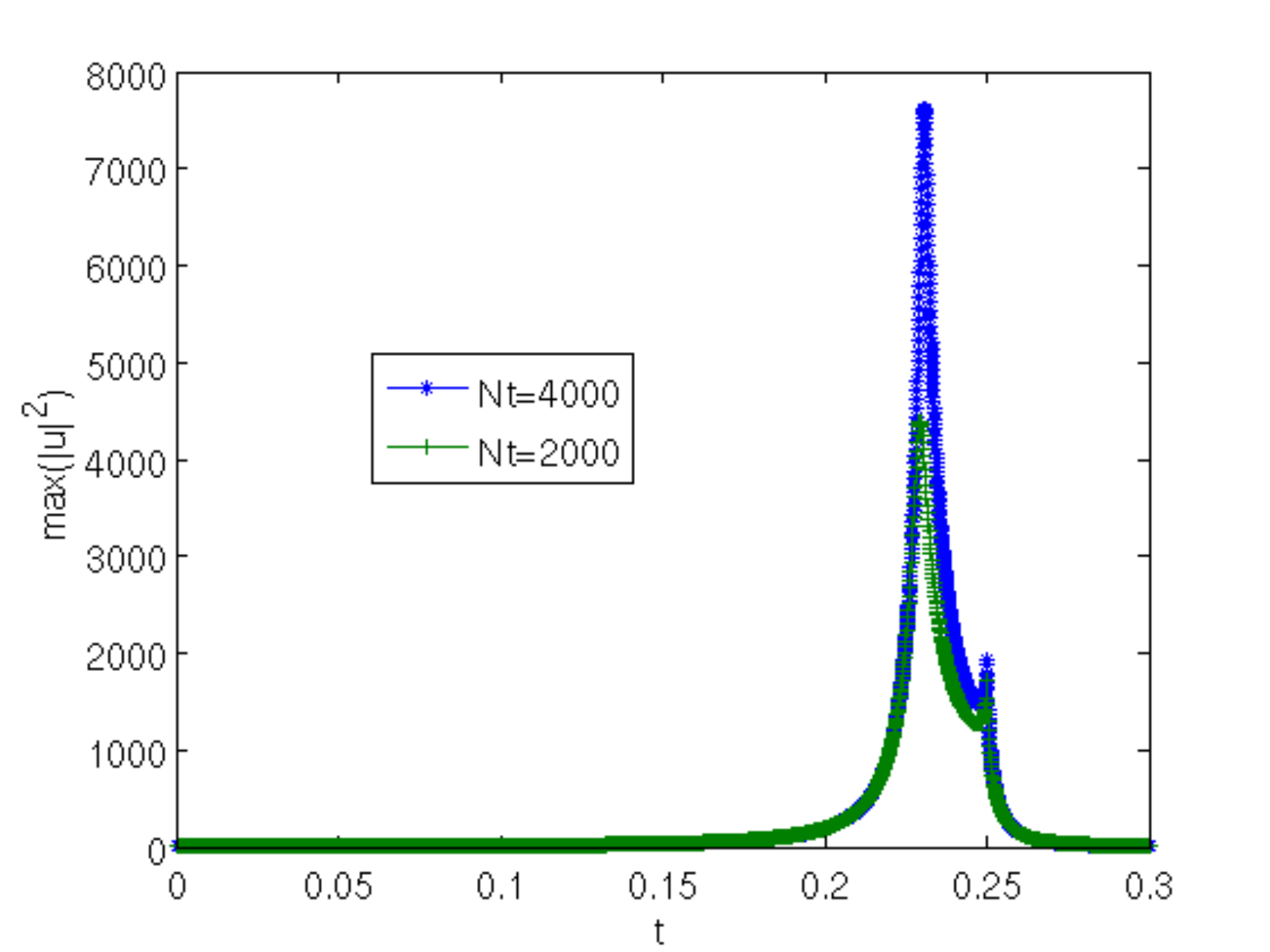}
\includegraphics[width=0.45\textwidth]{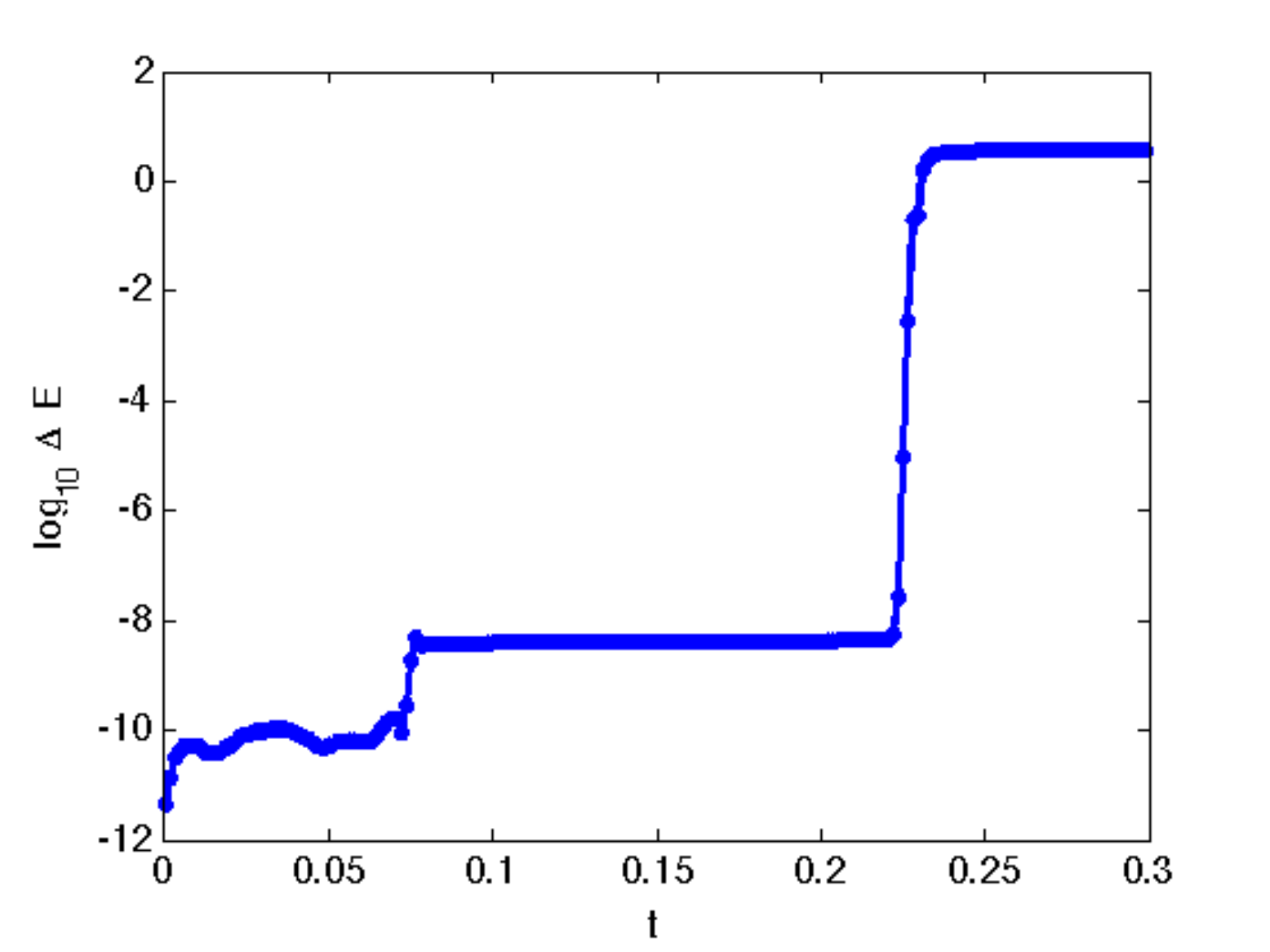}
\caption{Evolution of $max(|u|^{2})$ and the numerically computed energy
in dependence of time for the solution to the focusing DS II equation 
(\ref{DSII}) for an initial condition of the form (\ref{ozawagauss}) with $D=0.1$.}
\label{ampluozpg1}
\end{figure}

\begin{figure}[htb!]
\centering
\includegraphics[width=0.45\textwidth]{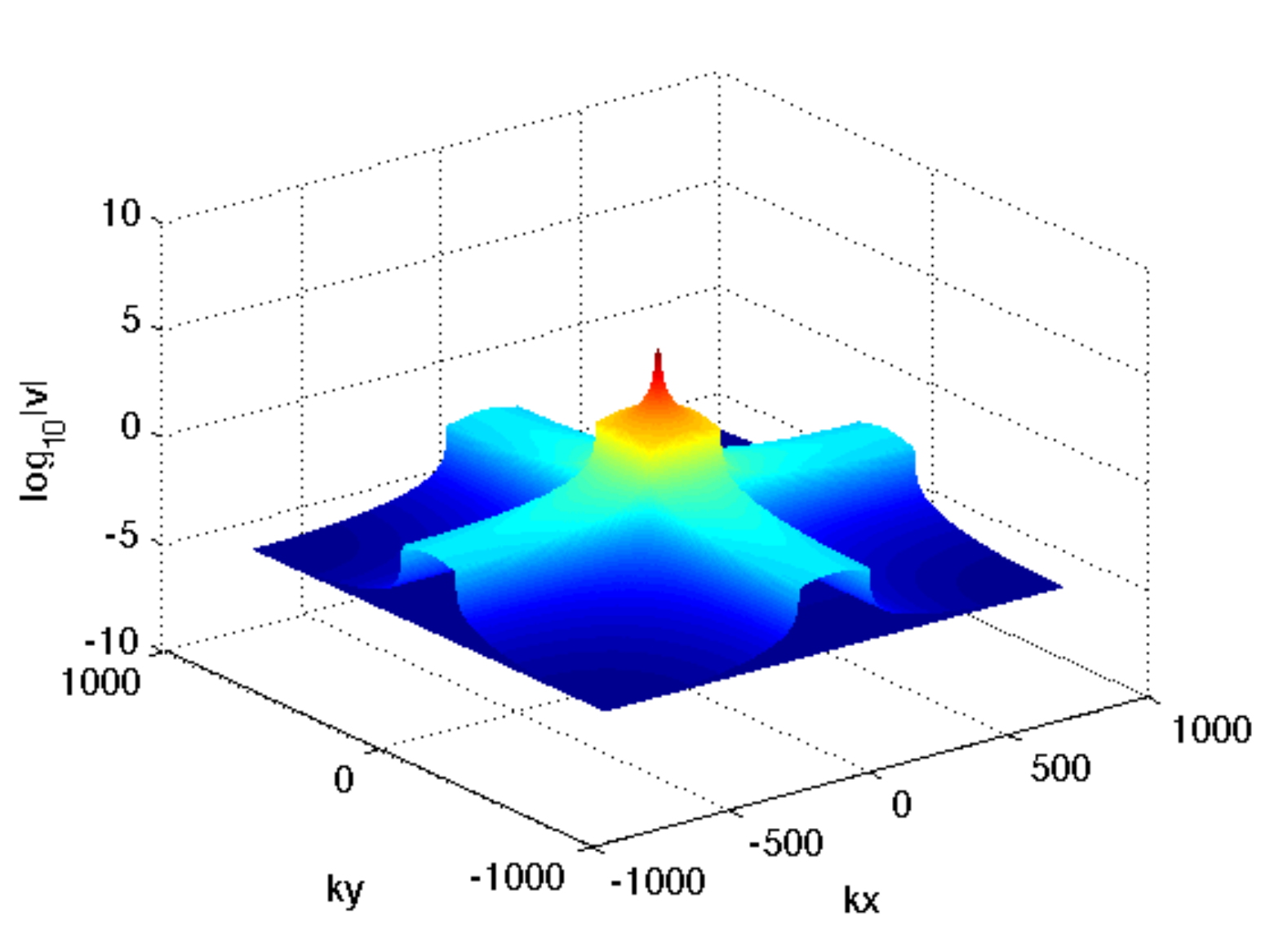}
\includegraphics[width=0.45\textwidth]{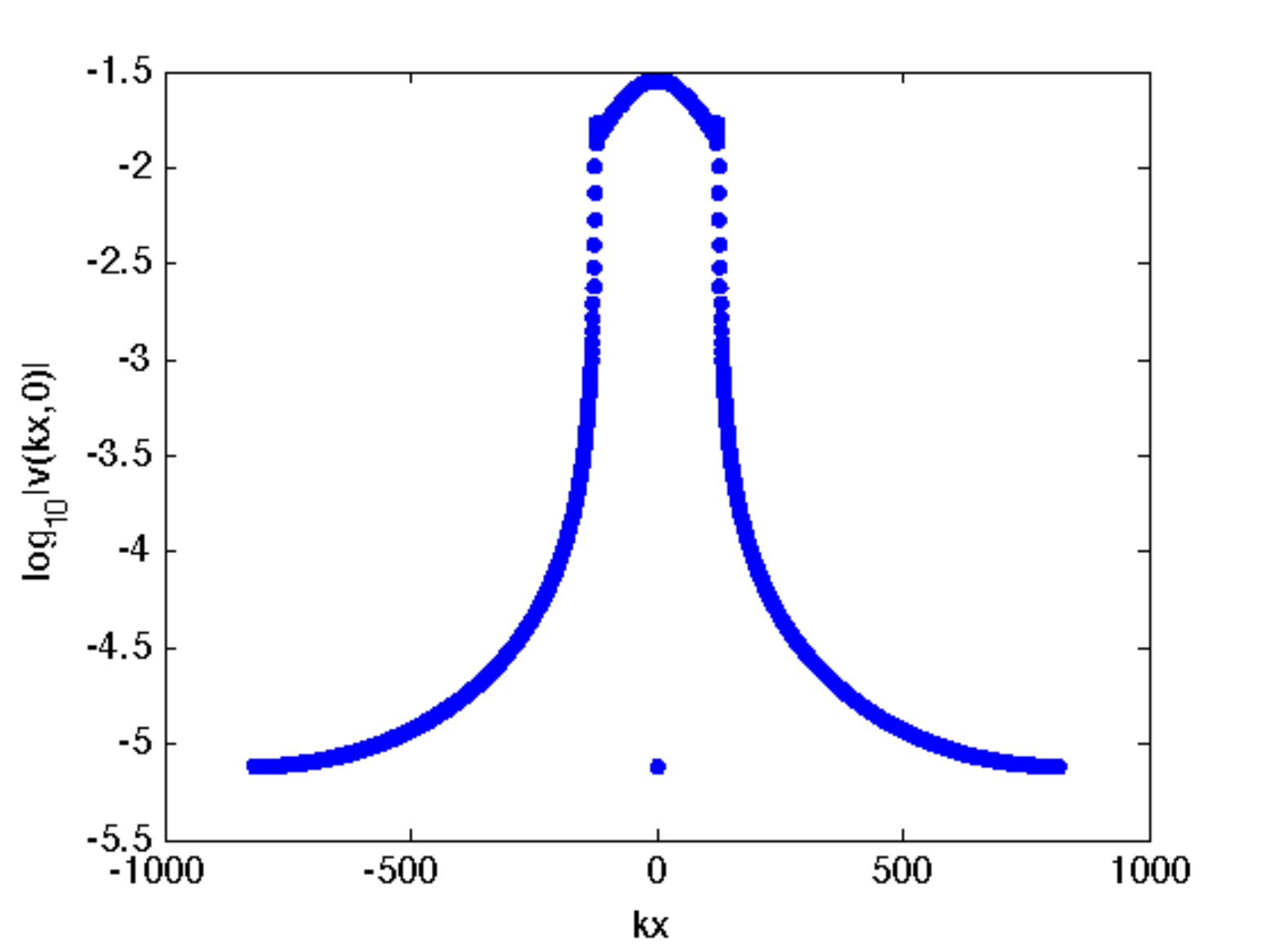}
\caption{Fourier coefficients of the solution to the focusing DS II equation (\ref{DSII}) at $t=0.15$
for an initial condition of the form (\ref{ozawagauss}) with $D=0.1$.} 
\label{oz01cf}
\end{figure}

\begin{figure}[htb!]
\centering
\includegraphics[width=0.45\textwidth]{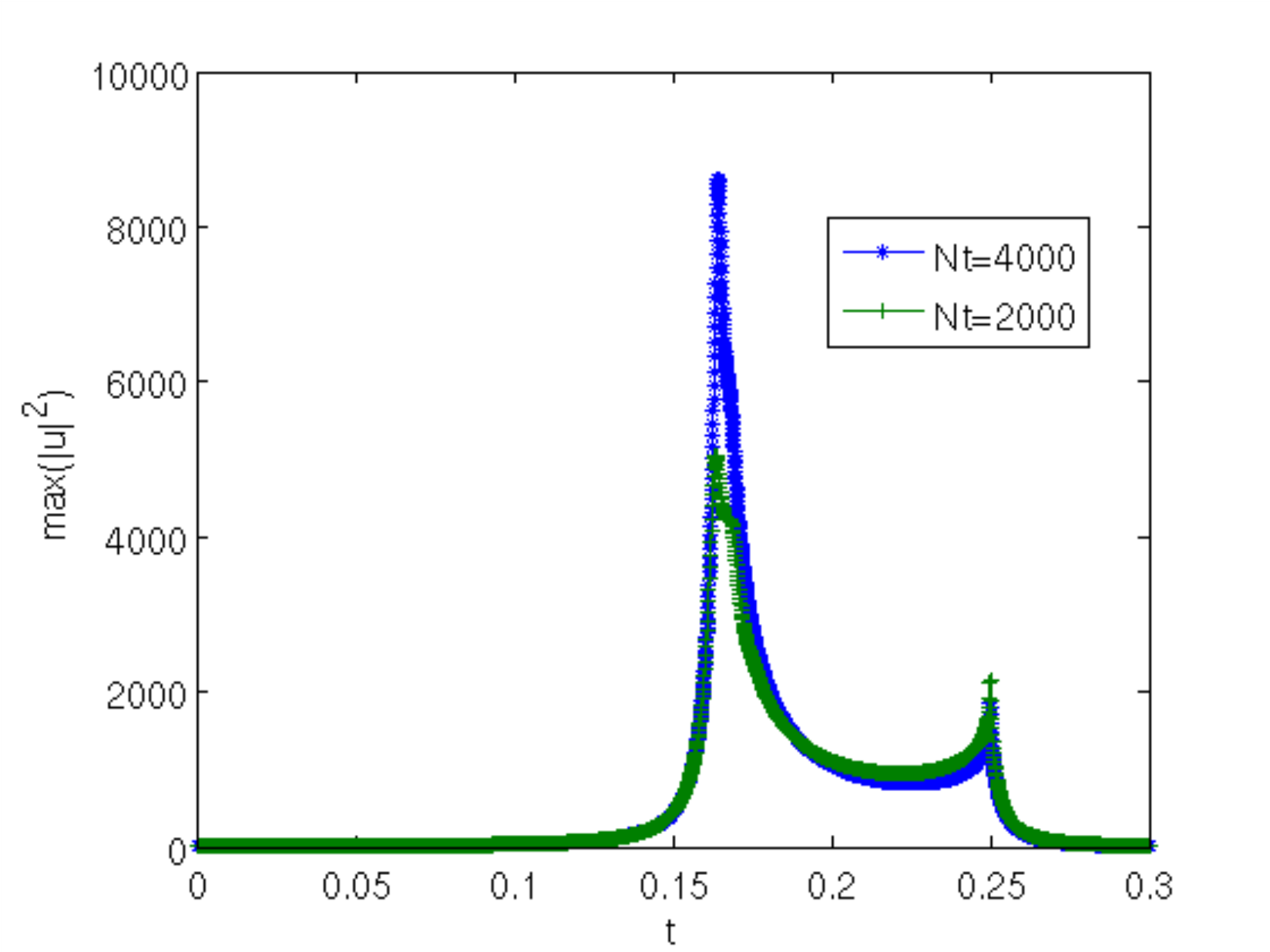}
\includegraphics[width=0.45\textwidth]{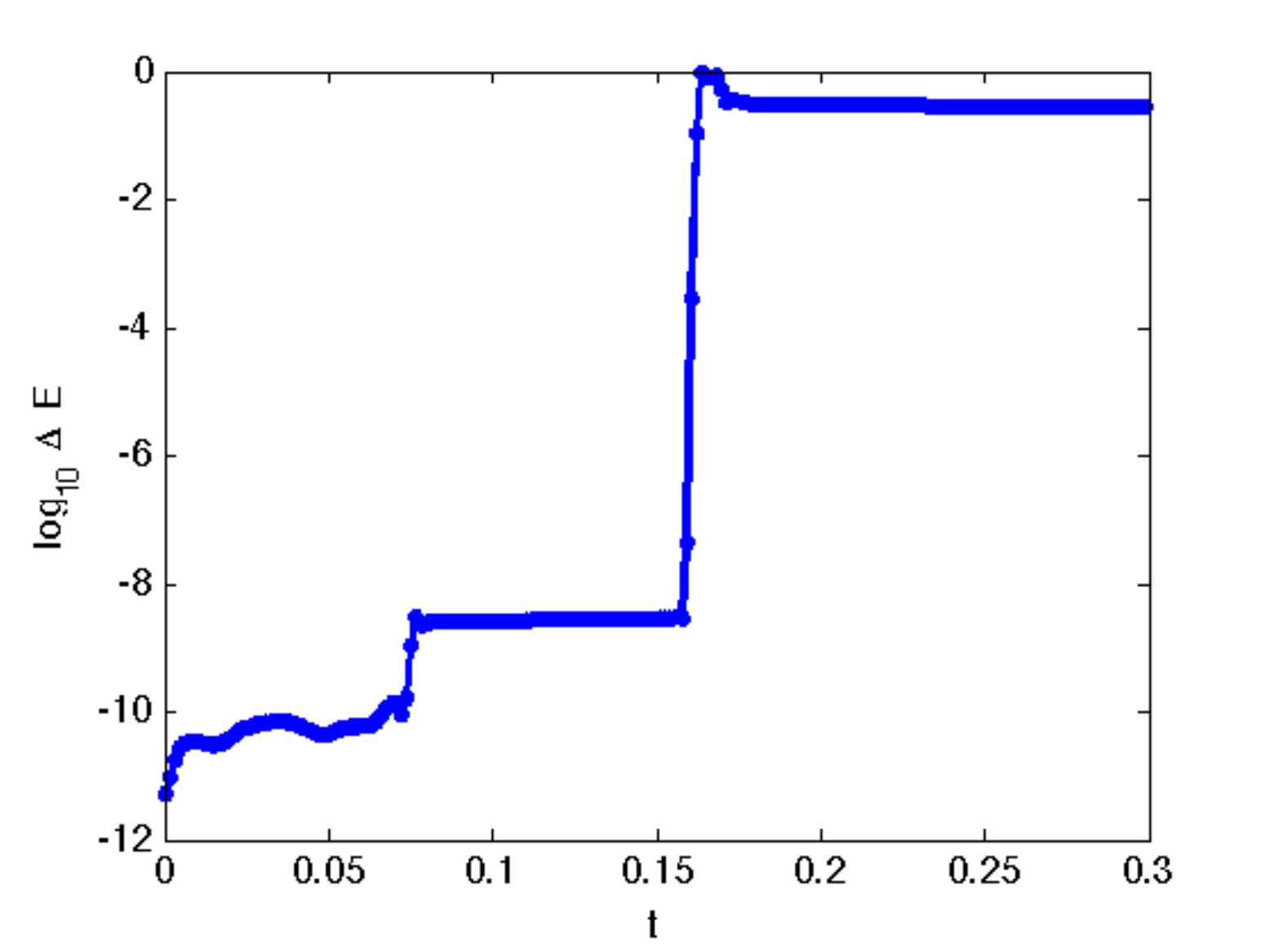}
\caption{Evolution of $max(|u|^{2})$ and the numerically computed 
energy for the solution to the focusing DS II equation (\ref{DSII})
 for an initial condition of the form (\ref{ozawagauss}) with $D=0.5$.}
\label{ampluozpg5}
\end{figure}

We consider an initial condition of the form 
\begin{equation}
    u(x,y,0) = 2\frac{ \exp \left(- i (x^2 - y^2)
    \right)}{1+x^2 +y^2} + D \exp(-(x^2+y^2))
    \label{ozawagauss}.
\end{equation}
For $D=0.1$ and
$N_t = 2000$, 
we show the behavior of $|u|^2$ at different times in Fig.~\ref{uozpg1}.
The time evolution of $\underset{x,y}{\max} |u(x,y,t)|^2$ is shown in 
Fig.~\ref{ampluozpg1}.
We observe a jump of the energy indicating blowup at the time $t_c\sim0.2332$.
The  Fourier coefficients at $t_c=0.15$ in 
Fig.~\ref{oz01cf} show that the wanted spatial resolution is achieved.
\\
\\
The same experiment with $D=0.5$ appears again to show blow up, but at an earlier time $t_c\sim0.1659$, see 
 Fig.~\ref{ampluozpg5}.
\\
\\
Thus the energy added by the perturbation of the form $D 
\exp(-(x^2+y^2))$ seems to lead to a blowup before the critical time of the 
Ozawa solution. This means that the blowup in the Ozawa solution is 
clearly a generic feature at least for initial data close to Ozawa for 
the focusing DS II equation.

 \subsection{Deformation of the Ozawa solution}
 
\begin{figure}[htb!]
\centering
\includegraphics[width=0.45\textwidth]{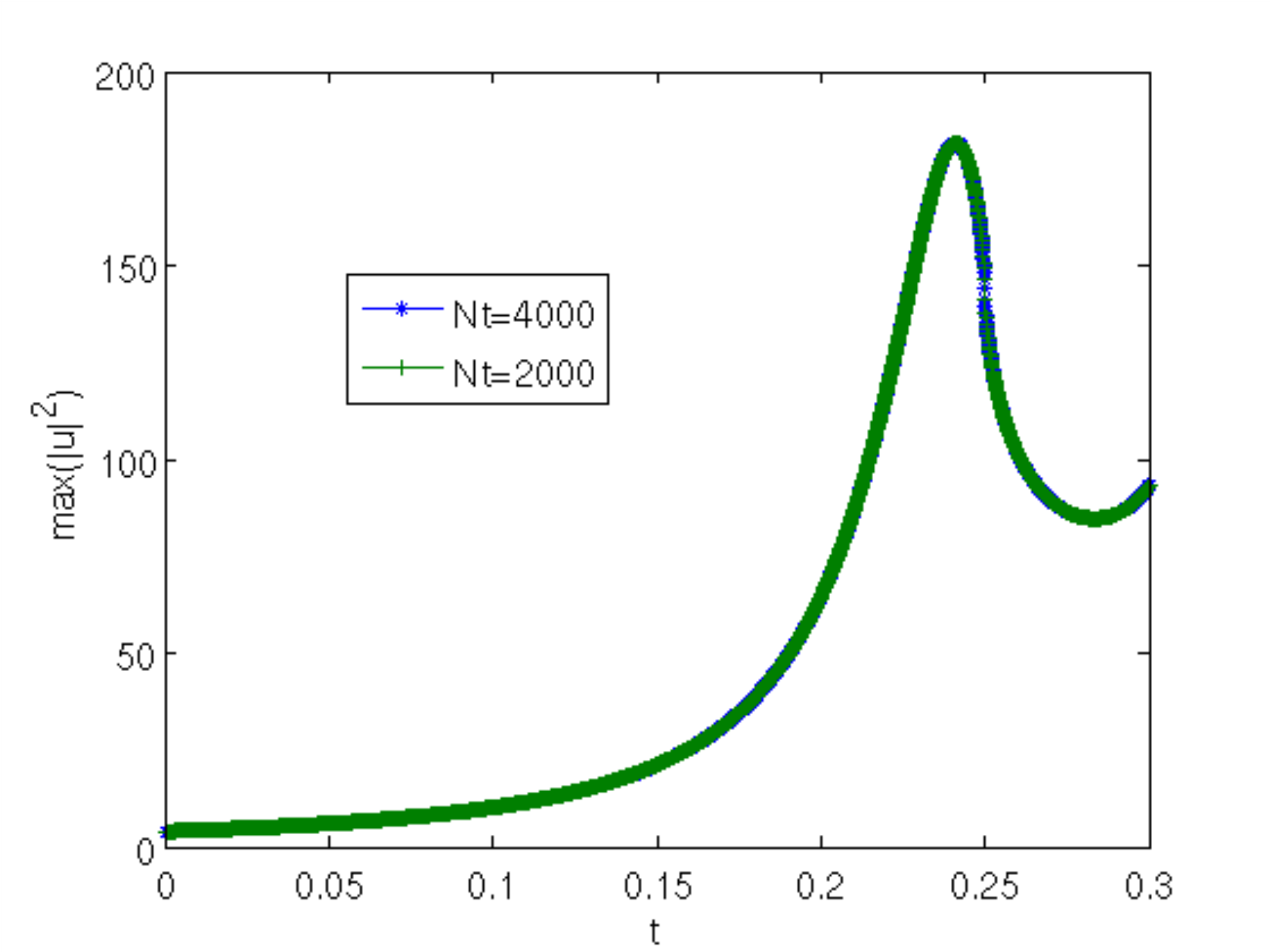}
\includegraphics[width=0.45\textwidth]{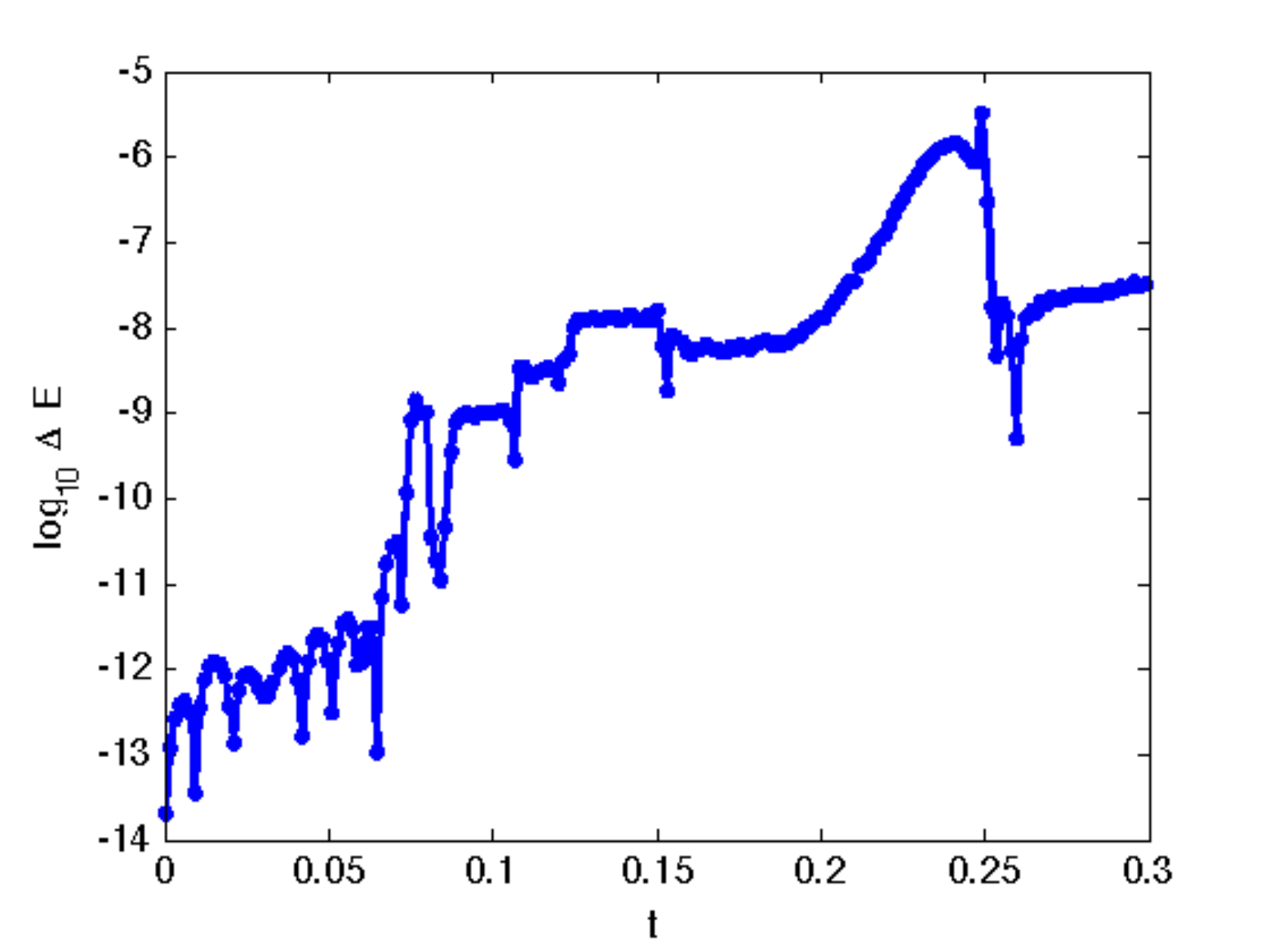}
\caption{Evolution of  $max(|u|^{2})$ and the numerically computed 
energy $E$ in dependence of time for  a solution to the
focusing DS II equation 
(\ref{DSII}) for an initial condition of the form (\ref{ozyd}) with $\nu=0.9$.}
\label{ampozy09}
\end{figure}

We study deformations of Ozawa initial data of the form
\begin{equation}\label{ozyd}
 u(x,y,0) = 2\frac{\exp \left( -i(x^{2}-(\nu 
 y)^{2})\right)}{1+x^{2}+(\nu y)^{2}},
\end{equation}
i.e.,  a deformation in the $y$-direction.
The computations are carried out with $N_{x}=N_{y}=2^{15}$ points for 
$x\times y \in [-20\pi, 20\pi] \times [-20\pi, 20\pi]$ and $t\in[0, 0.3]$.
\\
\\ 
For $\nu=0.9$, we observe a maximum of the solution at $t=0.2441$, see 
Fig. \ref{ampozy09}, followed by a second maximum, but there is no 
indication of a
blowup. Energy conservation is in principle high enough to 
indicate that the solution stays regular on the considered time scales.

%
\begin{figure}[htb!]
\centering
\includegraphics[width=0.45\textwidth]{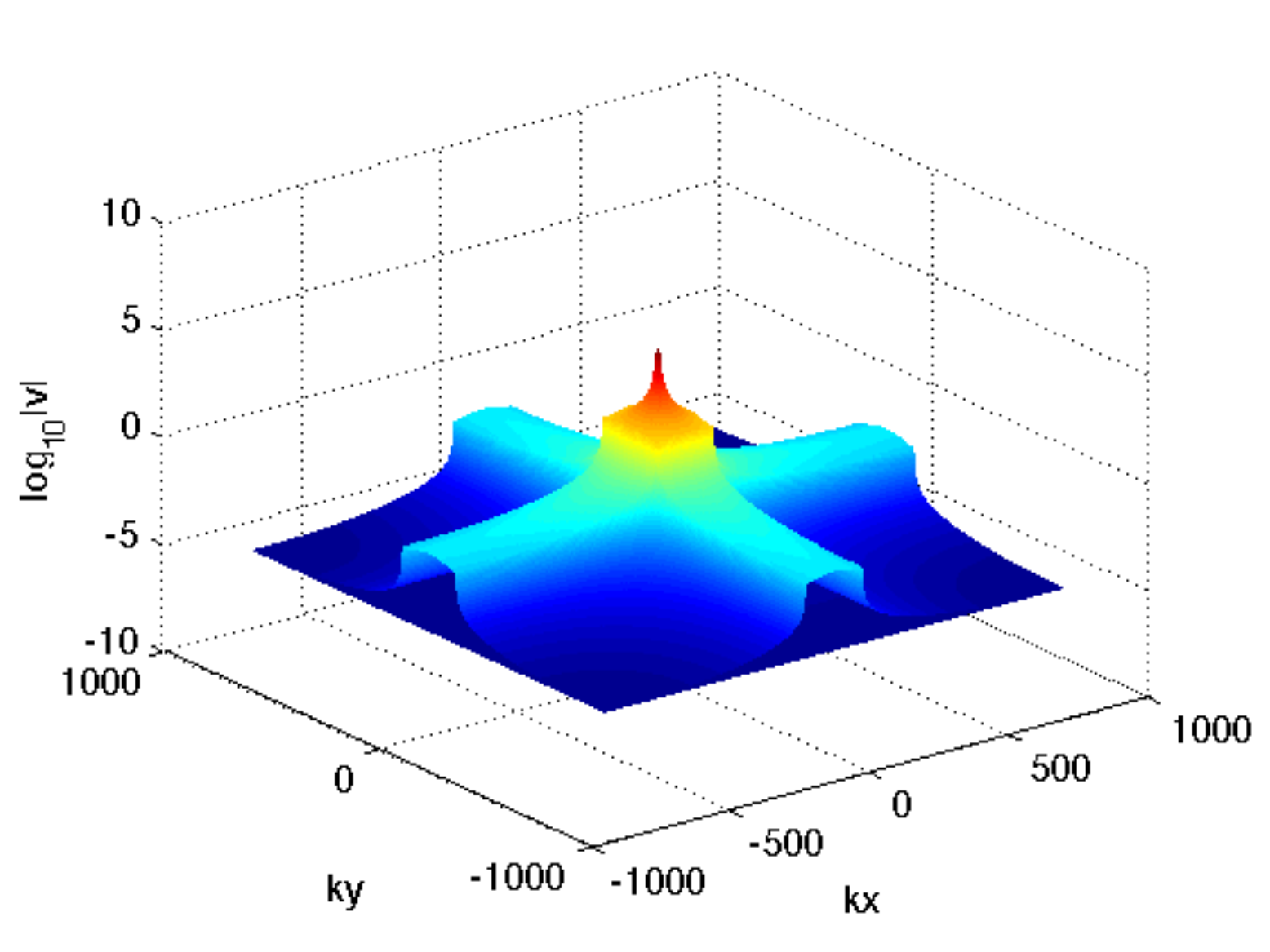}
\includegraphics[width=0.45\textwidth]{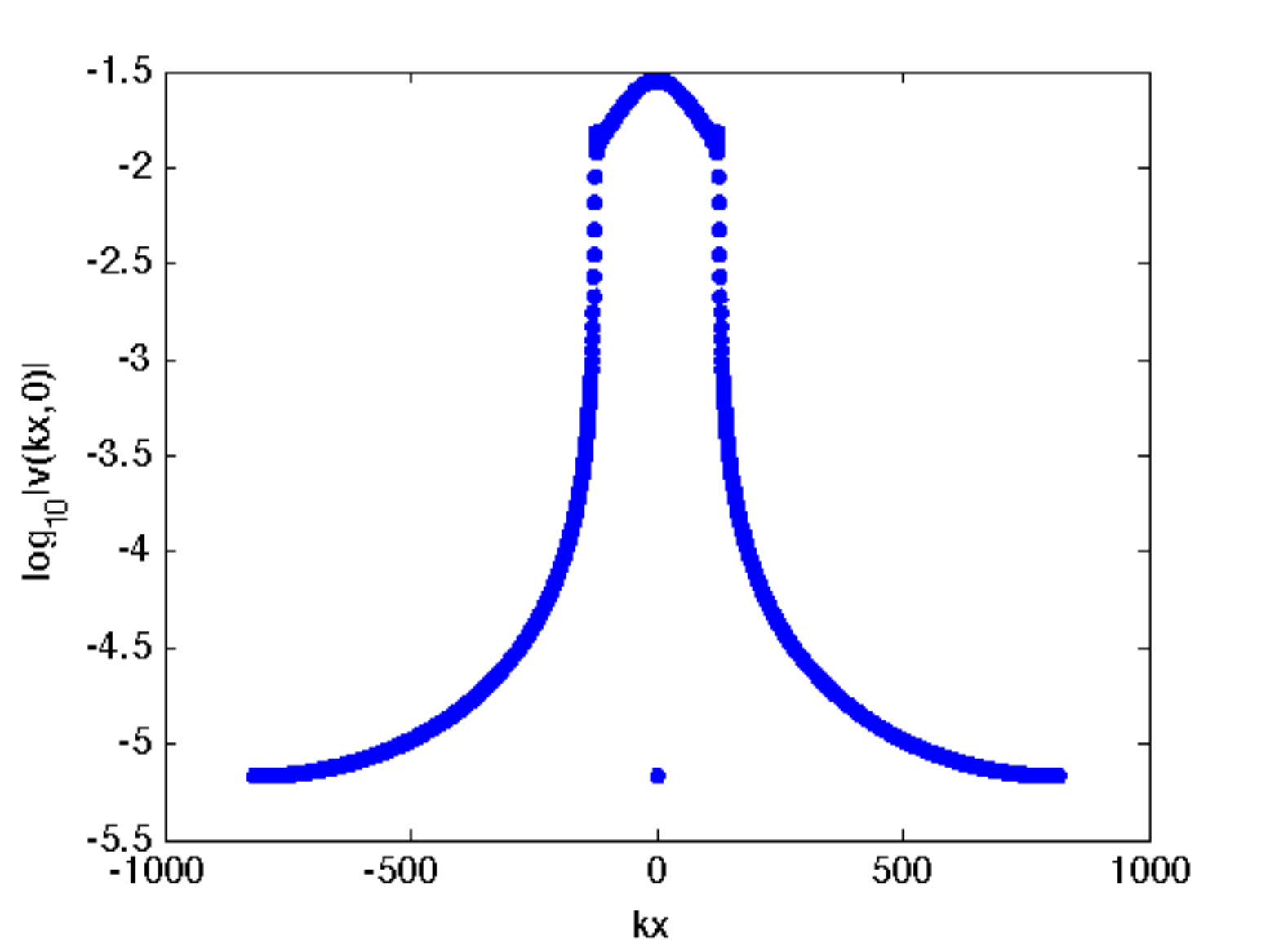}
\caption{Fourier coefficients of the solution to the focusing DS II equation 
(\ref{DSII}) 
for an initial condition of the form (\ref{ozyd}) with $\nu=0.9$  at $t=0$.}  
\label{ozy09cf}
\end{figure}

The  Fourier coefficients at $t=0.15$  in Fig.~\ref{ozy09cf} show  
the wanted spatial resolution.

\begin{figure}[htb!]
\centering
\includegraphics[width=0.45\textwidth]{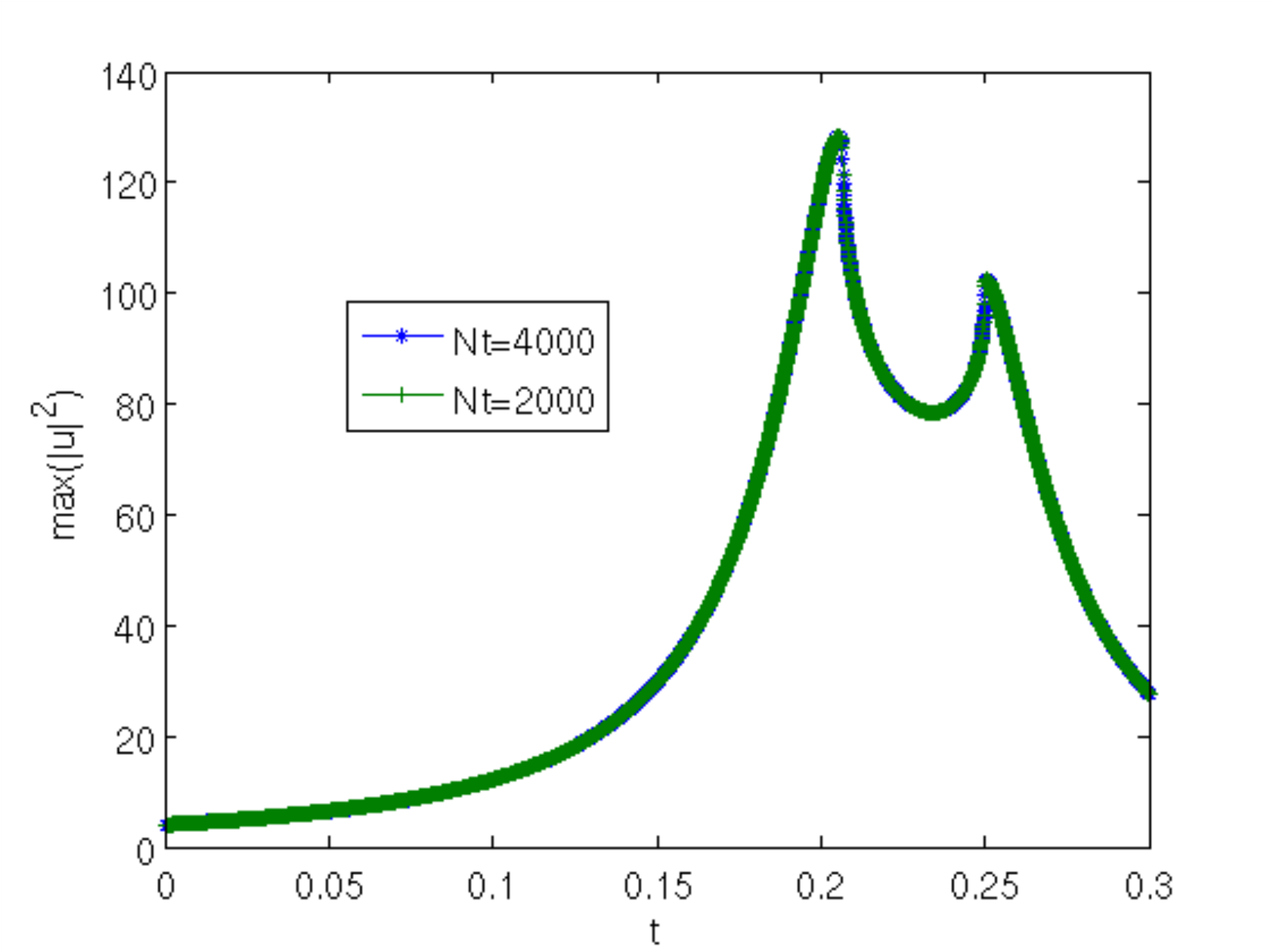}
\includegraphics[width=0.45\textwidth]{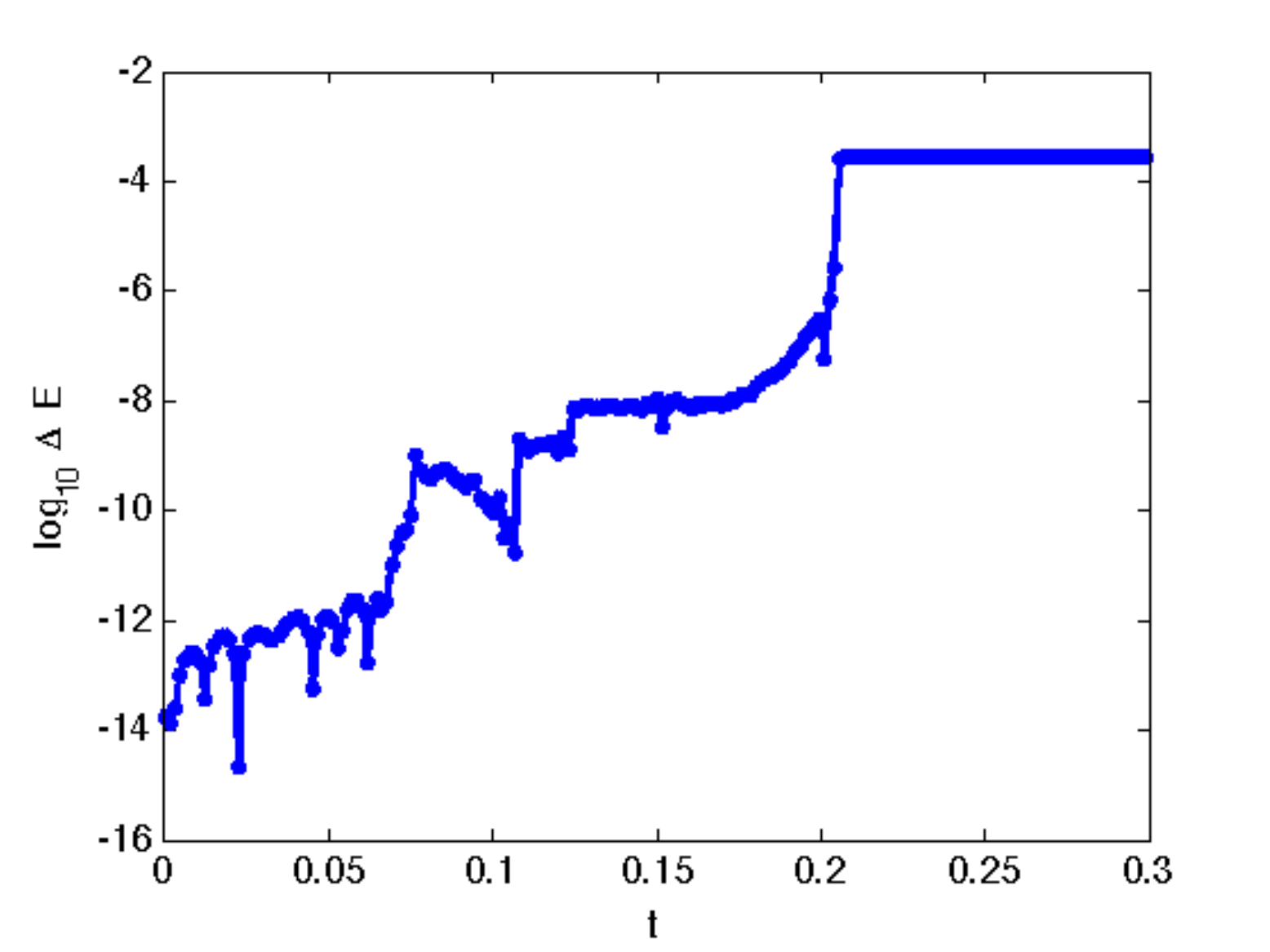}
\caption{Evolution of  $max(|u|^{2})$ and the numerically computed 
energy $E$ for a solution to the focusing DS II equation 
(\ref{DSII}) 
for an initial condition of the form (\ref{ozyd}) with $\nu=1.1$.}
\label{ampozy11}
\end{figure}

The situation is similar for $\nu=1.1$. The maximum of the solution is observed at $t=0.2254$, see 
Fig. \ref{ampozy11}, followed again by a second maximum. Energy 
conservation appears once more to rule out a blowup in this case.

\begin{figure}[htb!]
\centering
\includegraphics[width=0.45\textwidth]{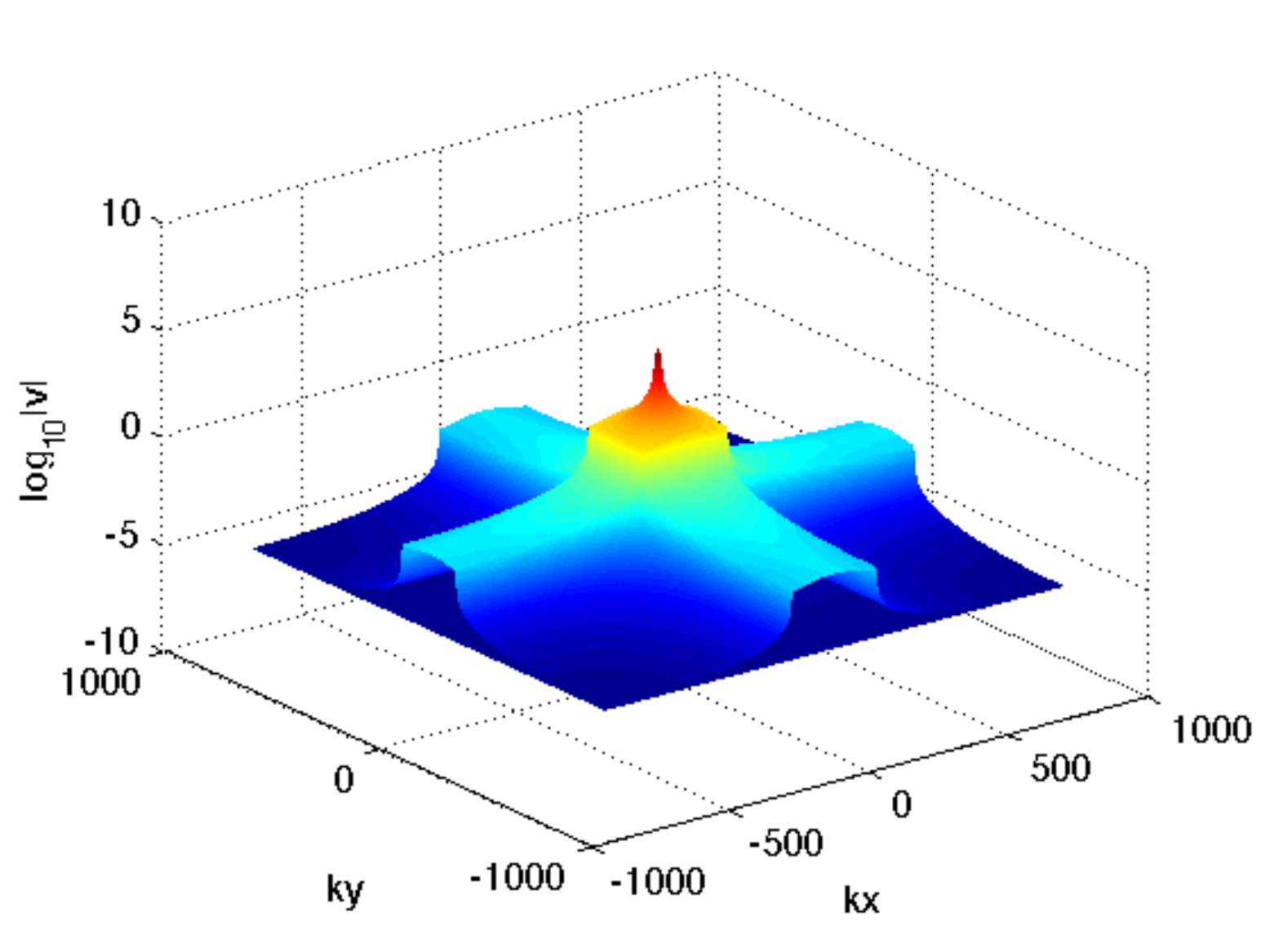}
 \includegraphics[width=0.45\textwidth]{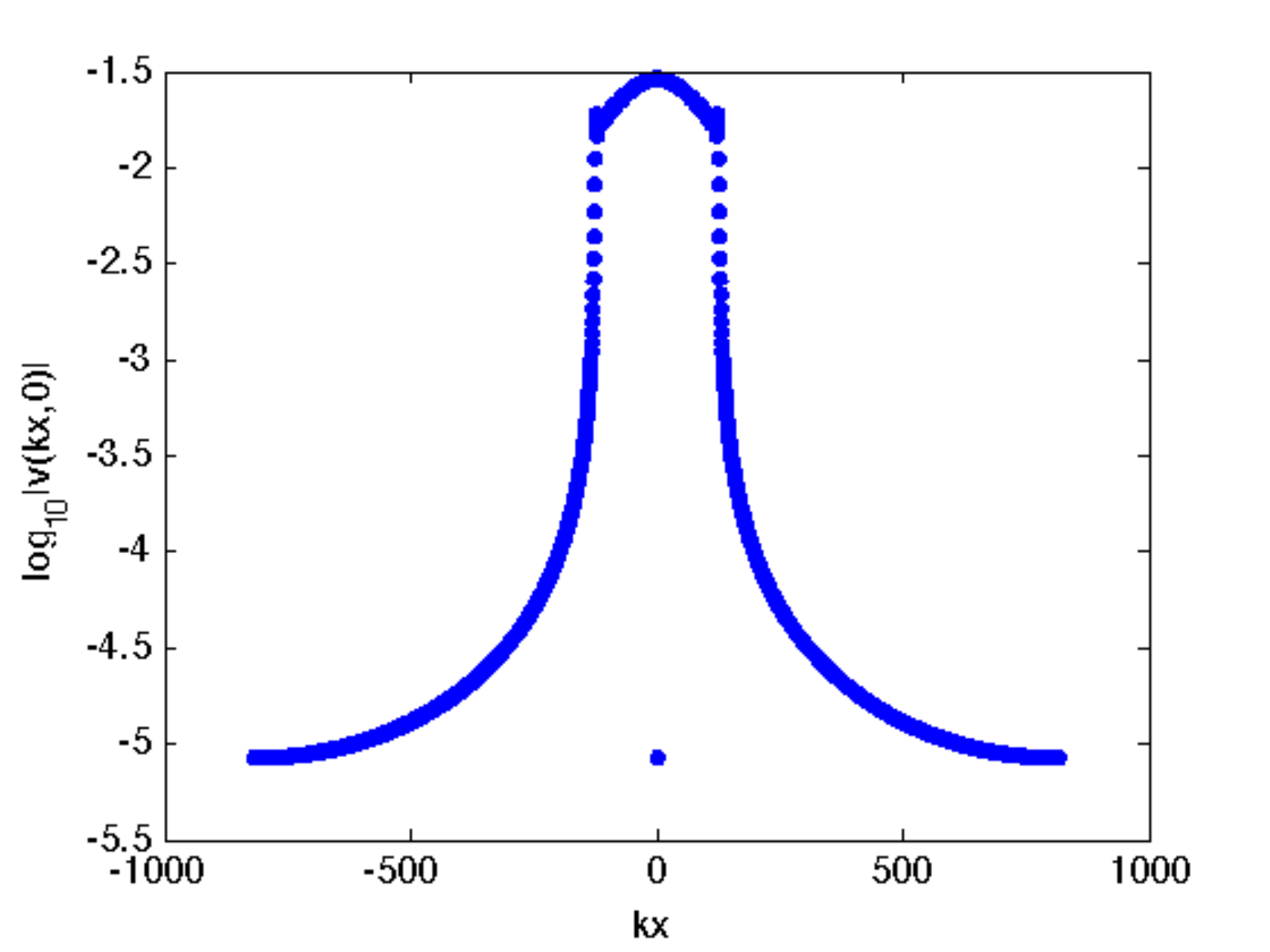}
\caption{Fourier coefficients of the solution to the focusing DS II equation 
(\ref{DSII}) 
for an initial condition of the form (\ref{ozyd}) with $\nu=1.1$  at $t=0.15$.}  
\label{ozy11cf}
\end{figure}

The Fourier coefficients at $t=0.15$ in Fig.~\ref{ozy11cf} again show the wanted spatial resolution.\\

\section{Conclusion}
In this paper we have numerically studied long time behavior and 
stability of exact solutions to the focusing DS II equation with an 
algebraic falloff towards infinity. We have shown that the necessary 
resolution can be achieved with a  parallelized version of a 
spectral code. The spatial resolution as seen at the Fourier 
coefficients was always well beyond typical plotting accuracies of 
the order of $10^{-3}$. For the time integration we used an 
unconditionally stable fourth order splitting scheme. As argued in 
\cite{ckkdvnls,KR}, the numerically computed energy of the solution 
gives a valid indicator of the accuracy for sufficient spatial 
resolution.  To ensure the latter, we always presented the Fourier 
coefficients of the solution at a time before a singularity appeared. 
In addition we show here that the numerically computed energy
indicates blowup by jumping to a different value in cases where the code runs beyond a singularity in 
time.

After testing the code for exact solutions, the lump and the blowup 
solution by Ozawa, we showed that both solutions are critical in the 
following sense: adding energy to it leads to a blowup for the lump, 
and an earlier blowup time for the Ozawa solution. For initial data 
with less energy, no blowup was observed in both cases, the initial 
data asymptotically just seem to be dispersed. This is in accordance 
with the conjecture in \cite{MFP} that solutions to the focusing DS II 
equations either blow up or disperse. In particular the lump is 
unstable against both blowup and dispersion, in contrast to the lump 
of the KP I equation that appears to be stable, see for instance 
\cite{PS}. Note that the 
perturbations we considered here test the nonlinear regime of the PDE 
for which so far no analytical results appear to be established.

\bibliographystyle{siam}
\bibliography{biblio}{}

\end{document}